\documentclass[a4paper,11pt]{amsart}

\usepackage{amsmath,amssymb}
\usepackage{amssymb}
\usepackage{latexsym}
\usepackage{amsmath}

\begin{document}

\newtheorem{df}{Definition}[section]
\newtheorem{ex}[df]{Example}
\newtheorem{thm}[df]{Theorem}
\newtheorem{lm}[df]{Lemma}
\newtheorem{pr}[df]{Proposition}
\newtheorem{cor}[df]{Corollary}
\newtheorem{rmk}[df]{Remark}

\newcommand{\FLm}{p_{i_j}^-}
\newcommand{\FLp}{p_{i_j}^+}

\newcommand{\FLmq}{\lfloor\frac{q_{i_j}-p_1}{2}\rfloor}
\newcommand{\FLpq}{\lfloor\frac{p_1+q_{i_j}-1}{2}\rfloor}

\newcommand{\cal}{\mathcal}

\title[The spherical growth series of amalgamated free products $G(p_1,p_2,\dots,p_n)$]
{The spherical growth series of amalgamated free products of infinite cyclic groups}

%
\author{Michihiko Fujii}
\address{Department of Mathematical Sciences, 
Faculty of Science, University of the Ryukyus, 
Nishihara-cho, Nakagami-gun, Okinawa 903-0213, Japan.}
\email{mfujii@sci.u-ryukyu.ac.jp}
\author{Takuya Sakasai}
\address{Graduate School of Mathematical Sciences, 
University of Tokyo, 3-8-1 Komaba, Meguro-ku,
Tokyo 153-8914, Japan.}
\email{sakasai@ms.u-tokyo.ac.jp}

\date{\today}

\thanks{
The authors were partially supported by KAKENHI 
(No.18K03283, No.24K06740), 
Japan Society for the Promotion of Science, Japan.
}


\begin{abstract}

Let $n$ be an integer greater than 1.
We consider a group presented as
$G(p_1,p_2,\dots,p_n)=$
\linebreak
$ \langle\ x_1,x_2,\dots, x_n\ |\ 
x_1^{p_1} =x_2^{p_2}=\cdots =x_n^{p_n} \
\rangle$,
with
integers $p_1,p_2,\dots,p_n$ satisfying $2 \leq p_1 \leq p_2 \leq \cdots \leq p_n$.
This group is an amalgamated free product of infinite cyclic groups
and is geometrically realized as the fundamental group 
of a Seifert fiber space over the 2-dimensional disk with $n$ cone points
whose associated cone angles are $\frac{2\pi}{p_1},\frac{2\pi}{p_2},\dots,\frac{2\pi}{p_n}$.
In this paper, 
we present a formula for the spherical growth series of the group $G(p_1,\dots,p_n)$
with respect to the generating set $\{x_1,\dots,x_n,x_1^{-1},\dots,x_n^{-1}\}$.
We show that from this formula,
a rational function expression
for the spherical growth series of $G(p_1,\dots,p_n)$ can be derived in concrete form
for given $p_1,\dots,p_n$.
In fact, we wrote an elementary computer program based on this formula
that yields an explicit form of a single rational fraction expression for the spherical growth series of $G(p_1,\dots,p_n)$.
We present such expressions for several tuples $(p_1,\dots,p_n)$.
In 1999, C. P. Gill obtained a similar formula for the same group
in the case $n=2$
and showed that there exists a rational function expression
for the spherical growth series of $G(p_1,\dots,p_n)$ for $n \geq 2$.
\end{abstract}

\maketitle

%
%
%


\section{Introduction}\label{introduction}

\noindent
Let $n$ be an integer greater than 1 and
$p_1,p_2, \dots,p_n$ be integers satisfying $2 \leq p_1 \leq p_2 \leq \cdots \leq p_n$.
We consider a group $G(p_1,p_2,\dots,p_n)$ presented as
\begin{eqnarray}\label{group-presentation}
G(p_1,p_2,\dots,p_n)\ =\ \langle x_1,x_2,\dots, x_n\ |\ 
x_1^{p_1} =x_2^{p_2}=\cdots =x_n^{p_n}\ \rangle
\end{eqnarray}
and its natural symmetric generating set 
$\Sigma:=\{x_1,\dots,x_n,x_1^{-1},\dots,x_n^{-1} \}$.
The group $G(p_1,\dots,p_n)$ is an amalgamated free product
of $n$ infinite cyclic groups
and 
is geometrically realized as the fundamental group 
of a Seifert fiber space over the 2-dimensional disk with $n$ cone points
whose associated cone angles are $\frac{2\pi}{p_1},\frac{2\pi}{p_2},\dots,\frac{2\pi}{p_n}$.
In this paper,
we investigate the combinatorial group structure of $G(p_1,\dots,p_n)$
with respect to the presentation (\ref{group-presentation}),
and in particular,
the spherical growth series ${\cal S}_{G(p_1,\dots,p_n)}(t)$ of $G(p_1,\dots,p_n)$ with respect to the generating set $\Sigma$.
(See (\ref{def-growth-series}) in Section 2 for the definition of ${\cal S}_{G(p_1,\dots,p_n)}(t)$.)
The aim of this paper is to obtain a formula for ${\cal S}_{G(p_1,\dots,p_n)}(t)$
from which a rational function expression 
can be derived.

For the particular case $n=2$ with $p_1=p_2$,
Edjvet and Johnson \cite{Ed-J} derived
a formula that yields a rational function expression
for ${\cal S}_{G(p_1,p_1)}(t)$.
Also, in the case $n=2$ with $(p_1,p_2)=(2,3)$,
Shapiro \cite{Shapiro} and 
Johnson, Kim and Song \cite{J-K-S}
independently
obtained the same rational function expression 
for ${\cal S}_{G(2,3)}(t)$.
Generalizing these results, 
Gill \cite{Gill} derived a formula yielding a rational function expression
for ${\cal S}_{G(p_1,p_2)}(t)$ for every pair $(p_1,p_2)$,
and moreover,  he showed that 
${\cal S}_{G(p_1\dots,p_n)}(t)$ possesses a rational function expression
for every $n$-tuple $(p_1,\dots,p_n)$,
where $2 \leq p_1 \leq p_2 \leq \cdots \leq p_n$.
Later, Fujii \cite{Fujii2} obtained another formula yielding a rational function expression
for ${\cal S}_{G(p_1,p_2)}(t)$ for every pair $(p_1,p_2)$.
In this paper, we generalize this all at once to arbitrary $n\ (\geq 2)$
and derive a formula that provides a rational function expression
for ${\cal S}_{G(p_1\dots,p_n)}(t)$ for every $n$-tuple $(p_1,\dots,p_n)$,
where $2 \leq p_1 \leq p_2 \leq \cdots \leq p_n$.
Using the formula given here,
we obtain an explicit rational function expression for ${\cal S}_{G(p_1,\dots,p_n)}(t)$.

The principal method employed in this paper is 
the so-called ``suitable-spread procedure."
It has been previously utilized for applications to the braid group of three strands 
(see \cite{Berger}), 
the Artin groups of dihedral type (see \cite{M-M}),
the pure Artin groups of dihedral type (see \cite{Fujii1}),
and the group $G(p_1,p_2)$ (see \cite{Fujii2}).
Furthermore, it can be applied to the group $G(p_1,\dots,p_n)$.
In this regard,
the fundamental strategy of the paper aligns with \cite{Fujii2}.
However, it is important to highlight that
when applying the procedure to the group $G(p_1,\dots,p_n)$,
a more detailed discussion is needed.
Through the discussion, we obtain the formula 
for the spherical growth series of $G(p_1,\dots,p_n)$
that includes all known cases of the group so far.

Before we apply the suitable-spread procedure to the group $G(p_1,\dots,p_n)$,
for each element of $G(p_1,\dots,p_n)$,
we must identify one unique, particular representative $\nu$ 
from among all of the representatives,
which is called the {\it modified normal form} of the element
 (see Proposition \ref{pr-mnf}).
Then, we show that
if an element $g$ of $G(p_1,\dots,p_n)$ is of a particular type (Type 3 defined in Section \ref{section-ns-cond}),
applying the suitable-spread procedure to its modified normal form $\nu$,
we obtain all of the geodesic representatives of $g$
(see Propositions \ref{pr-ss-2} and \ref{pr-sc-T3}).
For some elements of $G(p_1,\dots,p_n)$,
there exists just one geodesic representative,
but there are also elements
for which there exist multiple geodesic representatives.
However, even in the latter case,
it is shown that 
only one of these geodesic representatives
can be uniquely determined
(see Section \ref{section-ugr}).
This determination procedure for the case $n \geq 3$
becomes complicated in comparison with the case $n=2$.
In fact, when $n \geq 3$,
we have to pay attention to whether $p_k-p_1$ is even or odd for every 
$k \in \{1,\dots,n\}$.
For this reason,
before beginning the procedure,
for the sake of clearity
we re-express each integer $p_k$ according to whether $p_k-p_1$ is even or odd,
using $q_1,\dots,q_m$ and $r_{m+1},\dots,r_n$ in these two cases, respectively,
according to the assignment descrived in (\ref{pqr}).
In conjunction with this,
we re-express the elements $x_1,\dots, x_n$
using $y_1,\dots,y_m$ and $z_{m+1},\dots,z_n$,
according to the assignment descrived in (\ref{xyz}).
The representation of $G(p_1,\dots,p_n)$ in terms of 
$y_1,\dots,y_m$ and $z_{m+1},\dots,z_n$ is given in (\ref{new-presentation}).
With this notation, we describe a method that uniquely specifies 
a single output of the suitable-spread procedure.
This output is regarded as a {\it proper} geodesic representative
for the element in question,
because the sets consisting of the geodesic representatives
of the elements of $G(p_1,\dots,p_n)$
have simple expressions provided by the decomposition of $G(p_1,\dots,p_n)$.
In every case, for each element of $G(p_1,\dots,p_n)$,
a proper geodesic representative is uniquely determined.
Then, in Section \ref{section-growth},
using a lemma concerning rational function expressions for 
formal power series associated with sets consisting of words of certain kinds
 (Lemma \ref{lm-growth}),
we obtain a formula that provides a rational function expression 
for the spherical growth series
${\cal S}_{G(p_1,\dots,p_n)}(t)$
(see Theorem \ref{th-growth}).

As an important step in obtaining our results,
we define a monoid $G(p_1,\dots,p_n)^+$ associated with the group $G(p_1,\dots,p_n)$,
which is called the positive monoid.
Dehornoy and Paris \cite{D-P} showed that $G(p_1,\dots,p_n)^+$ is a Garside monoid,
which implies the unique existence of 
both the Garside normal form and the modified normal form 
for each element of $G(p_1,\dots,p_n)$
(see Propositions \ref{pr-nf} and \ref{pr-mnf}).
The uniqueness of these normal forms plays a crucial role in the proofs of
many of the fundamental results in this paper.
The most important of these are presented in 
Propositions \ref{pr-T1,2,3^+,3^-}, \ref{pr-ss-1} and \ref{pr-sc-T3}.

While we are interested in the generating set $\Sigma$ in this paper,
motivated by works of Shapiro \cite{Shapiro} and Mann \cite{Mann1,Mann2},
Nakagawa, Tamura and Yamashita \cite{N-T-Y} 
studied a distinct generating set $\Sigma'$
for the same group $G(p_1,\dots,p_n)$ and
obtained a rational function expression for
the spherical growth series of $G(p_1,\dots,p_n)$ with respect to $\Sigma'$.
The set $\Sigma'$ is the union of $\Sigma$ and $\{\Delta,\Delta^{-1}\}$,
where $\Delta =x_1^{p_1}$ and 
$\Delta^{-1} \ =x_1^{-p_1}$.
It is worth noting that the rational function expression for $G(p_1,\dots,p_n)$
with respect to $\Sigma'$ has a simpler form than that with respect to $\Sigma$.

As supplemental material,
we provide a computer program 
written in {\it Mathematica}, released by Wolfram Research,
which yields explicit forms of single rational fraction expressions
for the spherical growth series ${\cal S}_{G(p_1,\dots,p_n)}(t)$.
This program is available at the following website:
\begin{center}
http://hdl.handle.net/20.500.12000/0002019730
\end{center}

Before ending this section,
we establish some notation.
Throughout the paper,
we use the following notation to represent the floor function:
\begin{eqnarray}\label{floor}
\lfloor r \rfloor:={\rm max} \{ n \ ; \ n \leq r ,\  n \in {\bf Z} \} \quad
\mbox{for $r \in {\bf R}$}.
\end{eqnarray}
Also, for each $k \in \{1,\dots,n\}$,
we define the following:
\begin{eqnarray}\label{p-+}
p_k^-:=\lfloor\frac{p_k-p_1}2\rfloor,\quad
p_k^+:=\lfloor\frac{p_1+p_k-1}2\rfloor.
\end{eqnarray}
Then,
if $p_k-p_1$ is even,
we have
\begin{eqnarray}\label{pq-oo}
p_k^-=\frac{p_k-p_1}{2} \ \ 
\mbox{and} \ \ 
p_k^+=\frac{p_1+p_k-2}{2},
\end{eqnarray}
while if $p_k-p_1$ is odd,
we have
\begin{eqnarray}\label{pq-oe}
p_k^-=\frac{p_k-p_1-1}{2} \ \ 
\mbox{and} \ \ 
p_k^+=\frac{p_1+p_k-1}{2}.
\end{eqnarray}


\section{Garside normal form and modified normal form}\label{section-mnf}

Let $n$ be an integer greater than 1 and
$p_1,p_2, \dots,p_n$ be integers satisfying $2 \leq p_1 \leq p_2 \leq \cdots \leq p_n$.
Then,
let $G(p_1,\dots,p_n)$ be the group presented as in (\ref{group-presentation}).
We define the group $\overline{G}(p_1,\dots,p_n)$ by
$\overline{G}(p_1,\dots,p_n):=
\langle \ x_1,x_2,\dots, x_n\ |\ x_1^{p_1} =x_2^{p_2}=\cdots =x_n^{p_n}=1\  \rangle$.
The sequence
\begin{eqnarray}\label{exact-sq}
1 \rightarrow {\bf Z} \rightarrow G(p_1,\dots,p_n) \rightarrow \overline{G}(p_1,\dots,p_n) \rightarrow 1
\end{eqnarray}
is exact,
because the center of $G(p_1,\dots,p_n)$ is the infinite cyclic group
generated by the element $x_1^{p_1}\ (=x_2^{p_2}=\dots =x_n^{p_n})$.
The group $\overline{G}(p_1,\dots,p_n)$ is the so-called orbifold fundamental group
of the 2-dimensional disk with $n$ cone points,
 whose associated cone angles are
$2 \pi/p_1,\dots,2 \pi/p_n$.
It is thus seen that
$G(p_1,\dots,p_n)$ is geometrically realized as the fundamental group of a Seifert fiber space
over the 2-dimensional orbifold.

We define the sets
\[
  \Sigma^+:=\{x_1,\dots,x_n\},\ \ 
  \Sigma^-:=\{x_1^{-1},\dots, x_n^{-1}\}\ \ \mbox{and}\ \ 
  \Sigma:=\Sigma^+ \cup \Sigma^-.
\]
Then, let
$\Sigma^*$, $(\Sigma^+)^*$ and $(\Sigma^-)^*$ 
be the free monoids generated by $\Sigma$, $\Sigma^+$
and $\Sigma^-$, respectively.
We refer to $\Sigma$ as an {\it alphabet},
its elements as  {\it letters}, and 
elements of $\Sigma^*$ as {\it words}.
Elements of $\Sigma^+$ and $\Sigma^-$ are referred to as {\it positive letters}
and {\it negative letters}, respectively, 
while 
elements of $(\Sigma^+)^*$ and $(\Sigma^-)^*$ are referred to as {\it positive words}
and {\it negative words}, respectively.
The length of a word $w$ is the number of letters it contains,
which is denoted by $|w|$.
The length of the null word, $\varepsilon$, is zero.
The null word is the identity of each monoid.

We write
the canonical monoid homomorphism
as $\pi : \Sigma^* \rightarrow G(p_1,\dots,p_n)$.
If $u$ and $v$ are words, then $u=v$ means $\pi(u)=\pi(v)$,
and
$u \equiv v$ means that
$u$ and $v$ are identical letter by letter.
A word $w \in \pi^{-1}(g)$ is called a {\it representative} of $g$.
The length of a group element $g$ is defined by the quantity
\[
 \Vert g \Vert :=\mbox{min}\{k \ ; \ g=\pi(\sigma_1\cdots \sigma_k),\ \sigma_i \in \Sigma\}.
\]
A word $w \in \Sigma^*$ is {\it geodesic} if
$|w|=\Vert \pi(w) \Vert$.
A word $\sigma_1\cdots \sigma_k \in \Sigma^*$ is called a {\it reduced} word
if $\sigma_i \neq \sigma^{-1}_{i+1}$ for all $i \in \{1,\dots,k-1\}$.
A geodesic representative is a reduced word.

Now, we define
the {\it spherical growth series} of $G(p_1,\dots,p_n)$ with respect to $\Sigma$ as the following formal power series:
\begin{eqnarray}\label{def-growth-series}
\begin{array}{ccl}
{\cal S}_{G(p_1,\dots,p_n)}(t)&:=&
\displaystyle{\sum_{l=0}^{\infty}}\ 
\sharp
\{
g \in G(p_1,\dots,p_n) \ ; \
\Vert g \Vert =l
\}\ t^l.\\
\end{array} 
\end{eqnarray}
It is well known that the radius of convergence of the spherical growth series of 
any finitely generated group is positive
(cf. \cite{C-D-Pap, H}).
Thus, ${\cal S}_{G(p_1,\dots,p_n)}(t)$ is a holomorphic function near the origin, 0.

The group structure of $G(p_1,\dots,p_n)$ is quite simple,
as can be understood from the existence of the exact sequence (\ref{exact-sq})
and the fact that 
$\overline{G}(p_1,\dots,p_n)$ is isomorphic to ${\bf Z}/p_1 {\bf Z} *\cdots * {\bf Z}/p_n {\bf Z}$ as a group.
In this paper, however,
we wish to elucidate the combinatorial group structure of $G(p_1,\dots,p_n)$
with respect to the presentation (\ref{group-presentation}).
In particular,
we investigate the spherical growth series of $G(p_1,\dots,p_n)$
with respect to the generating set
$ \Sigma=\{x_1,\dots,x_n,x_1^{-1},\dots,x_n^{-1}  \}$
and derive a rational function expression for 
the spherical growth series ${\cal S}_{G(p_1,\dots,p_n)}(z)$ of
$G(p_1,\dots,p_n)$ with respect to $\Sigma$.
\\

In the remainder of this section,
we introduce a unique particular representative for an element $g$ of $G(p_1,\dots,p_n)$,
which we call the {\it modified normal form} of $g$.

First, we introduce the positive word
\begin{eqnarray}\label{Delta}
\Delta :=x_1^{p_1} \in (\Sigma^+)^*.
\end{eqnarray}
Then,
$\pi(\Delta)$ is a generator of the center of the group $G(p_1,\dots,p_n)$.
Specifically, we have
\begin{eqnarray}\label{central}
\sigma \Delta^{\pm 1} = \Delta^{\pm 1}  \sigma \quad \mbox{for} \quad \sigma \in \Sigma=\{x_1,\dots,x_n,x_1^{-1},\dots,x_n^{-1}\}.
\end{eqnarray}

\vspace{0.1cm}

Let $g$ be an element of $G(p_1,\dots,p_n)$.
Consider an arbitrary representative $w$ of $g$
and suppose that it contains $l_+$ instances of a positive word $u$
satisfying $u = \Delta$
and $l_-$ instances of a negative word $v$ satisfying $v = \Delta^{-1}$.
Then, from (\ref{central}),
we can obtain a distinct representative $w'$ of $g$
by moving each of these instances (in arbitrary order)
to the rightmost position of $w$.
Let us write the word $w'$ as $w' \equiv \sigma_1 \cdots \sigma_{k} \cdot \Delta^c$,
where $c:=l_+ - l_- (\in {\bf Z})$ and $\sigma_j \in \Sigma$ for each $1 \leq j \leq k$.
If there exists $j$ such that $\sigma_{j+1} \equiv \sigma_j^{-1}$,
then we reduce $\sigma_j \cdot \sigma_{j+1}$ to the null word, $\varepsilon$.
Repeating this reduction procedure,
we obtain a reduced word $w''$ satisfying $w'' = \sigma_1 \cdots \sigma_k$.
This yields a representative $\lambda$ given by 
$\lambda :=w'' \cdot  \Delta^c$.
We present $\lambda$ as
\begin{eqnarray}\label{representative}
\lambda  \equiv  
x_{i_1}^{a_{1}}\cdot x_{i_2}^{a_{2}}\cdots x_{i_\tau}^{a_{\tau}}
\cdot \Delta^c,
\end{eqnarray}
where 
\begin{eqnarray}\label{tau}
\left\{
\begin{array}{l}
\tau \in {\bf N} \cup \{0\},\\
i_j \in \{1,\dots,n\}\ (1\leq j \leq \tau), 
i_j \neq i_{j+1}\ (1\leq j \leq \tau-1),
\end{array}
\right.
\end{eqnarray}
and 
\begin{eqnarray}\label{condition-representative}
-(p_{i_j}-1) \leq  a_{j} \leq p_{i_j}-1 \ (1 \leq j \leq \tau),\
a_{j} \neq 0 \ (1 \leq  j \leq \tau).
\end{eqnarray}
If $\tau=0$,
we stipulate that 
 the sub-word $x_{i_1}^{a_{1}}\cdot x_{i_2}^{a_{2}}\cdots x_{i_\tau}^{a_{\tau}}$
is the null word, $\varepsilon$.
We note that $x_{i_1}^{a_{1}}\cdot x_{i_2}^{a_{2}}\cdots x_{i_\tau}^{a_{\tau}}$
is a reduced word, and that $|\lambda| \leq |w|$.

\vspace{0.2cm}

Next, we define the {\it positive monoid} 
$G(p_1,\dots,p_n)^+$
associated with the group $G(p_1,\dots,p_n)$. 
This monoid is presented by
\[
G(p_1,\dots,p_n)^+\ =\ \langle\ 
 \ x_1,\dots, x_n\ |\ x_1^{p_1} =\cdots =x_n^{p_n}\
 \rangle^+,
\]
where the right-hand side is the quotient of the free monoid 
$(\Sigma^+)^*$ by an equivalence relation on $(\Sigma^+)^*$ defined as follows: 
(i) two positive words $w,w'\in (\Sigma^+)^*$ are {\it elementarily equivalent}
if there exist positive words $u,v\in (\Sigma^+)^*$ 
such that $w \equiv u \cdot x_{i_j}^{p_j} \cdot v$ 
and $w' \equiv u \cdot x_{i_{j'}}^{p_{j'}} \cdot v$;
(ii) two positive words $w,w'\in (\Sigma^+)^*$ are equivalent 
if there exists a sequence $w_0 \equiv w, w_1,\dots, w_l \equiv w'$ 
for some $l \in {\bf N} \cup \{0\}$ 
such that $w_s$ is elementarily equivalent to $w_{s+1}$ 
for $s=0,\dots,l-1$.

Regarding the monoid $G(p_1,\dots,p_n)^+$, Dehornoy and Paris showed the following:

\begin{pr}[Example 4 in Section 5 of \cite{D-P}]\label{pr-DP}
$G(p_1,\dots,p_n)^+$ is a Garside monoid whose fundamental element is $\Delta$.
\end{pr}

\noindent
See \cite{D-P}  for the definition of Garside monoids
and background discussion concerning them.
Also, let us note that Garside monoids are defined differently in \cite{D},
where a weaker condition is employed.
Additional discussion can be found in \cite{C-M}.

\vspace{0.1cm}

If the exponent $a_{j}$ of 
the representative $\lambda$ of $g \in G(p_1,\dots,p_n)$ appearing in (\ref{representative}) is negative,
we replace $x_{i_j}^{a_{j}}$ with
$x_{i_j}^{p_{i_j}+a_{j}} \cdot \Delta^{-1}$.
Then, moving every $\Delta^{-1}$ to the rightmost position, we obtain a representative $\mu$ as in the following proposition.
The uniqueness of $\mu$ is implied by Proposition \ref{pr-DP}.

\begin{pr}\label{pr-nf}
For each $g \in G(p_1,\dots,p_n)$,
there exists a unique representative of $g$ given by
\begin{eqnarray}\label{nf}
\mu:=
x_{i_1}^{\alpha_{1}}\cdot x_{i_2}^{\alpha_{2}}\cdots x_{i_\tau}^{\alpha_{\tau}}
\cdot \Delta^{d},
\end{eqnarray}
where $d \in {\bf Z}$, and the condition in (\ref{tau}) and the following condition hold:
\begin{eqnarray}\label{cond-nf}
0 \leq \alpha_{j} \leq p_{i_j}-1 \ (1 \leq j \leq \tau),\
\alpha_i \neq 0 \ (1 \leq j \leq \tau).
\end{eqnarray}

\end{pr}

\vspace{0.1cm}

\noindent
The representative of $g \in G(p_1,\dots,p_n)$ in Proposition \ref{pr-nf}, $\mu$,
corresponds to the standard form (resp., the normal form) 
for an element of the braid group (resp., the Artin group of finite type)
(see \cite{Garside} and \cite{B-S}).
We call $\mu$ 
the {\it Garside normal form} of $g$
and 
$x_{i_1}^{\alpha_{1}}\cdot x_{i_2}^{\alpha_{2}}\cdots x_{i_\tau}^{\alpha_{\tau}}$
the {\it non-$\Delta$ part} of the Garside normal form.

Next, we obtain another representative of $g \in G(p_1,\dots,p_n)$
from the Garside normal form, $\mu$,
defined in Proposition \ref{pr-nf}.
First, recall the symbols $p_k^-$ and $p_k^+$, defined in (\ref{p-+}).
If there exists a sub-word $x_{i_j}^{\alpha_{j}}$ 
in the non-$\Delta$ part of $\mu$
satisfying
$p_{i_j}^++1 \leq \alpha_{j}$,
replacing all such sub-words $x_{i_j}^{\alpha_{j}}$ with $x_{i_j}^{\alpha_{j}-p_{i_j}} \Delta$
and moving all such $\Delta$ to the rightmost position of $\mu$,
we obtain the representative
\begin{eqnarray}\label{mnf}
\nu :=
x_{i_1}^{\overline{\alpha}_{1}}\cdot x_{i_2}^{\overline{\alpha}_{2}}\cdots x_{i_\tau}^{\overline{\alpha}_{\tau}}
\cdot \Delta^{d+\rho},
\end{eqnarray}
where $\rho$ is the number of quantities $\alpha_{j}$ satisfying
$p_{i_j}^++1 \leq \alpha_{j}$,
and $\overline{\alpha}_{j}$ is defined to be $\alpha_{j}-p_{i_j}$
in the case $p_{i_j}^+ +1 \leq \alpha_{j}$
and $\alpha_{j}$
in the case
$\alpha_{j} \leq p_{i_j}^+$.
Then, 
the exponent $\overline{\alpha}_{j}$ satisfies
the following condition:
\begin{eqnarray}\label{cond-mnf}
-p_{i_j}^- \leq \overline{\alpha}_{j} \leq  p_{i_j}^+ \quad(1 \leq j \leq \tau),\quad 
\overline{\alpha}_{j}\neq 0 \quad (1 \leq j \leq \tau).
\end{eqnarray}
We refer to the representative $\nu$ as the {\it modified normal form} of $g \in G(p_1,\dots,p_n)$
and 
to $x_{i_1}^{\overline{\alpha}_{1}}\cdot x_{i_2}^{\overline{\alpha}_{2}}\cdots x_{i_\tau}^{\overline{\alpha}_{\tau}}$
as the {\it non-$\Delta$ part} of  the modified normal form.
From the uniqueness of the Garside normal form,
we immediately obtain the following:

\begin{pr}\label{pr-mnf}
The modified normal form exists and is unique for each element of
$G(p_1,\dots,p_n)$.
\end{pr}

For the modified normal form $\nu$,
we define the following sets and quantities:
\begin{eqnarray}\label{def-R}
\left\{
\begin{array}{rcl}
R_\nu&:=&
\{j \ ; \ 
\FLm+1 \leq \overline{\alpha}_{j} \leq \FLp \ \},\\
{\bf r}_\nu&:=&\# R_\nu.
\end{array}
\right.
\end{eqnarray}

\begin{ex}\label{example-modified-normal-form}
{\rm
Let us consider the case $n=3$, with $p_1=3, p_2=6, p_3=7$.
Then, we have
\begin{eqnarray*}
p_1^-=0, \ p_1^+=2; \quad
p_2^-=1, \ p_2^+=4;\quad
p_3^-=2, \ p_3^+=4.
\end{eqnarray*}

Let us consider an element $g \in G(3,6,7)$ that has the following representative:
\begin{eqnarray*}
\begin{array}{rcl}
\lambda &\equiv&
x_2^2 \cdot x_3^{-3} \cdot x_1^{-2} \cdot x_3^4
 \cdot x_2^5 \cdot x_3^{-2} 
\cdot x_2^4  \cdot x_3^2 \cdot x_1^2 \cdot x_2^4 \cdot
x_3^4 \cdot x_2^2 \cdot \Delta^{-3}\\
\\
&\equiv&
x_2^{a_1} \cdot  x_3^{a_2} 
 \cdot x_1^{a_3} 
 \cdot x_3^{a_4} \cdot x_2^{a_5} 
\cdot x_3^{a_6}  \cdot x_2^{a_{7}} 
\cdot x_3^{a_{8}} 
\cdot x_1^{a_{9}} \cdot x_2^{a_{10}} 
\cdot x_3^{a_{11}} 
\cdot x_2^{a_{12}} \cdot \Delta^{c}.
\end{array}
\end{eqnarray*}
Then the Garside normal form of $g$ is given by
\begin{eqnarray*}
\begin{array}{rcl}
\mu &\equiv&
x_2^2 \cdot x_3^4 \cdot x_1 \cdot x_3^4
 \cdot x_2^5 \cdot x_3^5 
\cdot x_2^4  \cdot x_3^2 \cdot x_1^2 \cdot x_2^4 \cdot
x_3^4 \cdot x_2^2 \cdot \Delta^{-6}\\
\\
&\equiv&
x_2^{{\alpha}_1} \cdot  x_3^{{\alpha}_2} 
 \cdot x_1^{{\alpha}_3} 
 \cdot x_3^{{\alpha}_4} \cdot x_2^{{\alpha}_5} 
\cdot x_3^{{\alpha}_6}  \cdot x_2^{{\alpha}_{7}} 
\cdot x_3^{{\alpha}_{8}}  
\cdot x_1^{{\alpha}_{9}} \cdot x_2^{{\alpha}_{10}} 
\cdot x_3^{{\alpha}_{11}} 
\cdot x_2^{{\alpha}_{12}} \cdot \Delta^d,
\end{array}
\end{eqnarray*}
and the modified normal form of $g$ is given by
\begin{eqnarray*}
\begin{array}{rcl}
\nu &\equiv&
x_2^2 \cdot x_3^4  \cdot x_1 \cdot x_3^4
 \cdot x_2^{-1} \cdot x_3^{-2}
\cdot x_2^4  \cdot x_3^2 \cdot x_1^2 \cdot x_2^4 \cdot
x_3^4 \cdot x_2^2 \cdot \Delta^{-4}\\
\\
&\equiv&
x_2^{\overline{\alpha}_1} \cdot  x_3^{\overline{\alpha}_2} 
 \cdot x_1^{\overline{\alpha}_3} 
 \cdot x_3^{\overline{\alpha}_4} \cdot x_2^{\overline{\alpha}_5} 
\cdot x_3^{\overline{\alpha}_6}  \cdot x_2^{\overline{\alpha}_{7}} 
\cdot x_3^{\overline{\alpha}_{8}} 
\cdot x_1^{\overline{\alpha}_{9}} \cdot x_2^{\overline{\alpha}_{10}} 
\cdot x_3^{\overline{\alpha}_{11}} 
\cdot x_2^{\overline{\alpha}_{12}} \cdot \Delta^{d+\rho}.
\end{array}
\end{eqnarray*}
Thus, we have
\begin{eqnarray*}
R_\nu=\{1,2,3,4,7,9,10,11,12\}\quad \mbox{and}\quad
{\bf r}_\nu=9. \ \Box
\end{eqnarray*}

}
\end{ex}


\section{Necessary and sufficient conditions for words to be geodesic,
and the suitable-spread procedure}\label{section-ns-cond}

In this section,
we present necessary and sufficient conditions for words in $\Sigma^*$ to be geodesic.
We also introduce a procedure 
that yields geodesic representatives for elements of $G(p_1,\dots,p_n)$,
which we call the {\it suitable-spread procedure}.

First,
we present the following lemma concerning geodesic representatives:

\begin{lm}\label{lm-sr}
Let $g \in G(p_1,\dots,p_n)$ and $w$ be a geodesic representative of $g$.
Then
there exists a geodesic representative of $g$ expressed as 
\begin{eqnarray}\label{sr}
\lambda:=
x_{i_1}^{a_{1}}\cdot x_{i_2}^{a_{2}}\cdots x_{i_\tau}^{a_{\tau}}
 \cdot \Delta^{c},
\end{eqnarray}
where the condition in (\ref{tau}) holds,
and $\lambda$ satisfies one of the following three conditions (\ref{cond-T1-geod})-(\ref{cond-T3-geod}):
\begin{eqnarray}
&\bullet& 
\left\{
\begin{array}{l}
c \in {\bf Z},\ c  > 0,\\
-p_{i_j}^- \leq a_{j} \leq  p_{i_j}^+ 
\quad (1 \leq j \leq \tau) ,\quad
a_{j} \neq 0 \quad (1 \leq j \leq \tau)
\end{array}
\right.\label{cond-T1-geod}\\
&\bullet&
\left\{
\begin{array}{l}
c \in {\bf Z},\ c <0,\\
-p_{i_j}^+ \leq a_{j} \leq p_{i_j}^- 
\quad  (1 \leq j \leq \tau),\quad
a_{j} \neq 0 \quad (1 \leq j \leq \tau)
\end{array}
\right.\label{cond-T2-geod}\\
&\bullet& 
\left\{
\begin{array}{l}
c =0,\\
-p_{i_j}^+ \leq a_{j} \leq  p_{i_j}^+
\quad  (1 \leq j \leq \tau),\quad
a_{j} \neq 0 \quad (1 \leq j \leq \tau).
\end{array}
\right.\label{cond-T3-geod}
\end{eqnarray}
If $\lambda$ satisfies the condition (\ref{cond-T1-geod}) (resp., (\ref{cond-T2-geod})),
then $w$ contains at least one positive word $u \ (=\Delta)$ (resp., negative word $v \ (=\Delta^{-1})$)
and $\lambda$ is obtained from $w$ by moving every positive word $u\ (=\Delta)$ (resp., negative word $v\ (=\Delta^{-1})$)
in $w$ to the rightmost position of $w$.
\end{lm}

\noindent
{\it Proof.}
Consider a representative $\lambda$ of $g$ 
as specified by (\ref{representative}) with (\ref{condition-representative}).
Because $w$ is geodesic,
a reduction of words does not occur during the procedure carried out in Section 2 in which $\lambda$ is obtained from $w$.
Hence, $\lambda$ is also geodesic.

We only consider the case in which $p_1$ is even and $p_k$ $(k \in \{2,\dots,n\})$ is odd,
because the proof for the other cases can be carried out similarly.

Let $p_1=2a$ and $p_k=2b+1$, 
where $a, b \in {\bf N}$ and $a \leq b$.
It is readily seen that 
$x_k^{2b}=x_k^{-1} \Delta, 
x_k^{2b-1}=x_k^{-2} \Delta, \dots, x_k^{a+b+1} = x_k^{- (b-a)} \Delta$
and $x_k^{-2b}=x_k \Delta^{-1}, \dots, x_k^{-(a+b+1)} =x_k^{b-a} \Delta^{-1}$.
If $b \geq a+1$,
then we have
$|x_k^{2b}| > |x_k^{-1} \Delta|, \dots, |x_k^{a+b+1}| > |x_k^{-(b-a)} \Delta|$
and $|x_k^{-2b}| > |x_k \Delta^{-1}|, \dots, |x_k^{-(a+b+1)}| > |x_k^{b-a} \Delta^{-1}|$.
Contrastingly if $b=a$, 
we do not have such inequalities.
Hence, if $b \geq a+1$ and $\lambda$ contains
sub-words $x_k^{2b}, \dots, x_k^{a+b+1}$ and $x_k^{-2b}=x_k, \dots, x_k^{-(a+b+1)}$, 
replacing $x_k^{2b}, \dots$, $x_k^{a+b+1}$ with
$x_k^{-1} \Delta, \dots$, $x_k^{-(b-a)} \Delta$, 
and also replacing $x_k^{-2b}, \dots$, $x_k^{-(a+b+1)}$ with
$x_k \Delta^{-1}, \dots$, $x_k^{b-a} \Delta^{-1}$, 
we obtain a representative $\lambda'$
that satisfies $|\lambda'| < |\lambda|$.
This is a contradiction, because $\lambda$ is geodesic.
Hence,
$\lambda$ contains no sub-word $x_k^A$ such that 
$-2b \leq A \leq  -(a+b+1)$ or $a+b+1 \leq A \leq 2b$.
Therefore, $\lambda$ takes the form
$\lambda \equiv x_{i_1}^{a_{1}}\cdots x_{i_\tau}^{a_{\tau}}\cdot \Delta^{c}$,
where $c$ and $a_j$ satisfy the following condition:
\begin{eqnarray}\label{cond-sr-1}
\left\{
\begin{array}{l}
c \in {\bf Z},\\
-(a+b) \leq  a_j \leq  a+b \ (1 \leq j \leq \tau),\ 
a_j\neq 0 \ (1 \leq j \leq \tau).\quad \\
\end{array}
\right.
\end{eqnarray}

Now, consider the case $c>0$.
If $\lambda$ contains a sub-word $x_k^{A}$ such that
$-(a+b) \leq A \leq -(b-a)-1$,
we move $\Delta$ in $\lambda$ so that it is the immediate right-hand neighbor of 
$x_k^{A}$,
and then replace $x_k^{A} \Delta$ with $x_k^{2b+1+A}$.
We thereby obtain a shorter representative of $g$, because 
$|x_k^{A}\Delta| > |x_k^{2b+1+A}|$.
This is a contradiction.
Hence, 
$\lambda$ contains no sub-word $x_k^{A}$ with $-(a+b) \leq A \leq -(b-a)-1$.
Therefore, $\lambda$ satisfies the condition (\ref{cond-T1-geod}),
because $\FLp =a+b$ and $\FLm =b-a$ in the present case.
Moreover, from the construction of $\lambda$, it is immediately seen that
$w$ contains at least one positive word $u \ (=\Delta)$ and that
$\lambda$ is obtained from $w$ by moving every positive word $u\ (=\Delta)$
in $w$ to the rightmost position of $w$.

In the case $c < 0$, 
a similar argument shows
that $\lambda$ satisfies the condition (\ref{cond-T2-geod}).

In the case $c=0$,
$\lambda$ satisfies the condition (\ref{cond-sr-1}),
which is the same as (\ref{cond-T3-geod}), because $\FLp =a+b$ in the present case.
$\Box$

\vspace{0.2cm}

Next, we consider the following three types of words.


\noindent
{\bf Type {\boldmath $1$}.}
A word $\xi \in \Sigma^*$ that satisfies
the following condition:
\begin{eqnarray*}\label{T1}
\begin{array}{rcl}
&(i)& \mbox{$\xi$ contains at least one positive word $u\ (=\Delta)$}.\\
&(ii)& \mbox{By moving every positive word $u\ (=\Delta)$ in $\xi$ to the rightmost position of $\xi$},\\
&&\mbox{we obtain a word 
$x_{i_1}^{a_{1}}\cdot x_{i_2}^{a_{2}}\cdots x_{i_\tau}^{a_{\tau}}
\cdot \Delta^{c} \ (\in \pi^{-1}(\pi(\xi)))$
that satisfies 
the conditions}\\
&&\mbox{in  (\ref{tau}) and (\ref{cond-T1-geod})}.\\
\end{array}
\end{eqnarray*}
If $\xi$ is of Type 1,
then
we have 
$|\xi| = 
|x_{i_1}^{a_{1}}\cdot x_{i_2}^{a_{2}}\cdots x_{i_\tau}^{a_{\tau}}
\cdot \Delta^{c}|$, and
the word $x_{i_1}^{a_{1}}\cdot x_{i_2}^{a_{2}}\cdots x_{i_\tau}^{a_{\tau}}
\cdot \Delta^{c}$ is the modified normal form
of $\pi(\xi) \in G(p_1,\dots,p_n)$.

\vspace{0.2cm}

\noindent
{\bf Type {\boldmath $2$}.}
A word $\xi \in \Sigma^*$ that satisfies
the following condition:
\begin{eqnarray*}\label{T2}
\begin{array}{rcl}
&(i)& \mbox{$\xi$ contains at least one negative word $v\ (=\Delta^{-1})$}.\\
&(ii)& \mbox{By moving every negative word $v\ (=\Delta^{-1})$ in $\xi$ to the rightmost position of $\xi$},\\
&&\mbox{we obtain a word 
$x_{i_1}^{a_{1}}\cdot x_{i_2}^{a_{2}}\cdots x_{i_\tau}^{a_{\tau}}
\cdot \Delta^{c} \ (\in \pi^{-1}(\pi(\xi)))$
that satisfies the conditions}\\
&&\mbox{in (\ref{tau}) and (\ref{cond-T2-geod})}.\\
\end{array}
\end{eqnarray*}
If $\xi$ is of Type 2,
then
we have 
$|\xi| = 
|x_{i_1}^{a_{1}}\cdot x_{i_2}^{a_{2}}\cdots x_{i_\tau}^{a_{\tau}}
 \cdot \Delta^{c}|$.

\vspace{0.2cm}

\noindent
{\bf Type {\boldmath $3$}.}
A word $\xi \in \Sigma^*$ that is presented as
\begin{eqnarray}\label{WT3}
\xi
\equiv 
x_{i_1}^{a_{1}}\cdot x_{i_2}^{a_{2}}\cdots x_{i_\tau}^{a_{\tau}},
\end{eqnarray}
where the conditions in (\ref{tau}) and (\ref{cond-T3-geod}) hold.
Here, for each $k\in \{1,\dots,n\}$,
we define
\[
\begin{array}{l}
{\rm Pos}_{x_k} (\xi):=\mbox{max}\{ a_{j}\ ; \  x_{i_j}=x_k,\ a_{j} \geq 0,\ 1 \leq j \leq \tau \ \},\\
{\rm Neg}_{x_k} (\xi):=\mbox{max}\{ - a_{j} \ ; \ x_{i_j}=x_k,\  a_{i _j}\leq 0,\ 1 \leq j \leq \tau\ \}.\\
\end{array}
\]
These quantities satisfy the relations
\[ 
\begin{array}{l}
0 \leq {\rm Pos}_{x_k}(\xi) \leq p_k^+,\
0 \leq {\rm Neg}_{x_k}(\xi) \leq p_k^+.\\
\end{array}
\]

Furthermore, for words of Type 3,
we define the following three sub-types.

\vspace{0.2cm}

\noindent
{\bf Type {\boldmath $3^+$}.}
A word $\xi \in \Sigma^*$ that is presented as in (\ref{WT3})
for which $a_{j}$ satisfies the condition in (\ref{cond-T1-geod}).
If $\xi$ is of Type $3^+$, 
then $\xi$ itself is the modified normal form of $\pi(\xi) \in G(p_1,\dots,p_n)$.


\noindent
{\bf Type {\boldmath $3^-$}.}
A word $\xi \in \Sigma^*$ that is presented as in (\ref{WT3})
for which
$a_{j}$ satisfies the condition in (\ref{cond-T2-geod}).


\noindent
{\bf Type {\boldmath $3^0$}.}
A word $\xi$ of Type 3
that is 
neither of Type $3^+$ nor of Type $3^-$.
\\

Next,
for $I \in \{1,2,3,3^+,3^-,3^0\}$, we define
\begin{eqnarray*}
{\rm WT}_I:=\{
\xi \in \Sigma^* \ ;
\
\xi \ \mbox{is a word of Type}\ I \ \}
\end{eqnarray*}
and introduce the following notation:
\begin{eqnarray*}
{\rm WT}_{3^+ \cup 3^-}:={\rm WT}_{3^+} \cup {\rm WT}_{3^-},\ \ 
{\rm WT}_{3^- \setminus 3^+}:={\rm WT}_{3^-} \setminus {\rm WT}_{3^+}.
\end{eqnarray*}
We call an element of ${\rm WT}_{3^+ \cup 3^-}$ (resp., ${\rm WT}_{3^- \setminus 3^+}$)
a {\it word of Type} $3^+ \cup 3^-$ (resp., {\it of Type} $3^- \setminus 3^+$).
We define ${\rm WT}$ to be the following disjoint union:
\begin{eqnarray*}\label{WT}
{\rm WT}:={\rm WT}_1 \sqcup {\rm WT}_2 \sqcup {\rm WT}_3.
\end{eqnarray*}
Then we have
\begin{eqnarray}\label{WT-1}
\begin{array}{rcl}
{\rm WT}&=&
{\rm WT}_1 \sqcup {\rm WT}_2 \sqcup {\rm WT}_{3^+ \cup 3^-} 
\sqcup {\rm WT}_{3^0}\quad \mbox{(disjoint union)}\\
&=&{\rm WT}_1 \sqcup {\rm WT}_2 \sqcup {\rm WT}_{3^+} \sqcup {\rm WT}_{3^- \setminus 3^+} 
\sqcup {\rm WT}_{3^0}\quad \mbox{(disjoint union)}.\\
\end{array}
\end{eqnarray}
From Lemma \ref{lm-sr},
it is seen that
for each element $g \in G(p_1,\dots,p_n)$,
any geodesic representative of $g$ belongs to ${\rm WT}$.

\begin{rmk}\label{rmk-wt}
{\rm We have
\begin{eqnarray*}
\begin{array}{rcl}
{\rm WT}_{3^+} \cap {\rm WT}_{3^-}
=\{
x_{i_1}^{a_{1}}\cdot x_{i_2}^{a_{2}}\cdots x_{i_\tau}^{a_{\tau}}
 & ; &
\mbox{the condition in (\ref{tau}) holds},\\
&&{\textstyle -\FLm 
\leq a_{j} \leq
\FLm} \ (1\leq j \leq \tau),\
a_{j}\neq 0 \ (1 \leq j \leq \tau) \},
\end{array}
\end{eqnarray*}
which is not empty,
because it necessarily contains the identity,
${\rm id}_{\Sigma^*} (=\varepsilon)$.
}
\end{rmk}

Let $\xi$ be a word of Type $I$ ($I \in \{1,2,3^+,3^- \setminus 3^+, 3^0\}$),
and consider the modified normal form $\nu$ for the element $\pi(\xi) \in G(p_1,\dots,p_n)$.
Recall the definitions of $\rho$ (appearing just below (\ref{mnf})) and 
${\bf r}_\nu$ (appearing in (\ref{def-R})).
Then, we find that 
the exponent $d+\rho$ and the quantity ${\bf r}_\nu$
satisfy the following:
\begin{eqnarray}\label{T1,2,3}
\left\{
\begin{array}{l}
\mbox{(i) If $\xi$ is of Type 1, then $d+\rho >0$}.\\
\mbox{(ii) If $\xi$ is of Type 2, then $d+\rho <0$ and ${\bf r}_\nu < -(d+\rho)$}.\\
\mbox{(iii${}^+$) If $\xi$ is of Type $3^+$, then $d+\rho=0$}.\\
\mbox{(iii${}^-$ $\setminus$ iii${}^+$) If $\xi$ is of Type $3^- \setminus 3^+$,
then $d+\rho <0$ and ${\bf r}_\nu =  -(d+\rho)$}.\\
\mbox{(iii${}^0$) If $\xi$ is of Type $3^0$, then $d+\rho<0$ and
${\bf r}_\nu > -(d+\rho)$}.\\
\end{array}
\right.
\end{eqnarray}

\begin{lm}\label{lm-wt}
Let $I,J \in \{1,2,3^+,3^- \setminus 3^+, 3^0\}$ with $I \neq J$.
In this case, we have
$\pi({\rm WT}_I) \cap \pi({\rm WT}_J)= \emptyset$.
\end{lm}

\noindent
{\it Proof.}
Suppose that there exists an element 
$g \in \pi({\rm WT}_I) \cap \pi({\rm WT}_J)$
with $I \neq J$.
Then
$g$ has representatives $\xi_I \in {\rm WT}_I$
and $\xi_J \in {\rm WT}_J$.
Let $\nu_{\xi_I}$ and $\nu_{\xi_J}$ be the modified normal forms
for $\xi_I$ and $\xi_J$, respectively.
From (\ref{T1,2,3}), we conclude that
$\nu_{\xi_I}$ is not identical to $\nu_{\xi_J}$.
This and the uniqueness of the modified normal form (Proposition \ref{pr-mnf})
 together yield a contradiction.
$\Box$

\vspace{0.2cm}

With the above preparation, we obtain the following proposition concerning
words in 
${\rm WT}_1 \sqcup {\rm WT}_2 \sqcup
{\rm WT}_{3^+ \cup 3^-}$.

\begin{pr}\label{pr-T1,2,3^+,3^-}
Let $\xi$ be a word such that
$\xi \equiv x_{i_1}^{a_{1}}\cdot x_{i_2}^{a_{2}}\cdots x_{i_\tau}^{a_{\tau}}
\cdot \Delta^{c}$.
Assume that $c \geq 0$ (resp., $c \leq 0$)
and that
$a_{j}$ satisfies the condition in (\ref{cond-T1-geod})
(resp., (\ref{cond-T2-geod})).
Then $\xi$ is geodesic.
\end{pr}

\noindent
{\it Proof.}
We only consider the case in which
$c \geq 0$ and 
$a_{j}$ satisfies the condition in (\ref{cond-T1-geod}),
because the proof for the other case can be carried out similarly.

First, let us suppose that $c >0$.
In this case, $\xi$ is an element of ${\rm WT}_1$.
Thus, from Lemma \ref{lm-wt},
we know that $\pi(\xi)$ has no representative in ${\rm WT}_2 \sqcup {\rm WT}_3$.
Hence, any geodesic representative $\gamma$ of $\pi(\xi)$
belongs to ${\rm WT}_1$.
This implies that
there exist $a'_{j}$ and $c'$
such that
$\gamma =
 x_{i_1}^{a'_{1}}
\cdots
x_{i_{\tau'}}^{a'_{{\tau'}}} \Delta^{c'}$
and 
$|\gamma| =
 |x_{i_1}^{a'_{1}}
\cdots
x_{i_{\tau'}}^{a'_{\tau'}}  \Delta^{c'}|$.
Therefore, because the modified normal forms of $\pi(\xi)$ and $\pi(\gamma)$ 
are identical
(Proposition \ref{pr-mnf}), 
we obtain
$ x_{i_1}^{a_{1}}
\cdots
x_{i_{\tau}}^{a_{\tau}} \Delta^{c}
\equiv
 x_{i_1}^{a'_{1}}
\cdots
x_{\tau'}^{a'_{\tau'}} \Delta^{c'}$.
Hence,
we have
$|\xi|=|\gamma|$, and thus
$\xi$ is geodesic.

The assertion is easily proved in the case $c=0$
by using the result for the case $c>0$,
because a sub-word of a geodesic word is also geodesic.
$\Box$

\vspace{0.2cm}

Proposition \ref{pr-T1,2,3^+,3^-}
implies that all words 
in ${\rm WT}_1 \sqcup {\rm WT}_2 \sqcup
{\rm WT}_{3^+ \cup 3^-}$are geodesic.

\vspace{0.2cm}

Lemma \ref{lm-sr} and the definition of ${\rm WT}$ imply that
any element $g \in G(p_1,\dots,p_n)$ has a geodesic representative $\xi \in {\rm WT}$.
With this in mind, for $I \in \{1,2,3,3^+,3^-,3^+ \cup 3^-, 3^- \setminus 3^-,3^0\}$,
we define
\begin{eqnarray}
{G(p_1,\dots,p_n)}_I:=\{g \in G(p_1,\dots,p_n) \ ; \ g  \ \mbox{has a geodesic representative of Type} \ I \},
\end{eqnarray}
and we call 
$g \in G(p_1,\dots,p_n)_I$ an element of Type $I$.
Then, from Lemma \ref{lm-wt}, we have
\begin{eqnarray}\label{decomp-G}
\begin{array}{rcl}
G(p_1,\dots,p_n)
&=&{G(p_1,\dots,p_n)}_1 \sqcup {G(p_1,\dots,p_n)}_2 \sqcup {G(p_1,\dots,p_n)}_{3^+} \sqcup {G(p_1,\dots,p_n)}_{3^- \setminus 3^+}\\
&&\sqcup {G(p_1,\dots,p_n)}_{3^0}\ \ \mbox{(disjoint union)}\\
&=&{G(p_1,\dots,p_n)}_1 \sqcup {G(p_1,\dots,p_n)}_2 \sqcup {G(p_1,\dots,p_n)}_{3^+ \cup 3^-}\\
 &&\sqcup {G(p_1,\dots,p_n)}_{3^0}
\ \ \mbox{(disjoint union)}\\
&=&{G(p_1,\dots,p_n)}_1 \sqcup {G(p_1,\dots,p_n)}_2 \sqcup {G(p_1,\dots,p_n)}_{3}
\ \ \mbox{(disjoint union)}.\\
\end{array}
\end{eqnarray}

\vspace{0.2cm}

Next, we present a necessary condition for words of Type 3
to be geodesic.

\begin{pr}\label{pr-nc-T3}
Let $\xi$ be a word in ${\rm WT}_3$.
If $\xi$ is geodesic,
then 
for every $k, k' \in \{1, \dots,n\}$,
we have the following inequality:
\begin{eqnarray}\label{ncond-T3}
\begin{array}{l}
{\rm Pos}_{x_k}(\xi)+{\rm Neg}_{x_{k'}}(\xi) \leq {\displaystyle \frac{p_k+p_{k'}}2}.
\end{array}
\end{eqnarray}
\end{pr}

\noindent
{\it Proof.}
Suppose that the inequality is not satisfied for some $k$ and $k'$,
i.e., that
${\rm Pos}_{x_k}(\xi)+{\rm Neg}_{x_{k'}}(\xi) > \frac{p_k+p_{k'}}{2}$.
Then there exist sub-words $x_{i_j}^{a_j}$
and $x_{i_{j'}}^{a_{j'}}$ in $\xi$
such that the following hold:
\begin{eqnarray}\label{i-j-3}
x_{i_j}=x_k,\ x_{i_{j'}}=x_{k'},\ a_j >0,\ a_{j'} <0,\
a_j - a_{j'}> {\textstyle \frac{p_k+p_{k'}}{2}}.
\end{eqnarray}
Suppose that $j < j'$ and consider the
sub-word of $\xi$ given by
\[
\xi^*:=x_{i_j}^{a_j} \cdot v \cdot   x_{i_{j'}}^{a_{j'}}.
\]
Because
$x_{i_{j'}}^{a_{j'}}=x_{i_{j'}}^{p_{k'}+a_{j'}} \Delta^{-1}$ and $\Delta$ satisfies 
 (\ref{central}),
we have
\begin{eqnarray}\label{xi'-2}
\xi^*=x_{i_j}^{a_j}\cdot v \cdot x_{i_{j'}}^{p_{k'}+a_{j'}}\Delta^{-1}
=x_{i_j}^{a_j} \Delta^{-1} \cdot v \cdot x_{i_{j'}}^{p_{k'}+a_{j'}}
=x_{i_j}^{a_j-p_k} \cdot v \cdot x_{i_{j'}}^{p_{k'}+a_{j'}}.
\end{eqnarray}
Further, from (\ref{i-j-3}),
we have
$-(a_j-p_k)+(p_{k'}+a_{j'}) < \frac{p_k+p_{k'}}{2}$,
and this implies that
$-(a_j-p_k)+(p_{k'}+a_{j'}) < a_j-a_{j'}$.
Thus, from (\ref{xi'-2}), we find that
$\xi^*$ is not geodesic,
and hence neither is
$\xi$.
The case $j > j'$ can be treated similarly,
with the conclusion again that $\xi$ is not geodesic.
$\Box$
\\

We now introduce a procedure that yields geodesic representatives.
This procedure, called the {\it suitable-spread procedure},
is similar to that employed in \cite{Berger}, \cite{M-M}, \cite{Fujii1} and \cite{Fujii2}.
Following \cite{Fujii2},
we explain the procedure in the case of the group $G(p_1,\dots,p_n)$.

\vspace{0.2cm}

For an element $g \in G(p_1,\dots,p_n)$,
let $\mu :=
x_{i_1}^{\alpha_{1}}\cdots x_{i_\tau}^{\alpha_{\tau}}
\cdot \Delta^d$
be its Garside normal form.
Then, 
consider the modified normal form of $g$,
\begin{eqnarray}\label{mnf-2}
\nu :=
x_{i_1}^{\overline{\alpha}_{1}}\cdots x_{i_\tau}^{\overline{\alpha}_{\tau}}
\cdot \Delta^{d+\rho}.
\end{eqnarray}
The word $\nu$ is the input of this procedure.
Recall that the exponents ${\overline{\alpha}}_{j}$ satisfy
the condition (\ref{cond-mnf})
and that $R_\nu$ and ${\bf r}_\nu$ are defined in (\ref{def-R}).

\vspace{0.2cm}

\noindent
{\bf [Case 1 : $d+\rho \geq 0$]}
In this case, 
$\nu$ belongs to ${\rm WT}_1 \sqcup {\rm WT}_{3^+}$.
Hence,
by Proposition \ref{pr-T1,2,3^+,3^-},
$\nu$ itself is geodesic.
For this reason, we regard $\nu$ as the output of this procedure.

\vspace{0.2cm}

\noindent
{\bf [Case 2 : $d+\rho <0,\ {\bf r}_\nu \leq -(d+\rho)$]}
First, using $\delta := -(d+\rho) > 0$,
let us rewrite (\ref{mnf-2}) as
\[
\nu \equiv 
x_{i_1}^{\overline{\alpha}_{1}}\cdots x_{i_\tau}^{\overline{\alpha}_{\tau}}
 \cdot (\Delta^{-1})^{\delta}.
\]
Then, if $j \in R_\nu$,
moving one of the words $\Delta^{-1}$
so that it is the immediate right-hand neighbor
of $x_{i_j}^{\overline{\alpha}_{j}}$, 
we obtain the following representative of $g$:
\begin{eqnarray*}
\begin{array}{ccl}
&&
x_{i_1}^{\overline{\alpha}_{1}} \cdots (x_{i_j}^{\overline{\alpha}_{j}}\Delta^{-1}) \cdots
x_{i_\tau}^{\overline{\alpha}_{\tau}} 
 \cdot (\Delta^{-1})^{\delta-1}
\\
&&
=
x_{i_1}^{\overline{\alpha}_{1}} \cdots (x_{i_j}^{\overline{\alpha}_{j}-p_{i_j}}) \cdots
x_{i_\tau}^{\overline{\alpha}_{\tau}} 
 \cdot (\Delta^{-1})^{\delta-1}
=: \widetilde{\nu}_{\overline{\alpha}_{j}}.\\
\end{array}
\end{eqnarray*}
Note the equality
\begin{eqnarray}\label{equality-ss-1}
|\widetilde{\nu}_{\overline{\alpha}_{j}}|=|\nu|-(2\overline{\alpha}_{j}+p_1-p_{i_j}).
\end{eqnarray}

Because ${\bf r}_\nu \leq \delta$,
carrying out the above procedure for all $j \in R_\nu$,
we obtain a representative $\widetilde{\nu}$ of $g$. We write it as
\begin{eqnarray*}
\widetilde{\nu} \equiv 
x_{i_1}^{\overline{\alpha}'_{1}}\cdots x_{i_\tau}^{\overline{\alpha}'_{\tau}}
 \cdot
(\Delta^{-1})^{\delta - {\bf r}_\nu}.
\end{eqnarray*}
It is seen that $\overline{\alpha}'_{j}$ satisfies the condition
\begin{eqnarray*}\label{cond-case-2}
-\FLp \leq  \overline{\alpha}'_{j} \leq  \FLm \ (1 \leq j \leq \tau),
\quad
\overline{\alpha}'_{j} \neq 0 \ (1 \leq  j \leq \tau).
\end{eqnarray*}
The word $\widetilde{\nu}$ belongs to ${\rm WT}_2 \sqcup {\rm WT}_{3^- \setminus 3^+}$.
Hence,
from Proposition \ref{pr-T1,2,3^+,3^-}, 
we know that it is geodesic.
In this case, we regard $\widetilde{\nu}$ as the output of the procedure.

\vspace{0.2cm}

\noindent
{\bf [Case 3 : $d+\rho <0,\ {\bf r}_\nu > -(d+\rho)\ (=\delta)$]}
First, pick $\delta$ exponents
$\overline{\alpha}_{j_1},\dots, \overline{\alpha}_{j_\delta}$
from among the exponents of $\nu$
such that
$\sum_{t=1}^\delta (2 \overline{\alpha}_{j_t}+p_1-p_{i_{j_t}})$ realizes its maximal
value.
Then, applying the procedure in Case 2 to all $x_{i_{j_t}}^{\overline{\alpha}_{j_t}}$,
we obtain a representative $\widetilde{\nu}$ that contains
no $\Delta^{-1}$.
We present this representative as
\begin{eqnarray*}
\widetilde{\nu} \equiv
x_{i_1}^{\overline{\alpha}'_{1}}\cdots x_{i_\tau}^{\overline{\alpha}'_{\tau}}.
\end{eqnarray*}
Here, $\overline{\alpha}'_{j}$ satisfies the condition
\begin{eqnarray*}\label{cond-case-3}
-\FLp \leq  \overline{\alpha}'_{j} \leq  \FLp \ (1 \leq j \leq \tau),
\quad
\overline{\alpha}'_{j} \neq 0 \ (1 \leq j \leq \tau).
\end{eqnarray*}
The word $\widetilde{\nu}\ (\in {\rm WT}_3)$ is the output of the procedure.
In the present case,
from (\ref{equality-ss-1}), we obtain
\begin{eqnarray}\label{equality-ss-3}
|\widetilde{\nu}|=|\nu|- 
\sum_{t=1}^\delta (2 \overline{\alpha}_{j_t} +p_1-p_{i_{j_t}}).
\end{eqnarray}
Note that
the choice of the integers
$\overline{\alpha}_{j_1},\dots, \overline{\alpha}_{j_\delta}$
maximizing 
$\sum_{t=1}^\delta (2 \overline{\alpha}_{j_t}+p_1-p_{i_{j_t}})$
is not necessarily unique,
and
the output $\widetilde{\nu}$ depends on this choice.
We define the set of all such choices as follows:
\begin{eqnarray*}
\begin{array}{rl}
{{\rm CE}}_g :=
\{
(j_1,\dots,j_{\delta})\
;&
\nu \equiv 
x_{i_1}^{\overline{\alpha}_{1}}\cdots x_{i_\tau}^{\overline{\alpha}_{\tau}}
 \cdot \Delta^{-\delta}
\ (\mbox{modified normal form of}\ g),\
{\bf r}_\nu > \delta >0,\\
&
1 \leq j_1 < \cdots < j_{\delta} \leq \tau,\\
&
\mbox{and}\ 
\sum_{t=1}^\delta (2 \overline{\alpha}_{j_t}+p_1-p_{i_{j_t}})\
 \mbox{realizes its maximal value}
\}.
\end{array}
\end{eqnarray*}

\begin{ex}\label{example-SS}
{\rm 

Let us choose the element $\nu$ given in Example \ref{example-modified-normal-form},
which is the modified normal form of the element $g \in G(3,6,7)$.
This implies Case 3, because ${\bf r}_\nu=9$ and $\delta=4$.
We can verify the following:
\begin{eqnarray*}
{{\rm CE}}_g=\{
(2,4,7,10),
(2,7,10,11),
(4,7,10,11),
(2,7,9,10),
(4,7,9,10),
(7,9,10,11)
\}.
\end{eqnarray*}
Applying the suitable spread procedure in accordance with these choices,
we obtain
\begin{eqnarray*}
\begin{array}{rcl}
\nu&=&
x_2^2 \cdot (x_3^4\Delta^{-1}) \cdot x_1 \cdot 
(x_3^4\Delta^{-1})
 \cdot x_2^{-1} \cdot x_3^{-2} \cdot (x_2^4\Delta^{-1}) 
\cdot x_3^2 \cdot x_1^2 \cdot  (x_2^4\Delta^{-1}) \cdot
x_3^4 \cdot x_2^2
\\
&&=
x_2^2 \cdot x_3^{-3} \cdot x_1 \cdot x_3^{-3}
 \cdot x_2^{-1} \cdot x_3^{-2} 
\cdot x_2^{-2}  \cdot x_3^2 \cdot x_1^2 \cdot  x_2^{-2} \cdot
x_3^4 \cdot x_2^2
=:\widetilde{\nu}_1
\\
&=&
x_2^2 \cdot (x_3^4\Delta^{-1}) \cdot x_1 \cdot x_3^4
 \cdot x_2^{-1} \cdot x_3^{-2} \cdot (x_2^4\Delta^{-1}) 
\cdot x_3^2 \cdot x_1^2 \cdot  (x_2^4\Delta^{-1}) \cdot
 (x_3^4\Delta^{-1}) \cdot x_2^2
\\
&&=
x_2^2 \cdot x_3^{-3} \cdot x_1 \cdot x_3^4
 \cdot x_2^{-1} \cdot x_3^{-2} 
\cdot x_2^{-2}  \cdot x_3^2  \cdot x_1^2 \cdot x_2^{-2} \cdot
x_3^{-3} \cdot x_2^2
=:\widetilde{\nu}_2
\\
&=&
x_2^2 \cdot x_3^4 \cdot x_1 \cdot (x_3^4\Delta^{-1})
 \cdot x_2^{-1} \cdot x_3^{-2} \cdot (x_2^4\Delta^{-1}) 
\cdot x_3^2 \cdot x_1^2 \cdot (x_2^4\Delta^{-1}) \cdot
 (x_3^4\Delta^{-1}) \cdot x_2^2
\\
&&=
x_2^2 \cdot x_3^4 \cdot x_1 \cdot x_3^{-3}
 \cdot x_2^{-1} \cdot x_3^{-2} 
\cdot x_2^{-2}  \cdot x_3^2 \cdot x_1^2 \cdot x_2^{-2} \cdot
x_3^{-3} \cdot x_2^2
=:\widetilde{\nu}_3
\\
&=&
x_2^2 \cdot (x_3^4\Delta^{-1}) \cdot x_1 \cdot x_3^4
 \cdot x_2^{-1} \cdot x_3^{-2} \cdot (x_2^4\Delta^{-1}) 
\cdot x_3^2 \cdot (x_1^2\Delta^{-1}) \cdot (x_2^4\Delta^{-1}) \cdot
 x_3^4 \cdot x_2^2
\\
&&=
x_2^2 \cdot x_3^{-3} \cdot x_1 \cdot x_3^4
 \cdot x_2^{-1} \cdot x_3^{-2} 
\cdot x_2^{-2}  \cdot x_3^2 \cdot x_1^{-1} \cdot x_2^{-2} \cdot
x_3^4 \cdot x_2^2
=:\widetilde{\nu}_4
\\
&=&
x_2^2 \cdot x_3^4 \cdot x_1  \cdot (x_3^4\Delta^{-1})
 \cdot x_2^{-1} \cdot x_3^{-2} \cdot (x_2^4\Delta^{-1}) \cdot x_3^2 \cdot (x_1^2\Delta^{-1}) \cdot (x_2^4\Delta^{-1}) \cdot
 x_3^4 \cdot x_2^2
\\
&&=
x_2^2 \cdot x_3^4 \cdot x_1 \cdot x_3^{-3}
 \cdot x_2^{-1} \cdot x_3^{-2} 
\cdot x_2^{-2}  \cdot x_3^2 \cdot x_1^{-1} \cdot x_2^{-2} \cdot
x_3^4 \cdot x_2^2
=:\widetilde{\nu}_5
\\
&=&
x_2^2 \cdot x_3^4 \cdot x_1 \cdot x_3^4
 \cdot x_2^{-1} \cdot x_3^{-2} \cdot (x_2^4\Delta^{-1}) 
\cdot x_3^2 \cdot (x_1^2\Delta^{-1}) \cdot (x_2^4\Delta^{-1}) \cdot
 (x_3^4\Delta^{-1}) \cdot x_2^2
\\
&&=
x_2^2 \cdot x_3^4 \cdot x_1 \cdot x_3^4
 \cdot x_2^{-1} \cdot x_3^{-2} 
\cdot x_2^{-2}  \cdot x_3^2 \cdot x_1^{-1} \cdot x_2^{-2} \cdot
x_3^{-3} \cdot x_2^2
=:\widetilde{\nu}_6.
\\
\end{array}
\end{eqnarray*}
For this example,
the procedure given here implies 
the six outputs
$\widetilde{\nu}_1, \dots, \widetilde{\nu}_6$.
$\Box$

}
\end{ex}

\vspace{0.2cm}

The procedure described above is
the suitable-spread procedure.
Note that the output of this procedure is also a representative of $g$.
We denote by 
${\rm SS}_g$
the set consisting of all of the outputs of the suitable-spread procedure
applied to the modified normal form of $g$, $\nu$.
In Cases 1 and 2, we have
$\# {\rm SS}_g=1$.
In Case 3,
${\rm SS}_g$ is bijective with ${\rm CE}_g$.
The following provides the final piece for this procedure.


\begin{pr}\label{pr-ss-1}
 Let $g$ be an element of $G(p_1,\dots,p_n)$, and let
 $\mu \equiv
x_{i_1}^{\alpha_{1}}\cdots x_{i_\tau}^{\alpha_{\tau}}
\cdot \Delta^d$
and 
$\nu \equiv 
x_{i_1}^{\overline{\alpha}_{1}}\cdots x_{i_\tau}^{\overline{\alpha}_{\tau}}
\cdot \Delta^{d+\rho}$ 
be the Garside normal form and the modified normal form of $g$,
respectively.
Then, applying the suitable-spread procedure to $\nu$,
we obtain a geodesic representative of $g$.
\end{pr}

\noindent
{\it Proof.}
In the case $d+\rho \geq 0$ and in the case $d+\rho <0$ with 
${\bf r}_\nu \leq -(d+\rho)$,
it has been shown that the output of the suitable-spread procedure is geodesic.
Hence, it remains only to consider the case
$d+\rho < 0$ with ${\bf r}_\nu > -(d+\rho)$.

Now, note that there exists at least one geodesic representative of $g$.
Choose one of these, say $\gamma$.
Because $\gamma$ is geodesic,
from Lemma \ref{lm-sr} and (\ref{WT-1}),
one of the following holds:
(i) $\gamma$ is of Type 1;
(ii) $\gamma$ is of Type 2;
(iii${}^+$) $\gamma$ is of Type $3^+$;
(iii${}^-$ $\setminus$ iii${}^+$) $\gamma$ is of Type $3^- \setminus 3^+$;
(iii${}^0$) $\gamma$ is of Type $3^0$.
As in the proof of Lemma \ref{lm-wt},
using the uniqueness of the modified normal form
(Proposition \ref{pr-mnf}),
we can confirm that the modified normal form of $g=\pi(\gamma)$ satisfies the following condition
in each of the corresponding cases:
(i) $d+\rho >0$;
(ii) $d+\rho <0, {\bf r}_\nu < -(d+\rho)$;
(iii$^+$) $d+\rho = 0$;
(iii$^-$ $\setminus$ iii$^+$) $d+\rho < 0, {\bf r}_\nu = -(d+\rho)$;
(iii$^0$) $d+\rho < 0, {\bf r}_\nu > -(d+\rho)$.
Therefore, we need only treat Case (iii$^0$).

In Case (iii$^0$),
$\gamma$ takes the form
$\gamma \equiv x_{i_1}^{a_{1}}\cdots x_{i_T}^{a_{T}}$,
where $a_{j}$ satisfies the condition 
\begin{eqnarray*}
\begin{array}{l}
-p_{i_j}^+ \leq a_{j} \leq  p_{i_j}^+
\quad  (1 \leq j \leq T),\quad
a_{j} \neq 0 \quad (1 \leq j \leq T).
\end{array}
\end{eqnarray*}
Next, we obtain the Garside normal form, $\mu_\gamma$, for $g=\pi(\gamma)$
from the representative $\gamma$
as in Proposition \ref{pr-nf}.
We write $\mu_\gamma$ as
\begin{eqnarray*}
\mu_\gamma \equiv 
x_{i_1}^{A_{1}}\cdots x_{i_T}^{A_{T}}
\cdot \Delta^{D}.
\end{eqnarray*}
Here,
if $a_{j}<0$, 
then we have $A_{j}= p_{i_j}+a_{j}$,
and
if $a_{j} \geq 0$,
then we have $A_{j}=a_{j}$.
We also have $D=-({\bf r}^{(1)}_\gamma+{\bf r}^{(2)}_\gamma)$,
where we define
\begin{eqnarray*}
\left\{
\begin{array}{ll}
R^{(1)}_\gamma:= \{j \ ; \ -\FLm \leq a_{j} < 0 \},
& {\bf r}^{(1)}_\gamma:=\# R^{(1)}_\gamma,\\
R^{(2)}_\gamma:=\{j \ ; \ -\FLp \leq  a_{j} \leq -\FLm-1 \},
& {\bf r}^{(2)}_\gamma:=\# R^{(2)}_\gamma.\\
\end{array}
\right.
\end{eqnarray*}
By the uniqueness of the Garside normal form (Proposition \ref{pr-nf}),
for the quantities above,
we have $T=\tau$, $A_{j}=\alpha_{j}$ and $D=d$.
Also, because
$p_{i_j} -\FLm=\FLp+1$,
we have ${\bf r}^{(1)}_\gamma=\rho$.
Thus, the relation $-(d+\rho) = {\bf r}^{(2)}_\gamma$ holds,
and this implies that
$\delta={\bf r}^{(2)}_\gamma$.
Then, because $\gamma$ is of Type $3^0$,
we have ${\bf r}^{(2)}_\gamma >0$
and ${\bf r}_\nu > {\bf r}^{(2)}_\gamma$.
Hence, we find 
${\bf r}_\nu > \delta ={\bf r}^{(2)}_\gamma >0$.
It is readily seen that
the geodesic $\gamma$ is obtained from the modified normal form
$\nu$ by making the replacements
$x_{i_j}^{\overline{\alpha}_{j}} \rightarrow x_{i_j}^{\overline{\alpha}_{j}} \Delta^{-1}$
for all $j \in R^{(2)}_\gamma$
as in the suitable-spread procedure.
Moreover, it is clear that the quantity
\begin{eqnarray*}
\sum_{j \in R^{(2)}_\gamma} (2 \overline{\alpha}_{j}+p_1-p_{i_j})
\end{eqnarray*}
is less than or equal to the maximal value of the quantity
\begin{eqnarray*}
\sum_{t=1}^\delta (2 \overline{\alpha}_{j_t}+p_1-p_{i_{j_t}}).
\end{eqnarray*}
Thus, 
from (\ref{equality-ss-3}) and the assumption that $\gamma$ is geodesic,
it follows that they are equal.
Therefore,
we conclude that the suitable-spread procedure produces a word $\widetilde{\nu}$
of length $|\gamma|$.
Hence, this word is geodesic.
(Note that from the above argument,
it follows that 
if $\gamma$ is of Type $3^0$, then
$\gamma$ is obtained from the suitable-spread procedure applied to
$\nu$.)
$\Box$

\vspace{0.2cm}

Let $g \in G(p_1,\dots,p_n)_{3^+}$ (resp., $g \in G(p_1,\dots,p_n)_{3^- \setminus 3^+}$).
Then we have Case 1 (resp., Case 2),
and any geodesic representative $\gamma \in {\rm WT}_{3^+}$
(resp., $\gamma \in {\rm WT}_{3^- \setminus 3^+}$)
itself is the output of this procedure.
From this and
the argument given in the proof of Proposition \ref{pr-ss-1},
we obtain the following:

\begin{pr}\label{pr-ss-2}
Let $g$ be an element of $G(p_1,\dots,p_n)_3$ and
$\gamma \in {\rm WT}_3$ be a geodesic representative of $g$.
Then $\gamma$ is an output of the suitable-spread procedure.
\end{pr}

Next, we present a sufficient condition for words of Type 3 to be geodesic
and a necessary and sufficient condition for
elements of $G(p_1,\dots,p_n)_3$ to have a single geodesic representative.

\begin{pr}\label{pr-sc-T3}
\begin{enumerate}
\item
Let $\xi \equiv
x_{i_1}^{a_{1}}\cdots x_{i_\tau}^{a_{\tau}}$ be an element of ${\rm WT}_3$.
If $\xi$ satisfies the inequalities in (\ref{ncond-T3}),
then it is geodesic.
\item
Let $g$ be an element of $G(p_1,\dots,p_n)_3$.
Then
the set of all of the geodesic representatives of $g$
is identical to ${\rm SS}_g$.
\item
Let $g$ be an element of $G(p_1,\dots,p_n)_3$,
and let $\xi \equiv
x_{i_1}^{a_{1}}\cdots x_{i_\tau}^{a_{\tau}}
\in {\rm WT}_3$ be its geodesic representative.
If $\xi$ satisfies the inequalities
\begin{eqnarray}\label{inequality-T3}
{\rm Pos}_{x_k}(\xi)+{\rm Neg}_{x_{k'}}(\xi) < \frac{p_k+p_{k'}}{2}
\quad (1 \leq k \leq n, 1\leq k' \leq n),
\end{eqnarray}
then $\# {\rm SS}_g =1$.
\item
Let $g$ be an element of $G(p_1,\dots,p_n)_3$,
and let $\xi \equiv
x_{i_1}^{\alpha_{1}}\cdots x_{i_\tau}^{\alpha_{\tau}}
\in {\rm WT}_3$ be its geodesic representative.
Suppose that there appear some generators $x_k$ and $x_{k'}$ in $\xi$ such that the following equality holds:
\begin{eqnarray}\label{equality-T3}
{\rm Pos}_{x_k}(\xi)+{\rm Neg}_{x_{k'}}(\xi) = \frac{p_k+p_{k'}}{2}.
\end{eqnarray}
Then $\# {\rm SS}_g \geq 2$.
\end{enumerate}
\end{pr}

\noindent
{\it Proof.}
In the following,
for 
$\xi \equiv
x_{i_1}^{a_{1}}\cdots x_{i_\tau}^{a_{\tau}}
\in {\rm WT}_3$,
we use the notation $P_{x_k}:={\rm Pos}_{x_k}(\xi)$
and $N_{x_k}:={\rm Neg}_{x_k}(\xi)$ for each $k \in \{1,\dots,n\}$.

First, write the Garside normal form of $\pi(\xi)$ as
\begin{eqnarray*}
\mu \equiv 
x_{i_1}^{\alpha_{1}}\cdots x_{i_\tau}^{\alpha_{\tau}}
\cdot \Delta^{d}.
\end{eqnarray*}
Here, $d$ is given by
\begin{eqnarray}\label{d}
d=- (\# \{j \ ; \ a_{j} <0 \})
\leq 0,
\end{eqnarray}
and $\alpha_{j}$ 
satisfies the following:
\begin{eqnarray}\label{alpha-beta}
\left\{
\begin{array}{l}
a_{j} < 0 \Rightarrow \alpha_{j}=p_{i_j}+a_{j},\\
a_{j} \geq 0 \Rightarrow \alpha_{j}=a_{j}.\\
\end{array}
\right.
\end{eqnarray}
Next,
write the modified normal form of $\pi(\xi)$ as
\begin{eqnarray*}
\nu \equiv 
x_{i_1}^{\overline{\alpha}_{1}}\cdots x_{i_\tau}^{\overline{\alpha}_{\tau}}
\cdot \Delta^{d+\rho}.
\end{eqnarray*}
Here, $\rho$ is given by
\begin{eqnarray}\label{rho}
\rho=\# \{j \ ; {\textstyle \FLp}+1 \leq \alpha_{j}
\}
=\# \{j \ ; \ -{\textstyle \FLm} \leq a_{j} <0 \},
\end{eqnarray}
and
$\overline{\alpha}_{j}$ satisfies the following:
\begin{eqnarray}\label{beta'}
\left\{
\begin{array}{l}
-\FLp \leq a_{j} \leq 
-\FLm -1 \Rightarrow 
\overline{\alpha}_{j}=p_{i_j}+a_{j},\\
-\FLm \leq a_{j}\leq
\FLp \Rightarrow 
\overline{\alpha}_{j}=a_{j}.\\
\end{array}
\right.
\end{eqnarray}
(The second equality in (\ref{rho}) follows from the relation
$\FLp +1-p_{i_j}=-\FLm$ and the fact that
$\{j \ ; \ \FLp +1 \leq a_{j}  \}$ is empty.)
Moreover, recall that 
the set $R_\nu$ and the quantity ${\bf r}_\nu$ are defined in (\ref{def-R})
for the modified normal form $\nu$.

Next, note that the definition of $N_{x_k}$ implies the following:
\begin{eqnarray}\label{ab-1}
x_{i_j}=x_k,\ a_{j} < 0 \Rightarrow a_{j} \geq -N_{x_k}.
\end{eqnarray}
Further,
the definition of $P_{x_k}$ implies 
\begin{eqnarray}\label{ab-2}
x_{i_j}=x_k,\ a_{j} \geq 0 \Rightarrow \overline{\alpha}_{j} =
 \alpha_{j}= a_{j} \leq P_{x_k},
\end{eqnarray}
as follows from (\ref{cond-T3-geod})
and the relations $0 \leq a_{j}\leq \FLp$ and $\alpha_{j}=a_{j}$.
It is thus seen that 
from (\ref{d}) and (\ref{rho}),
we have
\begin{eqnarray}\label{d+rho}
-(d+\rho)=\# R^{(2)}_\xi \geq 0,
\end{eqnarray}
where we define 
\begin{eqnarray}\label{rab}
R^{(2)}_\xi:=
\{j \ ; \ -\FLp
\leq a_{j} \leq -\FLm -1 \}.
\end{eqnarray}
Moreover,  from (\ref{alpha-beta}), (\ref{beta'})
and the relation $\FLp+1-p_{i_j}=-\FLm$,
it is readily confirmed that
\begin{eqnarray}\label{ineq-rab}
R_\nu \supseteq R^{(2)}_\xi.
\end{eqnarray}
This and (\ref{d+rho}) imply
\begin{eqnarray}\label{ineq-r}
{\bf r}_\nu \geq -(d+\rho).
\end{eqnarray}

With the above preparation,
we now demonstrate the four assertions of the proposition one by one.

\vspace{0.2cm}

\noindent
1.\ 
Written in terms of $P_{x_{k'}}$ and $N_{x_{k}}$,
the inequality in (\ref{ncond-T3}) becomes
\begin{eqnarray}\label{ineq-Prop3.5}
P_{x_{k'}} + N_{x_{k} }\leq \frac{p_{k'}+p_{k}}{2}.
\end{eqnarray}
By Proposition \ref{pr-T1,2,3^+,3^-},
if $\xi \in {\rm WT}_{3^+ \cup 3^-}$,
then $\xi$ is geodesic.
Let us assume here that $\xi \notin  {\rm WT}_{3^+ \cup 3^-}$.
From this assumption, there exist $j$ and $j'$ such that
$-\FLp \leq a_{j} \leq -\FLm -1$
and $p_{i_{j'}}^- +1 \leq a_{j'} \leq 
p_{i_{j'}}^+$.
Let $x_{i_j}=x_k$ and $x_{i_{j'}}=x_{k'}$.
Then, from (\ref{ab-2}), (\ref{ineq-Prop3.5}), (\ref{ab-1}) and  (\ref{beta'}),
we obtain
\begin{eqnarray}\label{ineq-sub1}
\begin{array}{rcl}
2 {\overline{\alpha}}_{j'}+p_1-p_{k'} = 2 a_{j'} +p_1-p_{k'}
&\leq& 2 P_{x_{k'}}+p_1-p_{k'}\\
&\leq& p_1+p_k-2 N_{x_{k}}\\
&\leq& p_1+p_k+2 a_j\\
&=& 2 {\overline{\alpha}}_{j}+p_1-p_k.
\end{array}
\end{eqnarray}

Now, let us apply the suitable-spread procedure to 
the modified normal form $\nu$,
noting that the relations (\ref{d+rho}), (\ref{ineq-r})
and $\delta=\#R^{(2)}_\xi$ hold.
From (\ref{ineq-sub1}),
it is seen that the choice of all $x_{i_j}^{\overline{\alpha}_{j}}$ with $j \in R^{(2)}_\xi$
maximizes
$\sum_{t=1}^{\delta} (2 {\overline{\alpha}}_{j_t}+p_1-p_{i_{j_t}})$,
and thus this is one possible choice.
Proceeding with this choice,
we obtain $\xi$ exactly.
Thus, by Proposition \ref{pr-ss-1}, $\xi$ is geodesic.

\vspace{0.2cm}

\noindent
2.\ 
First, note that from (\ref{decomp-G}),
$g$ has no geodesic representative of Type 1 or Type 2.
Next, recall that
Proposition \ref{pr-ss-2} asserts that
all of the geodesic representatives of Type 3
can be obtained by applying the suitable-spread procedure to $\nu$.
We thus conclude that
any geodesic representative of $g$
is an element of ${\rm SS}_g$.
Conversely, by Proposition \ref{pr-ss-1},
any element of ${\rm SS}_g$ is a geodesic representative of $g$.
Therefore,
the set consisting of all of the geodesic representatives of $g$
is identical to ${\rm SS}_g$.

\vspace{0.2cm}

\noindent
3.\ 
If 
$\xi$ satisfies all the inequalities in (\ref{inequality-T3}),
we obtain
\begin{eqnarray}\label{ineq-sub2}
\begin{array}{rcl}
2 {\overline{\alpha}}_{j'}+p_1-p_{k'} = 2 a_{j'} +p_1-p_{k'}
&\leq& 2 P_{x_{k'}}+p_1-p_{k'}\\
&<& p_1+p_k-2 N_{x_{k}}\\
&\leq& p_1+p_k+2 a_j\\
&=& 2 {\overline{\alpha}}_{j}+p_1-p_k.
\end{array}
\end{eqnarray}
Hence, 
in contrast to the situation considered in the proof of the first assertion,
only the choice of all
$x_{i_j}^{{\overline{\alpha}}_j}$ with $j \in R^{(2)}_\xi$ 
maximizes
$\sum_{t=1}^{\delta} (2 {\overline{\alpha}}_{j_t}+p_1-p_{i_{j_t}})$.
Thus this is the only possible choice in this case, and therefore
$\xi$ is the unique geodesic representative obtained from the suitable-spread procedure.
Hence,
we have ${\rm SS}_g=\{\xi\}$ and $\# {\rm SS}_g =1$.

\vspace{0.2cm}

\noindent
4.\ 
By Proposition \ref{pr-nc-T3},
$\xi$ satisfies the inequalities in (\ref{ncond-T3}).
Hence, we can use the argument in the proof of the first assertion here too.
Now, suppose that $\xi$ satisfies one of the equalities in (\ref{equality-T3}).
We only consider the case that
$\xi$ satisfies $P_{x_k}+N_{x_k}=p_k$,
because the other cases can be treated similarly.
If $P_{x_k}=0$, then $N_{x_k}=p_k$,
and this implies that $\xi$ is not of Type 3.
Thus, in this case, we must have $P_{x_k} \geq 1$,
and similarly, $N_{x_k} \geq 1$.
Now, let $x_k^{a_N}$ and $x_k^{a_P}$ be sub-words of $\xi$
realizing the values $N_{x_k}$ and $P_{x_k}$, respectively.
Then, we have $a_N <0$ and $a_P >0$.
Also, using (\ref{alpha-beta}),
we obtain $\alpha_N=p_k+a_N=p_k-N_{x_k}$ and $\alpha_P=a_P=P_{x_k}$.
Thus, from $P_{x_k}+N_{x_k}=p_k$,
we obtain $\alpha_N=\alpha_P$.
Hence,
we can choose $x_k^{\alpha_P}$ 
in the suitable-spread procedure instead of $x_k^{\alpha_N}$.
Doing so,
we obtain another geodesic representative.
Hence,
$\xi$ is not the unique geodesic representative, 
and thus $\# {\rm SS}_g \geq 2$.
\noindent
$\Box$

%

\section{A unique geodesic representative}\label{section-ugr}

In this section, 
we propose a criterion that can be used to
uniquely specify a single proper geodesic representative for an element 
$g \in G(p_1,\dots,p_n)$,
even in the case that multiple geodesic representatives exist for 
that element. 
Using this criterion, we define a set ${\Gamma}$
such that the unique geodesic representative of
every element of $G(p_1,\dots,p_n)$ is an element of $\Gamma$.
This is done according to the decomposition of $G(p_1,\dots,p_n)$ given in (\ref{decomp-G}).

\subsection{Types {\boldmath $1$, $2$ \mbox{and} $3^+ \cup 3^-$}}\label{section-1,2,3^+,3^-}

\noindent
{\bf [Type {\boldmath $1$}]}
Consider an element $g \in {G(p_1,\dots,p_n)}_1$, 
and  let $\xi \in {\rm WT}_1$ be its geodesic representative.
By Proposition \ref{pr-T1,2,3^+,3^-},
the modified normal form $\nu$ of $g=\pi(\xi)$ is also geodesic.
We regard $\nu$ as a proper geodesic representative for this element $g$
and define
\begin{eqnarray}\label{def-gamma1}
\begin{array}{rl}
{\Gamma}_{1}:=
\{
\nu \in \Sigma^* \ 
;&
\nu \equiv 
x_{i_1}^{a_{1}} \cdots x_{i_\tau}^{a_{\tau}}
 \cdot \Delta^{c},\\
&\mbox{$\nu$ satisfies the conditions in (\ref{tau}) and (\ref{cond-T1-geod})}\}.\\
\end{array}
\end{eqnarray}

\noindent
{\bf [Type {\boldmath $2$}]}
Consider an element $g \in {G(p_1,\dots,p_n)}_2$.
The inverse of $g$, $g^{-1}$,
is contained in $G(p_1,\dots,p_n)_1$.
Let us consider the modified normal form $\nu$ of $g^{-1}$.
Then, we regard $\nu^{-1}$ as a proper geodesic representative 
for the element $g$
and define
\begin{eqnarray}\label{def-gamma2}
\begin{array}{rl}
{\Gamma}_{2}:=
\{
\nu^{-1} \in \Sigma^* \ 
;&
\nu \equiv 
x_{i_1}^{a_{1}} \cdots x_{i_\tau}^{a_{\tau}}
 \cdot \Delta^{c},\\
&\mbox{$\nu$ satisfies the conditions in (\ref{tau}) and
 (\ref{cond-T1-geod})}\}.
\end{array}
\end{eqnarray}

\noindent
{\bf [Type {\boldmath $3^+ \cup 3^-$}]}
Consider an element $g \in G(p_1,\dots,p_n)_{3^+ \cup 3^-}$,
and let $\xi \in {\rm WT}_{3^+ \cup 3^-}$ be
its geodesic representative.
Because $\xi$ satisfies (\ref{inequality-T3}),
then by Proposition \ref{pr-sc-T3}.3,
$\xi$ itself is the unique geodesic representative for $g=\pi(\xi)$.
Let us define 
\begin{eqnarray}\label{def-gamma3^+,3^-}
{\Gamma}_{3^+}:={\rm WT}_{3^+},\ \ 
{\Gamma}_{3^-}:={\rm WT}_{3^-}.
\end{eqnarray}
From Remark \ref{rmk-wt},
we have
\begin{eqnarray}\label{3^+and3^-}
\begin{array}{rcl}
 {\Gamma}_{3^+} \cap {\Gamma}_{3^-}
&=&{\rm WT}_{3^+} \cap {\rm WT}_{3^-}\\
&=&
\begin{array}{rcl}
\{\xi \in \Sigma^*
&;&  
\xi \equiv x_{i_1}^{a_{1}} \cdots x_{i_\tau}^{a_{\tau}},\
\xi \ \mbox{satisfies the condition in (\ref{tau})},\\
&&-\FLm \leq a_j \leq  \FLm \ (1 \leq j \leq \tau),\
a_j \neq 0 \ (1 \leq j \leq \tau)
\}.
\end{array}
\end{array}
\end{eqnarray}
We introduce the following notation:
\begin{eqnarray}\label{notation-gamma3^+,3^-}
{\Gamma}_{3^+ \cup 3^-}:=
{\Gamma}_{3^+} \cup {\Gamma}_{3^-},\ \ 
{\Gamma}_{3^+ \cap 3^-}:=
{\Gamma}_{3^+} \cap {\Gamma}_{3^-}.
\end{eqnarray}

\subsection{Type {\boldmath $3^0$}}

Consider an element $g \in G(p_1,\dots,p_n)_{3^0}$, and 
let $\xi \in {\rm WT}_{3^0}$ be its geodesic representative.
By Proposition \ref{pr-sc-T3}.4,
if $\xi$ satisfies one of the equalities in
(\ref{equality-T3}), 
there exist multiple geodesic representatives for
$g=\pi(\xi)$,
and moreover, by Proposition \ref{pr-sc-T3}.2,
these geodesic representatives are obtained from the suitable-spread procedure.
Here, we propose a method
that uniquely specifies a single proper output of the suitable-spread procedure
in such a case.

First, we rearrange $p_1,\cdots, p_n$ as follows:
(1) consider all $p_k$ such that $p_k-p_1$ are even;
(2) letting $m$ be the number of such $p_k$'s, 
rewrite them as $q_1,\dots,q_m$,
where the correspondence is chosen so as to satisfy
$q_1 \leq q_2 \leq \cdots \leq q_m$;
(3) rewrite the remaining as $r_{m+1},\dots, r_n$, 
again choosing the correspondence such that
$r_{m+1} \leq \cdots \leq r_n$.
Then, we have
\begin{eqnarray}\label{pqr}
\left\{
\begin{array}{rl}
\bullet&
q_1=p_1,\\
\bullet&
q_k-p_1\  \mbox{is even for each}\ k\in \{1,\dots,m\},\\
\bullet&
r_k-p_1\   \mbox{is odd for each}\ k\in \{m+1,\dots,n\},\\
\bullet&
q_1 \leq q_2 \leq \cdots \leq q_m,\
r_{m+1} \leq \cdots \leq r_n,\\
\bullet&
\{p_1,\dots,p_n\}=\{q_1,\dots,q_m\} \sqcup \{r_{m+1},\dots,r_n\}\
(\mbox{disjoint union}).
\end{array}
\right.
\end{eqnarray}

Next,
corresponding to the above rewriting of $p_1, \dots, p_n$,
we rename the generators, $x_1,\dots,x_n$, as
\begin{eqnarray}\label{xyz}
\begin{array}{rcl}
&\bullet &
\mbox{if $p_k$ is rewritten to $q_{k'}$, then $x_k$ is rename as $y_{k'}$},
\\
&\bullet &
\mbox{if $p_k$ is rewritten to $r_{k'}$, then $x_k$ is rename as $z_{k'}$}.
\end{array}
\end{eqnarray}
Then,
the group 
$G(p_1,p_2,\dots,p_n)$ is represented as
\begin{eqnarray}\label{new-presentation}
G(p_1,p_2,\dots,p_n)\ =\ \langle\ y_1,\dots, y_m,
z_{m+1}, \dots, z_n \ |\ 
y_1^{q_1} =\cdots=y_m^{q_m}=z_{m+1}^{r_{m+1}}=\cdots =z_n^{r_n}\ 
\rangle.
\end{eqnarray}

Also, we note the following:
\begin{eqnarray}\label{qr}
\begin{array}{rclll}
&\bullet&
q_k^-=\lfloor\frac{q_k-p_1}2\rfloor
=\frac{q_k-p_1}2,
&q_k^+=\lfloor\frac{p_1+q_k-1}2\rfloor
=\frac{p_1+q_k-2}2
&(k \in \{1,\dots,m\}),\\
&\bullet&
r_k^-=\lfloor\frac{r_k-p_1}2\rfloor
=\frac{r_k-p_1-1}2,
&r_k^+=\lfloor\frac{p_1+r_k-1}2\rfloor
=\frac{p_1+r_k-1}2
&(k \in \{m+1,\dots,n\}).\\
\end{array}
\end{eqnarray}

Now, by the argument given in the proof of Proposition \ref{pr-ss-1},
the modified normal form of $g=\pi(\xi)$ takes the form
\begin{eqnarray*}
\nu \equiv
x_{i_1}^{\overline{\alpha}_{1}}\cdots x_{i_\tau}^{\overline{\alpha}_{\tau}}
\cdot (\Delta^{-1})^{\delta},
\end{eqnarray*}
where
$\tau$ and $i_j$ satisfy the condition in (\ref{tau}),
${\overline{\alpha}}_j$ satisfies the condition in (\ref{cond-mnf}),
and $\delta$ satisfies the relations
$\delta=-(d+\rho)$
and 
${\bf r}_\nu > \delta >0$.
Next, for $k \in \{1,\dots,m\}$ and
$l \in \{ 0,1,\dots,p_1-2 \}$,
we define
\begin{eqnarray}\left\{
\begin{array}{rcl}
Y_k[l]&:=&\{x_{i_j}^{{\overline{\alpha}}_j}  \ ; \ x_{i_j}^{{\overline{\alpha}}_j} 
\ \mbox{is a sub-word of} \ \nu,\
x_{i_j}=y_k,\ 
\ {\overline{\alpha}}_j=q_k^+-l
\},\\
y_k[l]&:=&\# Y_k[l],\\
Y[l]&:=&\sqcup_{k=1}^{m} Y_k[l]\quad \mbox{(disjoint union)},\\
y[l]&:=&\# Y[l]=\sum_{k=1}^m y_k[l],\\
\end{array}
\right.
\end{eqnarray}
and  for $k \in \{m+1,\dots,n\}$ and
$l \in \{ 0,1,\dots,p_1-1 \}$,
 we define
\begin{eqnarray}\left\{
\begin{array}{rcl}
Z_k[l]&:=&\{x_{i_j}^{{\overline{\alpha}}_j}  \ ; \ x_{i_j}^{{\overline{\alpha}}_j} 
\ \mbox{is a sub-word of} \ \nu,\ 
x_{i_j}=z_k,\ 
{\overline{\alpha}}_j=r_k^+-l
\},\\
z_k[l]&:=&\# Z_k[l],\\
Z[l]&:=&\sqcup_{k=m+1}^{n} Z_k[l]\quad \mbox{(disjoint union)},\\
z[l]&:=&\# Z[l]=\sum_{k=m+1}^n z_k[l].
\end{array}
\right.
\end{eqnarray}
Then, we have
\begin{eqnarray}
&&
\begin{array}{rcl}\label{y-1}
&\bullet&y_{k}^{{\overline{\alpha}}_{j}} \in Y[l]\quad
\Rightarrow\quad
2 {\overline{\alpha}}_{j}+p_1-q_k =2 p_1-2-2 l,
\end{array}\\
&&
\begin{array}{rcl}\label{z-1}
&\bullet&z_k^{{\overline{\alpha}}_{j}} \in Z[l]\quad 
\Rightarrow\quad
2 {\overline{\alpha}}_{j}+p_1-r_k =2 p_1-1-2 l,
\end{array}
\end{eqnarray}
and hence also
\begin{eqnarray}
\begin{array}{rccl}\label{xy-1}
\bullet&z_k^{{\overline{\alpha}}_{j}} \in Z[l],\
z_{k'}^{{\overline{\alpha}}_{j'}} \in Z[l'] \quad
\mbox{and}\quad
l \leq l' \
&\Rightarrow&
2 {\overline{\alpha}}_{j}+p_1-r_k \geq 2 {\overline{\alpha}}_{j'}+p_1-r_{k'},
\end{array}
\\
\begin{array}{rccl}\label{xy-2}
\bullet&z_k^{{\overline{\alpha}}_{j}} \in Z[l],\
y_{k'}^{{\overline{\alpha}}_{j'}} \in Y[l'] \quad
\mbox{and}\quad
l \leq l' \
&\Rightarrow&
2 {\overline{\alpha}}_{j}+p_1-r_k > 2 {\overline{\alpha}}_{j'}+p_1-q_{k'},
\end{array}
\\
\begin{array}{rccl}\label{xy-3}
\bullet&
y_k^{{\overline{\alpha}}_{j}} \in Y[l],\
y_{k'}^{{\overline{\alpha}}_{j'}} \in Y[l'] \quad
\mbox{and}\quad
l \leq l' \
&\Rightarrow&
2 {\overline{\alpha}}_{j}+p_1-q_k \geq 2 {\overline{\alpha}}_{j'}+p_1-q_{k'},
\end{array}
\\
\begin{array}{rccl}\label{xy-4}
\bullet&
y_k^{{\overline{\alpha}}_{j}} \in Y[l],\
z_{k'}^{{\overline{\alpha}}_{j'}} \in Z[l'] \quad
\mbox{and}\quad
l < l' \
&\Rightarrow&
2 {\overline{\alpha}}_{j}+p_1-q_k> 2 {\overline{\alpha}}_{j'}+p_1-r_{k'}.
\end{array}
\end{eqnarray}
Moreover, we have
\begin{eqnarray}\label{xy-5}
{\bf r}_\nu=
\sum_{l=0}^{p_1-1} z[l]+\sum_{l=0}^{p_1-2} y[l],\label{r-1}
\end{eqnarray}
which follows from (\ref{pq-oo}) and (\ref{pq-oe}),
along with the definition of ${\bf r}_\nu$, given in (\ref{def-R}).



Note that from (\ref{y-1}),
in the suitable-spread procedure,
we can choose elements from $Y[l](=\sqcup_{k=1}^{m} Y_k[l])$
without considering which subsets they belong to.
Also, from (\ref{z-1}),
the same holds for $Z[l](=\sqcup_{k=m+1}^{n}Z_k[l])$.
With this observation, we adopt the following rule:
\begin{eqnarray}
\begin{array}{ll}\label{rule-y}
\bullet&
\mbox{We choose elements from 
$Y[l]=\sqcup_{k=1}^{m} Y_k[l]$
in descending order. In other words, }\\
&\mbox{we choose elements from $Y_k[l]$ before those from $Y_{k-1}[l]$
$(k=m, m-1, \dots, 2)$.}\\
\end{array}
\\
\begin{array}{ll}\label{rule-z}
\bullet&
\mbox{We choose elements from $Z[l]=\sqcup_{k=m+1}^{n}Z_k[l]$
in descending order. In other words,}\\
&\mbox{we choose elements from $Z_k[l]$ before those from $Z_{k-1}[l]$
$(k=n, n-1, \dots, m+2)$.}\\
\end{array}
\end{eqnarray}

\subsubsection{Case: {\boldmath $p_k-p_1$} is odd for some $k \in \{ 1, \dots, n\}$}\label{section:even and odd}

In this section, we consider the case $m<n$, i.e.,
the case in which
there exists $k \in \{1, \dots,n\}$ such that $p_k-p_1$ is odd.
Then, from (\ref{xy-5}), we have one of the following cases for $\nu$.
(These cases are separated according to where
$\delta$ is located in the open interval $(0,{\bf r}_\nu)$,
as seen in Figure 1.)
\vspace{0.3cm}

\noindent
$\bullet$ Case  $(N,2M+1)$, in which $0 \leq N \leq p_1-1$ and
$0 \leq  M \leq n-m-1$:
\begin{eqnarray*}
\begin{array}{rcl}
{\displaystyle \sum_{l=0}^{N-1} \left(z[l]+y[l]\right)}
&+&{\displaystyle \sum_{k=n-(M-1)}^{n} z_k[N]
\quad < \quad \delta} \\
&<& {\displaystyle \sum_{l=0}^{N-1}\left(z[l]+y[l]\right)
+\sum_{k=n-M}^{n} z_k[N]};
\end{array}
\end{eqnarray*}

\noindent
$\bullet$ Case $(N,2M+2)$, in which $0 \leq N \leq p_1-2$ and $0 \leq  M \leq n-m-1$,
or $N=p_1-1$ and $0 \leq M \leq n-m-2$:
\begin{eqnarray*}
{\displaystyle \delta = \sum_{l=0}^{N-1}\left(z[l]+y[l]\right)
+\sum_{k=n-M}^{n} z_k[N]}\quad
\mbox{with}\quad  z_{n-M}[N] >0;
\end{eqnarray*}

\noindent
$\bullet$ Case  $(N,2M+1)$, in which $0 \leq N \leq p_1-2$ and 
$n-m \leq M \leq n-1$:
\begin{eqnarray*}
\begin{array}{rcl}
{\displaystyle \sum_{l=0}^{N-1} \left(z[l]+y[l]\right)
+z[N]}
&+&{\displaystyle \sum_{k=n-(M-1)}^{m} y_k[N]}
\quad < \quad \delta \\
&<& {\displaystyle \sum_{l=0}^{N-1}\left(z[l]+y[l]\right)
+z[N]
+\sum_{k=n-M}^{m} y_k[N]};
\end{array}
\end{eqnarray*}

\noindent
$\bullet$ Case  $(N,2M+2)$, in which $0 \leq N \leq p_1-2$ and $n-m \leq M \leq n-1$:
 \begin{eqnarray*}
{\displaystyle \delta = \sum_{l=0}^{N-1}\left(z[l]+y[l]\right)
+z[N]
+\sum_{k=n-M}^{m} y_k[N]}\quad
\mbox{with} \quad  y_{n-M}[N]>0.
\end{eqnarray*}

\vspace{0.5cm}

\setlength{\unitlength}{0.7mm}
\begin{picture}(180,300)(-30,20)

\put(-4,282){\makebox(50,10){$z_n[0]$}}
\put(45,282){\makebox(50,10){$z_n[0]+z_{n-1}[0]$}}
\put(105,282){\makebox(50,10){$z_n[0]+z_{n-1}[0]+z_{n-2}[0]$}}
\put(-32,282){\makebox(10,10){$0$}}
\put(-30,279){\line(1,0){170}}
\put(148,278.2){$\dots \dots$}
\put(15,276.5){\makebox(5,5){$\bullet$}}
\put(66,276.5){\makebox(5,5){$\bullet$}}
\put(-30,276.5){\makebox(5,5){$\bullet$}}
\put(128,276.5){\makebox(5,5){$\bullet$}}
\put(-10,264){\makebox(10,10){\mbox{Case}\ $(0,1)$}}
\put(17,257){\makebox(10,10){\mbox{Case}\ $(0,2)$}}
\put(40,264){\makebox(10,10){\mbox{Case}\ $(0,3)$}}
\put(67,257){\makebox(10,10){\mbox{Case}\ $(0,4)$}}
\put(100,264){\makebox(10,10){\mbox{Case}\ $(0,5)$}}
\put(132,257){\makebox(10,10){\mbox{Case}\ $(0,6)$}}
\put(-8,273){\vector(0,1){5}}
\put(17.4,266){\vector(0,1){10}}
\put(42,273){\vector(0,1){5}}
\put(68.7,266){\vector(0,1){10}}
\put(102,273){\vector(0,1){5}}
\put(130.5,266){\vector(0,1){10}}
\put(-21,240){\makebox(50,10){$z[0]$}}
\put(50,240){\makebox(50,10){$z[0]+y_m[0]$}}
\put(135,240){\makebox(50,10){$z[0]+y_m[0]+y_{m-1}[0]$}}
\put(0,237){\line(1,0){165}}
\put(-30,236.2){$\dots \dots$}
\put(175,236.2){$\dots \dots$}
\put(0,234.5){\makebox(5,5){$\bullet$}}
\put(71,234.5){\makebox(5,5){$\bullet$}}
\put(153,234.5){\makebox(5,5){$\bullet$}}
\put(-5,214){\makebox(10,10){\mbox{Case}\ $(0,2n-2m)$}}
\put(35,222){\makebox(10,10){\mbox{Case}\ $(0,2n-2m+1)$}}
\put(65,214){\makebox(10,10){\mbox{Case}\ $(0,2n-2m+2)$}}
\put(110,222){\makebox(10,10){\mbox{Case}\ $(0,2n-2m+3)$}}
\put(150,214){\makebox(10,10){\mbox{Case}\ $(0,2n-2m+4)$}}
\put(2.4,223){\vector(0,1){10}}
\put(35,230){\vector(0,1){5}}
\put(73.5,223){\vector(0,1){10}}
\put(118,230){\vector(0,1){5}}
\put(155.5,223){\vector(0,1){10}}
\put(-20,198){\makebox(50,10){$\sum_{l=0}^{N-1} (z[l]+y[l])$}}
\put(40,198){\makebox(50,10){$\sum_{l=0}^{N-1} (z[l]+y[l])+z_n[N]$}}
\put(125,198){\makebox(50,10){$\sum_{l=0}^{N-1} (z[l]+y[l])+z_n[N]+z_{n-1}[N]$}}
\put(0,195){\line(1,0){165}}
\put(-30,194.2){$\dots \dots$}
\put(175,194.2){$\dots \dots$}
\put(5,192.5){\makebox(5,5){$\bullet$}}
\put(71,192.5){\makebox(5,5){$\bullet$}}
\put(153,192.5){\makebox(5,5){$\bullet$}}
\put(0,172){\makebox(10,10){\mbox{Case}\ $(N-1,2n)$}}
\put(35,180){\makebox(10,10){\mbox{Case}\ $(N,1)$}}
\put(65,172){\makebox(10,10){\mbox{Case}\ $(N,2)$}}
\put(105,180){\makebox(10,10){\mbox{Case}\ $(N,3)$}}
\put(150,172){\makebox(10,10){\mbox{Case}\ $(N,4)$}}
\put(7.4,181){\vector(0,1){10}}
\put(38,188){\vector(0,1){5}}
\put(73.5,181){\vector(0,1){10}}
\put(108,188){\vector(0,1){5}}
\put(155.5,181){\vector(0,1){10}}
\put(5,170){\makebox(50,10){}}
\put(-23,153){\makebox(50,10){$\sum_{l=0}^{N} z[l]+\sum_{l=0}^{N-1}y[l]$}}
\put(55,153){\makebox(50,10){$\sum_{l=0}^{N} z[l]+\sum_{l=0}^{N-1}y[l]+y_m[N]$}}
\put(125,163){\makebox(50,10){$\sum_{l=0}^{N} z[l]+\sum_{l=0}^{N-1}y[l]+y_m[N]+y_{m-1}[N]$}}
\put(160,163){\vector(0,-1){10}}
\put(0,150){\line(1,0){165}}
\put(-30,149.2){$\dots \dots$}
\put(175,149.2){$\dots \dots$}
\put(5,147.5){\makebox(5,5){$\bullet$}}
\put(86,147.5){\makebox(5,5){$\bullet$}}
\put(158,147.5){\makebox(5,5){$\bullet$}}
\put(3,127){\makebox(10,10){\mbox{Case}\ $(N,2n-2m)$}}
\put(43,135){\makebox(10,10){\mbox{Case}\ $(N,2n-2m+1)$}}
\put(80,127){\makebox(10,10){\mbox{Case}\ $(N,2n-2m+2)$}}
\put(120,135){\makebox(10,10){\mbox{Case}\ $(N,2n-2m+3)$}}
\put(161,127){\makebox(10,10){\mbox{Case}\ $(N,2n-2m+4)$}}
\put(7.4,136){\vector(0,1){10}}
\put(46,143){\vector(0,1){5}}
\put(88.5,136){\vector(0,1){10}}
\put(123,143){\vector(0,1){5}}
\put(161.5,136){\vector(0,1){10}}
\put(5,123){\makebox(50,10){}}
\put(-20,103){\makebox(50,10){$\sum_{l=0}^{N} (z[l]+y[l])$}}
\put(41,113){\makebox(50,10){$\sum_{l=0}^{N} (z[l]+y[l])+z_n[N+1]$}}
\put(120,103){\makebox(50,10){$\sum_{l=0}^{N} (z[l]+y[l])+z_n[N+1]+z_{n-1}[N+1]$}}
\put(0,100){\line(1,0){165}}
\put(-30,99.2){$\dots \dots$}
\put(175,99.2){$\dots \dots$}
\put(5,97.5){\makebox(5,5){$\bullet$}}
\put(65,97.5){\makebox(5,5){$\bullet$}}
\put(130,97.5){\makebox(5,5){$\bullet$}}
\put(3,77){\makebox(10,10){\mbox{Case}\ $(N,2n)$}}
\put(35,85){\makebox(10,10){\mbox{Case}\ $(N+1,1)$}}
\put(65,77){\makebox(10,10){\mbox{Case}\ $(N+1,2)$}}
\put(88,85){\makebox(10,10){\mbox{Case}\ $(N+1,3)$}}
\put(125,77){\makebox(10,10){\mbox{Case}\ $(N+1,4)$}}
\put(7.4,86){\vector(0,1){10}}
\put(38,93){\vector(0,1){5}}
\put(67.5,113){\vector(0,-1){10}}
\put(67.5,86){\vector(0,1){10}}
\put(96.5,93){\vector(0,1){5}}
\put(133.5,86){\vector(0,1){10}}
%
\put(5,73){\makebox(50,10){}}
\put(0,53){\makebox(50,10){$\sum_{l=0}^{p_1-2} (z[l]+y[l])+z_n[p_1-1]+\cdots+z_{m+2}[p_1-1]$}}
\put(146,63){\makebox(50,10){$\sum_{l=0}^{p_1-1}z[l]+\sum_{l=0}^{p_1-2}y[l]$}}
\put(148,49){\makebox(50,10){${\bf r}_\nu$}}
\put(171,58){\line(0,1){5}}
\put(173,58){\line(0,1){5}}
\put(0,50){\line(1,0){180}}
\put(-30,49.2){$\dots \dots$}
\put(5,47.5){\makebox(5,5){$\bullet$}}
\put(170,47.5){\makebox(5,5){$\bullet$}}
\put(3,27){\makebox(10,10){\mbox{Case}\ $(p_1-1,2n-2m-2)$}}
\put(88,35){\makebox(10,10){\mbox{Case}\ $(p_1-1,2n-2m-1)$}}
\put(7.4,36){\vector(0,1){10}}
\put(96.5,43){\vector(0,1){5}}

\end{picture}
\begin{center}
Figure 1.
\end{center}

\vspace{0.2cm}

\begin{ex}\label{example-unique-geodesic}
{\rm
Let $n=3$, with $p_1=3, p_2=6$ and $p_3=7$.
Then, we consider the same group as in Example \ref{example-modified-normal-form},
$G(3,6,7)\ =\ \langle\ x_1,x_2, x_3 \ |\ x_1^{3} =x_2^{6}=x_3^{7}\ \rangle$.
In this case, $p_2-p_1$ is odd and $p_3-p_1$ is even.
Thus, we have
$m=2$, 
$q_1=3$, $q_2=7$ and $r_3=6$, and $G(3,6,7)$ is represented as
$G(3,6,7)\ =\ \langle\ y_1,y_2, z_3 \ |\ y_1^{3} =y_2^{7}=z_3^{6}\ \rangle$.
Also, we have
\begin{eqnarray*}
q_1^-=0, \ q_1^+=2; \quad
q_2^-=2, \ q_2^+=4; \quad
r_3^-=1, \ r_3^+=4.
\end{eqnarray*}
Now, let us consider an element $h \in G(3,6,7)$
whose modified normal form is given by
\begin{eqnarray*}
\begin{array}{rcl}
\nu &\equiv&
z_3^2 \cdot y_2^4  \cdot y_1 \cdot y_2^4
 \cdot z_3^{-1} \cdot y_2^{-2} 
\cdot z_3^4  \cdot y_2^2 \cdot y_1^2 \cdot z_3^4 \cdot
y_2^4 \cdot z_3^2 \cdot (\Delta^{-1})^\delta\\
\\
&\equiv&
z_3^{\overline{\alpha}_1} \cdot  y_2^{\overline{\alpha}_2} 
 \cdot y_1^{\overline{\alpha}_3} 
 \cdot y_2^{\overline{\alpha}_4} \cdot z_3^{\overline{\alpha}_5} 
\cdot y_2^{\overline{\alpha}_6}  \cdot z_3^{\overline{\alpha}_{7}} 
\cdot y_2^{\overline{\alpha}_{8}} 
\cdot y_1^{\overline{\alpha}_{9}} \cdot z_3^{\overline{\alpha}_{10}} 
\cdot y_2^{\overline{\alpha}_{11}} 
\cdot z_3^{\overline{\alpha}_{12}} \cdot (\Delta^{-1})^\delta.
\end{array}
\end{eqnarray*}
Then, we have
\begin{eqnarray*}
\left\{
\begin{array}{lll}
Y_1[0]=\{y_1^{\overline{\alpha}_{9}}\},\
&Y_1[1]=\{y_1^{\overline{\alpha}_{3}}\},
&\\
Y_2[0]=
\{y_2^{\overline{\alpha}_{2}}, y_2^{\overline{\alpha}_{4}}, y_2^{\overline{\alpha}_{11}}\},\
&Y_2[1]=\emptyset,
&\\
Z_3[0]=
\{z_3^{\overline{\alpha}_{7}}, z_3^{\overline{\alpha}_{10}}\},\
&Z_3[1]=\emptyset,\
&Z_3[2]=\{z_3^{\overline{\alpha}_{1}},z_3^{\overline{\alpha}_{12}}\},\\
\end{array}
\right.
\end{eqnarray*}
\begin{eqnarray*}
R_\nu=\{1,2,3,4,7,9,10,11,12\}
\end{eqnarray*}
and
\begin{eqnarray*}
\left\{
\begin{array}{l}
z[0]=z_3[0]=2,\\
z[0]+y_2[0]=z_3[0]+y_2[0]=5,\\
z[0]+y[0]=z_3[0]+(y_2[0]+y_1[0])=6,\\
(z[0]+y[0])+z[1]=(z[0]+y[0])+z_3[1]=6,\\
(z[0]+y[0])+z[1]+y_2[1]=(z[0]+y[0])+z_3[1]+y_2[1]=6,\\
(z[0]+y[0])+(z[1]+y[1])=(z[0]+y[0])+(z_3[1]+(y_2[1]+y_1[1]))=7,\\
{\bf r}_\nu=
\sum_{l=0}^{p_1-1} z[l]+\sum_{l=0}^{p_1-2} y[l]
=(z[0]+y[0])+(z[1]+y[1])+z_3[2]=9.
\end{array}
\right.
\end{eqnarray*}
Hence, if $0 < \delta  < 9$,
then $h$ is of Type $3^0$.
\\

\noindent
For this example, we have the following:
\\

\noindent
$\bullet$ If $\delta=1$,
then $0 < \delta < z_3[0]$, and hence we have
Case $(0,1)$.

\noindent
$\bullet$ If $\delta=2$,
then $\delta =z_3[0]$ and $z_3[0] >0$,
and hence we have Case $(0,2)$.

\noindent
$\bullet$ If $\delta=3\ \mbox{or}\ 4$,
then $z[0] < \delta < z[0]+y_2[0]$, and hence we have
Case $(0,3)$.

\noindent
$\bullet$ If $\delta=5$,
then $\delta = z[0]+y_2[0]$ and $y_2[0] >0$,
and hence we have Case $(0,4)$.

\noindent
$\bullet$ If $\delta=6$,
then $\delta =z[0]+y_2[0]+y_1[0]$ and $y_1[0] >0$,
and hence we have Case $(0,6)$.

\noindent
$\bullet$ If $\delta=7$,
then $ \delta =(z[0]+y[0])+z[1]+y_2[1]+y_1[1] $ and $y_1[1] >0$,
and hence we have Case $(1,6)$.

\noindent
$\bullet$ If $\delta=8$,
then $ (z[0]+y[0])+(z[1]+y[1]) < \delta =(z[0]+y[0])+(z[1]+y[1])+z_3[2]$,
and hence we have Case $(2,1)$.
$\Box$
\\


}
\end{ex}

Next,
we describe the proposed method yielding a unique proper output $\widehat{\nu}$
of the suitable-spread procedure.
First, we present an explicit expression for $\widehat{\nu}$, and then
we define the set consisting of all $\widehat{\nu}$.
This is carried out separately in each of the cases considered in Figure 1.

\vspace{0.2cm}

\noindent
{\bf [Case {\boldmath $(N,2M+1)$}, in which {\boldmath $0 \leq N \leq p_1-1$} and {\boldmath $0 \leq M \leq n-m-1$}]}

\noindent
In this case,
we perform the suitable-spread procedure 
in accordance with the following rules.
\begin{enumerate}
\item
All elements of $\left(\sqcup_{l=0}^{N-1} (Z[l] \sqcup Y[l])\right)
\sqcup (\sqcup_{k=n-(M-1)}^n Z_k[N])$
are chosen,
because we have the relations in (\ref{xy-1}), (\ref{xy-2}) and  (\ref{xy-4}),
and we adopt the rule (\ref{rule-z}).
\item
Define
$\phi:=\delta-\left(\sum_{l=0}^{N-1} (z[l]+y[l])
+\sum_{k=n-(M-1)}^n z_k[N]\right)$
(see Figure 2).
In this case, we have $0 < \phi < z_{n-M}[N]$.
Hence, from (\ref{xy-1}), (\ref{xy-2}) and (\ref{rule-z}), 
we choose $\phi$ sub-words from $Z_{n-M}[N]$.
Here, for conciseness, we write $r_{n-M}$,
$r^+_{n-M}$ and $r^-_{n-M}$ as $r$, $r^+$ and $r^-$.
\item
When choosing the $\phi$ sub-words from $Z_{n-M}[N]$,
we always choose the leftmost of these.
We represent all the sub-words belonging to $Z_{n-M}[N]$
by 
$z_{n-M}^{{\overline{\alpha}}_{j_1}}, \dots, 
z_{n-M}^{{\overline{\alpha}}_{j_{z_{n-M}[N]}}}$.
Then, the sub-words chosen through this procedure are
$z_{n-M}^{{\overline{\alpha}}_{j_1}}, \dots, z_{n-M}^{{\overline{\alpha}}_{j_\phi}}$.

\end{enumerate}

\setlength{\unitlength}{0.7mm}
\begin{picture}(200,40)(-30,170)

\put(-4,193){\makebox(50,10){$\sum_{l=0}^{N-1} (z[l]+y[l])
+\sum_{k=n-(M-1)}^n z_k[N]$}}
\put(63,193){\makebox(50,10){$\delta$}}
\put(130,193){\makebox(50,10){$\sum_{l=0}^{N-1} (z[l]+y[l])
+\sum_{k=n-M}^n z_k[N]$}}
\put(0,190){\line(1,0){165}}

\put(15,187.5){\makebox(5,5){$\bullet$}}
\put(86,187.5){\makebox(5,5){$\bullet$}}
\put(153,187.5){\makebox(5,5){$\bullet$}}
\put(48,175){\makebox(10,10){$\phi$}}
\put(90,167){\makebox(10,10){$z_{n-M}[N]$}}
\put(17.4,170){\line(0,1){19}}
\put(42.4,180){\vector(-1,0){25}}
\put(63,180){\vector(1,0){25}}
\put(88.5,176){\line(0,1){14}}
\put(78,173){\vector(-1,0){60}}
\put(110,173){\vector(1,0){45}}
\put(155.5,170){\line(0,1){19}}
%

\end{picture}
\begin{center}
Figure 2.
The location of $\delta$ and the definition of $\phi$.
\end{center}

\noindent
In this way,
we obtain an output $\widehat{\nu}$ containing 
$z_{n-M}^{{\overline{\alpha}}_{j_1}-r}, \dots, 
z_{n-M}^{{\overline{\alpha}}_{j_\phi}-r}$
and $z_{n-M}^{{\overline{\alpha}}_{j_{\phi+1}}},\dots, 
z_{n-M}^{{\overline{\alpha}}_{j_{z_{n-M}[N]}}}$.
Here, note that
$z_{n-M}^{{\overline{\alpha}}_{j_t}-r} \equiv 
z_{n-M}^{-(r^-+N+1)}$ ($t \in \{1,\dots,\phi\}$),
because $(r^+-N)-r=-(r^-+N+1)$.
Then, choosing $z_-$ (resp., $z_+$) to be the leftmost element,
$z_{n-M}^{{\overline{\alpha}}_{j_1}-r}$ (resp., $z_{n-M}^{{\overline{\alpha}}_{j_{\phi+1}}}$),
we can write the output $\widehat{\nu}$ as 
\begin{eqnarray}\label{1-1}
\widehat{\nu}
\equiv
\xi^{(1)} \cdot z_-
\cdot \xi^{(2)} \cdot z_+
\cdot \xi^{(3)},
\end{eqnarray}
where $\xi^{(1)}$, $\xi^{(2)}$ and $\xi^{(3)}$ are given by
\begin{eqnarray}\label{1-2}
\left\{
\begin{array}{l}
\xi^{(1)} :=
x_{i_1}^{a_1} x_{i_2}^{a_2} \cdots
 x_{i_{\kappa_1}}^{a_{\kappa_1}},\\
\xi^{(2)} :=
x_{i_1}^{b_1} x_{i_2}^{b_2} \cdots
x_{i_{\kappa_2}}^{b_{\kappa_2}},\\
\xi^{(3)} :=
x_{i_1}^{c_1} x_{i_2}^{c_2} \cdots
x_{i_{\kappa_3}}^{c_{\kappa_3}},\\
\end{array}
\right.
\end{eqnarray}
with the following conditions:

\begin{eqnarray}
&\xi^{(1)}&
\left\{
\begin{array}{rl}\label{1-3}
\bullet& \kappa_1\in {\bf N} \cup \{0\},\\
\bullet& x_{i_{\kappa_1}}\neq z_{n-M},\\
\bullet&
x_{i_j}=z_k\ \mbox{with}\ n-M+1 \leq k \leq n\\
&\Longrightarrow
\begin{array}{l}
-(r_k^- +N+1) \leq a_j \leq  
r_k^+ -N-1,\  
a_j \neq 0,\\
\end{array}
\\
\bullet&
x_{i_j}=z_k\ \mbox{with}\ k=n-M\\
&\Longrightarrow 
\begin{array}{l}
-(r^- +N) \leq a_j \leq  
r^+ -N-1,\ 
a_j \neq 0,\\
\end{array}
\\
\bullet&
x_{i_j}=z_k\ \mbox{with}\ m+1 \leq k \leq n-M-1\\
&\Longrightarrow 
\begin{array}{l}
-(r_k^- +N) \leq a_j \leq  
r_k^+ -N,\
a_j \neq 0,\\
\end{array}
\\
\bullet&
x_{i_j}=y_k\ \mbox{with}\ 1 \leq k \leq m\\
&\Longrightarrow 
\begin{array}{l}
-(q_k^- +N) \leq a_j \leq  
q_k^+ -N,\
a_j \neq 0,\\
\end{array}
\\
\end{array}
\right.\\
&\xi^{(2)}&
\left\{
\begin{array}{rl}\label{1-4}
\bullet& \kappa_2\in {\bf N},\\
\bullet& x_{i_1} \neq z_{n-M}, x_{i_{\kappa_2}}\neq z_{n-M},\\
\bullet&
x_{i_j}=z_k\ \mbox{with}\ n-M+1 \leq k \leq n\\
&\Longrightarrow 
\begin{array}{l}
-(r_k^- +N+1) \leq b_j \leq  
r_k^+ -N-1,\ 
b_j \neq 0,\\
\end{array}
\\
\bullet&
x_{i_j}=z_k\ \mbox{with}\ k=n-M\\
&\Longrightarrow 
\begin{array}{l}
-(r^- +N+1) \leq b_j \leq  
r^+ -N-1,\  
b_j \neq 0,\\
\end{array}
\\
\bullet&
x_{i_j}=z_k\ \mbox{with}\ m+1 \leq k \leq n-M-1\\
&\Longrightarrow 
\begin{array}{l}
-(r_k^- +N) \leq b_j \leq  
r_k^+ -N,\ 
b_j \neq 0,\\
\end{array}
\\
\bullet&
x_{i_j}=y_k\ \mbox{with}\ 1 \leq k \leq m\\
&\Longrightarrow 
\begin{array}{l}
-(q_k^- +N) \leq b_j \leq  
q_k^+ -N,\
b_j \neq 0,\\
\end{array}
\\
\end{array}
\right.\\
&\xi^{(3)}&
\left\{
\begin{array}{rl}\label{1-5}
\bullet& \kappa_3 \in {\bf N} \cup \{0\},\\
\bullet& x_{i_1} \neq z_{n-M},\\
\bullet&
x_{i_j}=z_k\ \mbox{with}\ n-M+1 \leq k \leq n\\
&\Longrightarrow 
\begin{array}{l}
-(r_k^- +N+1) \leq c_j \leq  
r_k^+ -N-1,\ 
c_j \neq 0,\\
\end{array}
\\
\bullet&
x_{i_j}=z_k\ \mbox{with}\ k=n-M\\
&\Longrightarrow 
\begin{array}{l}
-(r^- +N) \leq c_j \leq  
r^+ -N,\ 
c_j \neq 0,\\
\end{array}
\\
\bullet&
x_{i_j}=z_k\ \mbox{with}\ m+1 \leq k \leq n-M-1\\
&\Longrightarrow 
\begin{array}{l}
-(r_k^- +N) \leq c_j \leq  
r_k^+ -N,\ 
c_j \neq 0,\\
\end{array}
\\
\bullet&
x_{i_j}=y_k\ \mbox{with}\ 1 \leq k \leq m\\
&\Longrightarrow 
\begin{array}{l}
-(q_k^- +N) \leq c_j \leq  
q_k^+ -N,\
c_j \neq 0.\\
\end{array}
\\
\end{array}
\right.
\end{eqnarray}
Then, for each pair $(N,M)$ 
satisfying $0 \leq N \leq p_1-1$ and $0 \leq M \leq n-m-1$,
we define the following set:
\begin{eqnarray*}
\begin{array}{rl}
{\Gamma}^{(N,2M+1)}_{3^0}:=
\{
\widehat{\nu} \in \Sigma^* \ ;
&
\widehat{\nu} 
\equiv
\xi^{(1)} \cdot z_{n-M}^{-(r^-+N+1)}
\cdot \xi^{(2)} \cdot z_{n-M}^{r^+ -N}
\cdot \xi^{(3)},\\
&\xi^{(1)}, \xi^{(2)} \ \mbox{and}\ \xi^{(3)} 
\ \mbox{are given in (\ref{1-2})} 
\ \mbox{with the conditions in (\ref{1-3})--(\ref{1-5})} \}.
\end{array}
\end{eqnarray*}

\vspace{0.2cm}

\begin{ex}
{\rm

We consider here the element $\nu$ given in Example \ref{example-unique-geodesic}
with $\delta=1$.
Then, $0 < \delta < 2=z_3[0]$, and we have
Case $(N,2M+1)$ with $N=0$ and $M=0$.
From (\ref{xy-1}) and (\ref{xy-2}),
we choose one element from $Z[0](=Z_3[0])$.
There are two choices, because $z_3[0]=2$.
Then, applying the suitable spread procedure in accordance with these choices,
we obtain
\begin{eqnarray*}
\begin{array}{rcl}
\nu&=&
z_3^2 \cdot y_2^4 \cdot y_1 \cdot y_2^4
 \cdot z_3^{-1} \cdot y_2^{-2} \cdot (z_3^4\Delta^{-1})
 \cdot y_2^2 \cdot y_1^2 \cdot z_3^4 \cdot
y_2^4 \cdot z_3^2
\\
&&=
z_3^2 \cdot y_2^4 \cdot y_1  \cdot y_2^4
 \cdot z_3^{-1} \cdot y_2^{-2} \cdot z_3^{-2}
\cdot y_2^2 \cdot y_1^2 \cdot z_3^4 \cdot
y_2^4 \cdot z_3^2
=:\widetilde{\nu}_1
\\
&=&
z_3^2 \cdot y_2^4 \cdot y_1 \cdot y_2^4
 \cdot z_3^{-1} \cdot y_2^{-2} \cdot z_3^4 
\cdot y_2^2 \cdot y_1^2 \cdot (z_3^4\Delta^{-1}) \cdot
y_2^4 \cdot z_3^2
\\
&&=
z_3^2 \cdot y_2^4 \cdot y_1 \cdot y_2^4
 \cdot z_3^{-1} \cdot y_2^{-2} \cdot z_3^4 
\cdot y_2^2 \cdot y_1^2 \cdot  z_3^{-2} \cdot
y_2^4 \cdot z_3^2
=:\widetilde{\nu}_2.
\\
\end{array}
\end{eqnarray*}
In this way, from Proposition \ref{pr-sc-T3},
 we obtain all the geodesic representatives of $\pi(\nu)$,
which are denoted as $\widetilde{\nu}_1$ and $\widetilde{\nu}_2$.
If we choose the leftmost element from $Z_3[0]$,
then we obtain
the word, $\widetilde{\nu}_1$, 
which is the output $\widehat{\nu}$ that we choose here.
\begin{eqnarray*}
\begin{array}{rcl}
\widehat{\nu} &\equiv&
\xi^{(1)} \cdot z_-
\cdot \xi^{(2)} \cdot z_+
\cdot \xi^{(3)}\\
&\equiv&
(z_3^2 \cdot y_2^4 \cdot y_1 \cdot y_2^4 \cdot z_3^{-1} \cdot y_2^{-2} ) \cdot z_3^{-2} \cdot
 (y_2^2 \cdot y_1^2) \cdot z_3^4 \cdot
(y_2^4 \cdot z_3^2). \ \Box
\end{array}
\end{eqnarray*}

}
\end{ex}

\begin{ex}
{\rm

We consider here the element $\nu$ given in Example \ref{example-unique-geodesic}
with $\delta=8$.
Then, $(z[0]+y[0])+(z[1]+y[1]) =7 < \delta < 9=(z[0]+y[0])+(z[1]+y[1])+z_3[2]$,
and we have
Case $(N,2M+1)$ with $N=2$ and $M=0$.
From (\ref{xy-1}) and (\ref{xy-4}),
we choose all the elements of $Z[0]\sqcup Y[0] \sqcup Z[1]\sqcup Y[1]$ and one element from $Z[2](=Z_3[2])$.
There are two choices, because $z_3[2]=2$.
Applying the suitable spread procedure in accordance with these choices,
we obtain
\begin{eqnarray*}
\begin{array}{rcl}
\nu&=&
(z_3^2\Delta^{-1}) \cdot (y_2^4\Delta^{-1}) \cdot
 (y_1\Delta^{-1}) \cdot 
(y_2^4\Delta^{-1})
 \cdot z_3^{-1} \cdot y_2^{-2} \cdot (z_3^4\Delta^{-1}) \\
&&  \cdot y_2^2 \cdot (y_1^2\Delta^{-1})
\cdot (z_3^4\Delta^{-1})
\cdot
  (y_2^4\Delta^{-1}) \cdot z_3^2
\\
&&=
z_3^{-4} \cdot y_2^{-3} \cdot y_1^{-2} \cdot y_2^{-3}
 \cdot z_3^{-1} \cdot y_2^{-2}
\cdot z_3^{-2}  \cdot y_2^2 \cdot y_1^{-1} \cdot z_3^{-2} \cdot
 y_2^{-3} \cdot z_3^2
=:\widetilde{\nu}_1
\\
&=&
z_3^2 \cdot (y_2^4\Delta^{-1}) \cdot
 (y_1\Delta^{-1}) \cdot 
(y_2^4\Delta^{-1})
 \cdot z_3^{-1} \cdot y_2^{-2} \cdot (z_3^4\Delta^{-1}) \\
&&  \cdot y_2^2 \cdot (y_1^2\Delta^{-1})
\cdot (z_3^4\Delta^{-1})
\cdot
  (y_2^4\Delta^{-1}) \cdot (z_3^2\Delta^{-1})
\\
&&=
z_3^2 \cdot y_2^{-3} \cdot y_1^{-2} \cdot y_2^{-3}
 \cdot z_3^{-1} \cdot y_2^{-2}
\cdot z_3^{-2}  \cdot y_2^2 \cdot y_1^{-1} \cdot z_3^{-2} \cdot
 y_2^{-3} \cdot z_3^{-4}
=:\widetilde{\nu}_2.
\\
\end{array}
\end{eqnarray*}
In this way,  from Proposition \ref{pr-sc-T3},
we obtain all the geodesic representatives of $\pi(\nu)$,
which are denoted as $\widetilde{\nu}_1$ and $\widetilde{\nu}_2$.
If we choose the leftmost element from $Z_3[2]$,
then we obtain $\widetilde{\nu}_1$, which is the output $\widehat{\nu}$ in this case:
\begin{eqnarray*}
\begin{array}{rcl}
\widehat{\nu} &\equiv&
\xi^{(1)} \cdot z_-
\cdot \xi^{(2)} \cdot z_+
\cdot \xi^{(3)}\\
&\equiv&
(\varepsilon) \cdot z_3^{-4} \cdot (y_2^{-3} \cdot y_1^{-2} \cdot y_2^{-3}
 \cdot z_3^{-1} \cdot y_2^{-2} 
\cdot z_3^{-2}  \cdot y_2^2 \cdot y_1^{-1} \cdot z_3^{-2} \cdot
 y_2^{-3}) \cdot z_3^2 \cdot(\varepsilon).  \ \Box
\end{array}
\end{eqnarray*}

}
\end{ex}

\noindent
 {\bf [Case {\boldmath $(N,2M+2)$}, in which {\boldmath $0 \leq N \leq p_1-2$} and {\boldmath $0 \leq M \leq n-m-1$},
or {\boldmath $N=p_1-1$ and $0 \leq M \leq n-m-2$}]}

\noindent
In this case,
from (\ref{xy-1}), (\ref{xy-2}), (\ref{xy-4}) and
 (\ref{rule-z}),
all the elements of 
$\left(\sqcup_{l=0}^{N-1} (Z[l] \sqcup Y[l])\right)
 \sqcup (\sqcup_{k=n-M}^n Z_k[N])$
are chosen.
Then we obtain a unique output $\widehat{\nu}$,
which has the form
\begin{eqnarray}\label{2-1}
\widehat{\nu}
\equiv
x_{i_1}^{a_1}  x_{i_2}^{a_2} \cdots
x_{i_{\kappa}}^{a_\kappa},
\end{eqnarray}
with the following conditions:
\begin{eqnarray}\label{2-2}
\left\{
\begin{array}{rl}
\bullet& \kappa \in {\bf N} \cup \{0\},\\
\bullet&
x_{i_j}=z_k\ \mbox{with}\ n-M \leq k \leq n\\
&\Longrightarrow 
\begin{array}{l}
-(r_k^- +N+1) \leq a_j \leq  
r_k^+ -N-1,\ 
a_j \neq 0,\\
\end{array}
\\
\bullet&
x_{i_j}=z_k\ \mbox{with}\ m+1 \leq k \leq n-M-1\\
&\Longrightarrow 
\begin{array}{l}
-(r_k^- +N) \leq a_j \leq  
r_k^+ -N,\ 
a_j \neq 0,\\
\end{array}
\\
\bullet&
x_{i_j}=y_k\ \mbox{with}\ 1 \leq k \leq m\\
&\Longrightarrow 
\begin{array}{l}
-(q_k^- +N) \leq a_j \leq  
q_k^+ -N,\
a_j \neq 0,\\
\end{array}
\\
\end{array}
\right.
\end{eqnarray}
and
\begin{eqnarray}
&\bullet &
\mbox{${}^\exists j \in \{1,\dots,\kappa\}$ s.t., 
$a_j=-(r^-+N+1)$},\label{2-3}\\
&\bullet & 
\left\{
\begin{array}{l}\label{2-4}
\mbox{$\widehat{\nu}$ contains at least one of the following sub-words:}\\
\begin{array}{l}

z_k^{r_k^-+1},
z_k^{r_k^-+2},
\dots,
z_k^{r_k^+-N-1}
\quad (k \in \{m+1,\dots, n\}),\\

z_k^{r_k^+-N}
\quad (k \in \{m+1,\dots,n-M-1\}),\\

y_k^{q_k^-+1},
y_k^{q_k^-+2},
\dots,
y_k^{q_k^+-N}
\quad (k \in \{1,\dots,m\}).\\
\end{array}
\end{array}
\right.
\end{eqnarray}
Here, note the following: 
(i) here too we write $r_{n-M}$,
$r^+_{n-M}$ and $r^-_{n-M}$ as $r$, $r^+$ and $r^-$,
and $\widehat{\nu}$ must satisfy the condition (\ref{2-3}),
because here we have $z_{n-M}[N]>0$ and $(r^+-N)-r
=-(r^- +N+1)$;
(ii) $\widehat{\nu}$ must satisfy the condition (\ref{2-4}),
because, from the assumption ${\bf r}_\nu > \delta$,
there exists at least one element in
$(\sqcup_{l=N+1}^{p_1-1} Z[l])\sqcup (\sqcup_{k=m+1}^{n-M-1} Z_k[N]) 
\sqcup (\sqcup_{l=N}^{p_1-2} Y[l])$,
and we have $r_k^+-(p_1-1)
=r_k^- +1$ $(k \in \{m+1,\dots,n\})$
and
$q_k^+-(p_1-2)
=q_k^- +1$ $(k \in \{1,\dots,m\})$.

Then, for each pair $(N,M)$ 
with $0 \leq N \leq p_1-2$ and $0 \leq M \leq n-m-1$
or $N =p_1-1$ and $0 \leq M \leq n-m-2$,
we define the following set:
\begin{eqnarray*}
{\Gamma}^{(N,2M+2)}_{3^0}:=
\{
\widehat{\nu} \in \Sigma^* \ ; \
\mbox{$\widehat{\nu}$ is given in (\ref{2-1})
with the conditions in (\ref{2-2})--(\ref{2-4})} \}.
\end{eqnarray*}

\vspace{0.2cm}

\begin{ex}
{\rm

Here, we choose the element $\nu$ given in Example \ref{example-unique-geodesic}
with $\delta=2$.
Then, $\delta = z_3[0]$ and $z_3[0] >0$,
and we have Case $(N,2M+2)$ with $N=0$ and $M=0$.
From (\ref{xy-1}) and (\ref{xy-2}),
we choose all the elements from $Z_3[0]$.
Applying the suitable spread procedure in accordance with this choice,
we obtain
\begin{eqnarray*}
\begin{array}{rcl}
\nu&=&
z_3^2 \cdot y_2^4  \cdot y_1 \cdot y_2^4
 \cdot z_3^{-1} \cdot y_2^{-2} \cdot (z_3^4\Delta^{-1}) 
\cdot y_2^2  \cdot y_1^2 \cdot  (z_3^4\Delta^{-1}) \cdot
 y_2^4 \cdot z_3^2
\\
&&=
z_3^2 \cdot y_2^4 \cdot y_1 \cdot y_2^4
 \cdot z_3^{-1} \cdot y_2^{-2} 
\cdot z_3^{-2}  \cdot y_2^2 \cdot y_1^2 \cdot  z_3^{-2} \cdot
y_2^4 \cdot z_3^2
=:\widetilde{\nu}_1.
\\
\end{array}
\end{eqnarray*}
From Proposition \ref{pr-sc-T3},
for $\pi(\nu)$,
there is a single geodesic representative, which is denoted as $\widetilde{\nu}_1$ here
and is the output $\widehat{\nu}$ in this case.  $\Box$

}
\end{ex}

\noindent
{\bf [Case {\boldmath $(N,2M+1)$}, in which {\boldmath $0 \leq N \leq p_1-2$} and {\boldmath $n-m \leq  M \leq n-1$}]}

\noindent
In this case,
we perform the suitable-spread procedure 
in accordance with the following rules.
\begin{enumerate}
\item
All the elements of 
$
(\sqcup_{l=0}^{N} Z[l])
\sqcup
(\sqcup_{l=0}^{N-1} Y[l]) \sqcup (\sqcup_{k=n-(M-1)}^{m} Y_k[N])
$
are chosen,
because we have the relations in (\ref{xy-2})--(\ref{xy-4})
and adopt the rule  (\ref{rule-y}).
\item
Define the quantity
$\psi:=\delta-\left(\sum_{l=0}^{N} z[l]+\sum_{l=0}^{N-1}y[l]
+\sum_{k=n-(M-1)}^{m}y_k[N]\right)$
(see Figure 3).
In this case, we have $0 < \psi < y_{n-M}[N]$.
Hence, from (\ref{xy-3}), (\ref{xy-4}) and  (\ref{rule-y}), 
we choose $\psi$ sub-words from $Y_{n-M}[N]$.
Here, for conciseness, we write $q_{n-M}$,
$q^+_{n-M}$ and $q^-_{n-M}$ as $q$, $q^+$ and $q^-$.
\item
When choosing the $\psi$ sub-words from $Y_{n-M}[N]$,
we always choose the leftmost of these.
We represent all the sub-words belonging to $Y_{n-M}[N]$
by $y_{n-M}^{\overline{\alpha}_{j_1}}, \dots, y_{n-M}^{\overline{\alpha}_{j_{y_{n-M}[N]}}}$.
Then,
the sub-words chosen through this procedure are
$y_{n-M}^{\overline{\alpha}_{j_1}}, \dots, y_{n-M}^{\overline{\alpha}_{j_\psi}}$.
\end{enumerate}

\setlength{\unitlength}{0.7mm}
\begin{picture}(200,40)(-30,170)

\put(-4,193){\makebox(50,10){$\sum_{l=0}^{N} z[l]+\sum_{l=0}^{N-1}y[l]
+\sum_{k=n-(M-1)}^{m}y_k[N]$}}
\put(63,193){\makebox(50,10){$\delta$}}
\put(130,193){\makebox(50,10){$\sum_{l=0}^{N} z[l]+\sum_{l=0}^{N-1}y[l]
+\sum_{k=n-M}^{m}y_k[N]$}}
\put(0,190){\line(1,0){165}}

\put(15,187.5){\makebox(5,5){$\bullet$}}
\put(86,187.5){\makebox(5,5){$\bullet$}}
\put(153,187.5){\makebox(5,5){$\bullet$}}
\put(48,175){\makebox(10,10){$\psi$}}
\put(90,167){\makebox(10,10){$y_{n-M}[N]$}}
\put(17.4,170){\line(0,1){19}}
\put(42.4,180){\vector(-1,0){25}}
\put(63,180){\vector(1,0){25}}
\put(88.5,176){\line(0,1){14}}
\put(78,173){\vector(-1,0){60}}
\put(110,173){\vector(1,0){45}}
\put(155.5,170){\line(0,1){19}}
%

\end{picture}
\begin{center}
Figure 3.
The location of $\delta$ and the definition of $\psi$.
\end{center}

\noindent
In this way,
we obtain an output $\widehat{\nu}$ containing 
$y_{n-M}^{\overline{\alpha}_{j_1}-q}, \dots, y_{n-M}^{\overline{\alpha}_{j_\psi}-q}$
and $y_{n-M}^{\overline{\alpha}_{j_{\psi+1}}},\dots, 
y_{n-M}^{\overline{\alpha}_{j_{y_{n-M}[N]}}}$.
Here, note that
$y_{n-M}^{\overline{\alpha}_{j_s}-q} \equiv
 y_{n-M}^{-(q^- +N+1)}$ ($s \in \{1,\dots,\psi\}$),
because $(q^+ -N)-q=
-(q^- +N+1)$.
Then, choosing $y_-$ (resp., $y_+$) to be the leftmost element,
$y_{n-M}^{\overline{\alpha}_{j_1}-q}$ (resp., $y_{n-M}^{\overline{\alpha}_{j_{\psi+1}}}$),
we can write the output $\widehat{\nu}$ as 
\begin{eqnarray}\label{3-1}
\widehat{\nu} \equiv 
\xi^{(4)} \cdot y_-
\cdot \xi^{(5)} \cdot y_+ \cdot \xi^{(6)},
\end{eqnarray}
where $\xi^{(4)}$, $\xi^{(5)}$ and $\xi^{(6)}$ are given by
\begin{eqnarray}\label{3-2}
\left\{
\begin{array}{l}
\xi^{(4)} :=
x_{i_1}^{a_1} x_{i_2}^{a_2} \cdots
 x_{i_{\kappa_1}}^{a_{\kappa_1}},\\
\xi^{(5)} :=
x_{i_1}^{b_1} x_{i_2}^{b_2} \cdots
x_{i_{\kappa_2}}^{b_{\kappa_2}},\\
\xi^{(6)} :=
x_{i_1}^{c_1} x_{i_2}^{c_2} \cdots
x_{i_{\kappa_3}}^{c_{\kappa_3}},\\
\end{array}
\right.
\end{eqnarray}
with the following conditions:

\begin{eqnarray}
&\xi^{(4)}&
\left\{
\begin{array}{rl}\label{3-3}
\bullet& \kappa_1\in {\bf N} \cup \{0\},\\
\bullet& x_{i_{\kappa_1}}\neq y_{n-M},\\
\bullet&
x_{i_j}=z_k\ \mbox{with}\ m+1 \leq k \leq n\\
&\Longrightarrow
\begin{array}{l}
-(r_k^- +N+1) \leq a_j \leq  
r_k^+ -N-1,\  
a_j \neq 0,\\
\end{array}
\\
\bullet&
x_{i_j}=y_k\ \mbox{with}\ n-M+1 \leq k \leq m\\
&\Longrightarrow 
\begin{array}{l}
-(q_k^- +N+1) \leq a_j \leq  
q_k^+ -N-1,\ 
a_j \neq 0,\\
\end{array}
\\
\bullet&
x_{i_j}=y_k\ \mbox{with}\ k=n-M\\
&\Longrightarrow 
\begin{array}{l}
-(q^- +N) \leq a_j \leq  
q^+ -N-1,\
a_j \neq 0,\\
\end{array}
\\
\bullet&
x_{i_j}=y_k\ \mbox{with}\ 1 \leq k \leq n-M-1\\
&\Longrightarrow 
\begin{array}{l}
-(q_k^- +N) \leq a_j \leq  
q_k^+ -N,\
a_j \neq 0,\\
\end{array}
\\
\end{array}
\right.\\
&\xi^{(5)}&
\left\{
\begin{array}{rl}\label{3-4}
\bullet& \kappa_2\in {\bf N},\\
\bullet& x_{i_1} \neq y_{n-M}, x_{i_{\kappa_2}}\neq y_{n-M},\\
\bullet&
x_{i_j}=z_k\ \mbox{with}\ m+1 \leq k \leq n\\
&\Longrightarrow 
\begin{array}{l}
-(r_k^- +N+1) \leq b_j \leq  
r_k^+ -N-1,\ 
b_j \neq 0,\\
\end{array}
\\
\bullet&
x_{i_j}=y_k\ \mbox{with}\ n-M+1 \leq k \leq m\\
&\Longrightarrow 
\begin{array}{l}
-(q_k^- +N+1) \leq b_j \leq  
q_k^+ -N-1,\  
b_j \neq 0,\\
\end{array}
\\
\bullet&
x_{i_j}=y_k\ \mbox{with}\ k=n-M\\
&\Longrightarrow 
\begin{array}{l}
-(q^- +N+1) \leq b_j \leq  
q^+ -N-1,\ 
b_j \neq 0,\\
\end{array}
\\
\bullet&
x_{i_j}=y_k\ \mbox{with}\ 1 \leq k \leq n-M-1\\
&\Longrightarrow 
\begin{array}{l}
-(q_k^- +N) \leq b_j \leq  
q_k^+ -N,\
b_j \neq 0,\\
\end{array}
\\
\end{array}
\right.\\
&\xi^{(6)}&
\left\{
\begin{array}{rl}\label{3-5}
\bullet& \kappa_3 \in {\bf N} \cup \{0\},\\
\bullet& x_{i_1} \neq y_{n-M},\\
\bullet&
x_{i_j}=z_k\ \mbox{with}\ m+1 \leq k \leq n\\
&\Longrightarrow 
\begin{array}{l}
-(r_k^- +N+1) \leq c_j \leq  
r_k^+ -N-1,\ 
c_j \neq 0,\\
\end{array}
\\
\bullet&
x_{i_j}=y_k\ \mbox{with}\ n-M+1 \leq k \leq m\\
&\Longrightarrow 
\begin{array}{l}
-(q_k^- +N+1) \leq c_j \leq  
q_k^+ -N-1,\ 
c_j \neq 0,\\
\end{array}
\\
\bullet&
x_{i_j}=y_k\ \mbox{with}\ k=n-M\\
&\Longrightarrow 
\begin{array}{l}
-(q^- +N) \leq c_j \leq  
q^+ -N,\ 
c_j \neq 0,\\
\end{array}
\\
\bullet&
x_{i_j}=y_k\ \mbox{with}\ 1 \leq k \leq n-M-1\\
&\Longrightarrow 
\begin{array}{l}
-(q_k^- +N) \leq c_j \leq  
q_k^+ -N,\
c_j \neq 0.\\
\end{array}
\\
\end{array}
\right.
\end{eqnarray}
Then, 
for each pair $(N,M)$ 
satisfying $0 \leq N \leq p_1-2$ and $n-m \leq M \leq n-1$,
we define
\begin{eqnarray*}
\begin{array}{rl}
{\Gamma}^{(N,2M+1)}_{3^0}:=
\{
\widehat{\nu} \in \Sigma^* \ ; 
&
\widehat{\nu} \equiv
\xi^{(4)} \cdot y_{n-M}^{-(q^- +N+1)}
\cdot \xi^{(5)} \cdot 
y_{n-M}^{q^+ -N}
\cdot \xi^{(6)},\\
&\xi^{(4)}, \xi^{(5)} \ \mbox{and}\ \xi^{(6)} 
\ \mbox{are given in (\ref{3-2})} 
\ \mbox{with the conditions in (\ref{3-3})--(\ref{3-5})} \}.
\end{array}
\end{eqnarray*}

\begin{ex}
{\rm 

Here, we choose the element $\nu$ given in Example \ref{example-unique-geodesic}
with $\delta=4$.
Then, $z[0]=2 < \delta < 5=z[0]+y_2[0]$,
and we have
Case $(N,2M+1)$ with $N=0$ and $M=1$.
From (\ref{xy-2})--(\ref{xy-4}) and  (\ref{rule-y}),
we choose all the elements of $Z[0]$ and two elements from $Y_2[0]$.
There are three choices, because $y_2[0]=3$.
Applying the suitable spread procedure in accordance with these choices,
we obtain
\begin{eqnarray*}
\begin{array}{rcl}
\nu&=&
z_3^2 \cdot (y_2^4\Delta^{-1}) \cdot y_1 \cdot 
(y_2^4\Delta^{-1})
 \cdot z_3^{-1} \cdot y_2^{-2} \cdot (z_3^4\Delta^{-1}) 
\cdot y_2^2 \cdot y_1^2 \cdot  (z_3^4\Delta^{-1}) \cdot
y_2^4 \cdot z_3^2
\\
&&=
z_3^2 \cdot y_2^{-3} \cdot y_1 \cdot y_2^{-3}
 \cdot z_3^{-1} \cdot y_2^{-2} 
\cdot z_3^{-2}  \cdot y_2^2 \cdot y_1^2 \cdot  z_3^{-2} \cdot
y_2^4 \cdot z_3^2
=:\widetilde{\nu}_1
\\
&=&
z_3^2 \cdot (y_2^4\Delta^{-1}) \cdot y_1 \cdot y_2^4
 \cdot z_3^{-1} \cdot y_2^{-2} \cdot (z_3^4\Delta^{-1}) 
\cdot y_2^2 \cdot y_1^2 \cdot  (z_3^4\Delta^{-1}) \cdot
 (y_2^4\Delta^{-1}) \cdot z_3^2
\\
&&=
z_3^2 \cdot y_2^{-3} \cdot y_1 \cdot y_2^4
 \cdot z_3^{-1} \cdot y_2^{-2} 
\cdot z_3^{-2} \cdot y_2^2 \cdot y_1^2 \cdot z_3^{-2} \cdot
y_2^{-3} \cdot z_3^2
=:\widetilde{\nu}_2
\\
&=&
z_3^2 \cdot y_2^4 \cdot y_1  \cdot (y_2^4\Delta^{-1})
 \cdot z_3^{-1} \cdot y_2^{-2} \cdot (z_3^4\Delta^{-1})
\cdot y_2^2  \cdot y_1^2 \cdot (z_3^4\Delta^{-1}) \cdot
 (y_2^4\Delta^{-1}) \cdot z_3^2
\\
&&=
z_3^2 \cdot y_2^4 \cdot y_1 \cdot y_2^{-3}
 \cdot z_3^{-1} \cdot y_2^{-2} 
\cdot z_3^{-2}  \cdot y_2^2 \cdot y_1^2 \cdot z_3^{-2} \cdot
y_2^{-3} \cdot z_3^2
=:\widetilde{\nu}_3.
\\
\end{array}
\end{eqnarray*}
If we choose the leftmost elements from $Y_2[0]$,
then we obtain $\widetilde{\nu}_1$, which is the output $\widehat{\nu}$ in this case:
\begin{eqnarray*}
\begin{array}{rcl}
\widehat{\nu} &\equiv&
\xi^{(4)} \cdot y_-
\cdot \xi^{(5)} \cdot y_+
\cdot \xi^{(6)}\\
&\equiv&
(z_3^2) \cdot y_2^{-3} \cdot (y_1 \cdot y_2^{-3}
 \cdot z_3^{-1} \cdot y_2^{-2} 
\cdot z_3^{-2}  \cdot y_2^2 \cdot y_1^2 \cdot  z_3^{-2}) \cdot
y_2^4 \cdot (z_3^2).
\end{array}
\end{eqnarray*}
Note that $\nu$, $\widetilde{\nu}_1$, $\widetilde{\nu}_2$ and $\widetilde{\nu}_3$
are the same elements considered in Example \ref{example-SS}
if we simply replace $y_1$ with $x_1$, $y_2$ with $x_3$ and $z_3$ with $x_2$.
From Proposition \ref{pr-sc-T3},
the set of all geodesic representatives of $\pi(\nu)$
consists of $\widetilde{\nu}_1,\dots,\widetilde{\nu}_6$,
which are obtained in Example \ref{example-SS}. $\Box$

}
\end{ex}

\vspace{0.2cm}

\noindent
{\bf [Case  {\boldmath $(N,2M+2)$}, in which {\boldmath $0 \leq N \leq p_1-2$} and {\boldmath $n-m \leq M \leq n-1$}]}

\noindent
In this case,
all the elements of 
$\left(\sqcup_{l=0}^{N}Z[l]\right)
\sqcup
\left(\sqcup_{k=0}^{N-1} Y[l]\right)
\sqcup
\left(\sqcup_{k=n-M}^{m} y_k[N]\right)$
 are chosen,
because we have the relations in (\ref{xy-2})--(\ref{xy-4})
and  (\ref{rule-y}).
Then we obtain a unique output $\widehat{\nu}$,
which has the form
\begin{eqnarray}\label{4-1}
\widehat{\nu}
\equiv
x_{i_1}^{a_1}  x_{i_2}^{a_2} \cdots
x_{i_{\kappa}}^{a_\kappa},
\end{eqnarray}
with the following conditions:
\begin{eqnarray}\label{4-2}
\left\{
\begin{array}{rl}
\bullet& \kappa \in {\bf N} \cup \{0\},\\
\bullet&
x_{i_j}=z_k\ \mbox{with}\ m+1 \leq k \leq n\\
&\Longrightarrow 
\begin{array}{l}
-(r_k^- +N+1) \leq a_j \leq  
r_k^+ -N-1,\ 
a_j \neq 0,\\
\end{array}
\\
\bullet&
x_{i_j}=y_k\ \mbox{with}\ n-M \leq k \leq m\\
&\Longrightarrow 
\begin{array}{l}
-(q_k^- +N+1) \leq a_j \leq  
q_k^+ -N-1,\ 
a_j \neq 0,\\
\end{array}
\\
\bullet&
x_{i_j}=y_k\ \mbox{with}\ 1 \leq k \leq n-M-1\\
&\Longrightarrow 
\begin{array}{l}
-(q_k^- +N) \leq a_j \leq  
q_k^+ -N,\
a_j \neq 0,\\
\end{array}
\\
\end{array}
\right.
\end{eqnarray}
and
\begin{eqnarray}
&\bullet &
\mbox{${}^\exists j \in \{1,\dots,\kappa\}$ s.t., 
$a_j=-(q^-+N+1)$},\label{4-3}\\
&\bullet & 
\left\{
\begin{array}{l}\label{4-4}
\mbox{$\widehat{\nu}$ contains at least one of the following sub-words:}\\
\begin{array}{l}

z_k^{r_k^-+1},
z_k^{r_k^-+2},
\dots,
z_k^{r_k^+-N-1}
\quad (k \in \{m+1,\dots, n\}),\\

y_k^{q_k^-+1},
y_k^{q_k^-+2},
\dots,
y_k^{q_k^+-N-1}
\quad (k \in \{1,\dots,m\}),\\

y_k^{q_k^+-N}
\quad (k \in \{1,\dots,n-M-1\}).\\
\end{array}
\end{array}
\right.
\end{eqnarray}
Here, note the following: 
(i) here too we write $q_{n-M}$,
$q^+_{n-M}$ and $q^-_{n-M}$ as $q$, $q^+$ and $q^-$,
and $\widehat{\nu}$ must satisfy the condition (\ref{4-3}),
because here we have $y_{n-M}[N]>0$ and $(q^+-N)-q
=-(q^- +N+1)$;
(ii) $\widehat{\nu}$ must satisfy the condition (\ref{4-4}),
because, from the assumption ${\bf r}_\nu > \delta$,
there exists at least one element in
$(\sqcup_{l=N+1}^{p_1-1} Z[l]) 
\sqcup (\sqcup_{l=N}^{p_1-2} Y[l])
\sqcup (\sqcup_{k=1}^{n-M-1} Y_k[N])$,
and we have $r_k^+-(p_1-1)
=r_k^- +1$ $(k \in \{m+1,\dots,n\})$
and
$q_k^+-(p_1-2)
=q_k^- +1$ $(k \in \{1,\dots,m\})$.

Then, for each pair $(N,M)$ 
with $0 \leq N \leq p_1-2$ and $n-m \leq M \leq n-1$,
we define the following set:
\begin{eqnarray*}
{\Gamma}^{(N,2M+2)}_{3^0}:=
\{
\widehat{\nu} \in \Sigma^* \ ; \ 
\mbox{$\widehat{\nu}$ is given in (\ref{4-1})
with the conditions in (\ref{4-2})--(\ref{4-4})} \}.
\end{eqnarray*}

\vspace{0.2cm}

\begin{ex}
{\rm

Here, we choose the element $\nu$ given in Example \ref{example-unique-geodesic}
with $\delta=5$.
Then, $\delta = z[0]+y_2[0]$ and $y_2[0] >0$,
and we have
Case $(N,2M+2)$ with $N=0$ and $M=1$.
From (\ref{xy-2})--(\ref{xy-4}) and  (\ref{rule-y}),
we choose all the elements from $Z[0] \sqcup Y_2[0]$.
Applying the suitable spread procedure in accordance with this choice,
we obtain
\begin{eqnarray*}
\begin{array}{rcl}
\nu&=&
z_3^2 \cdot (y_2^4\Delta^{-1}) \cdot y_1 \cdot 
(y_2^4\Delta^{-1})
 \cdot z_3^{-1} \cdot y_2^{-2} \cdot (z_3^4\Delta^{-1}) 
\cdot y_2^2 \cdot y_1^2 \cdot (z_3^4\Delta^{-1}) \cdot
  (y_2^4\Delta^{-1}) \cdot z_3^2
\\
&&=
z_3^2 \cdot y_2^{-3} \cdot y_1 \cdot y_2^{-3}
 \cdot z_3^{-1} \cdot y_2^{-2} 
\cdot z_3^{-2}  \cdot y_2^2 \cdot y_1^2 \cdot z_3^{-2} \cdot
y_2^{-3} \cdot z_3^2
=:\widetilde{\nu}_1.
\\
\end{array}
\end{eqnarray*}
The geodesic representative, $\widetilde{\nu}_1$,
is the output $\widehat{\nu}$ in this case.
Note that $\nu$ contains an element $y_1^2 \in Y_1[0] \ (\subset Y[0])$
and we have $y[0]=4$.
Therefore, there are four ways to choose three elements from $Y[0]$.
Thus, from Propotition \ref{pr-sc-T3},
there are four geodesic representatives of $\nu$.
$\Box$
\\

}
\end{ex}

\begin{ex}
{\rm

Here, we choose the element $\nu$ given in Example \ref{example-unique-geodesic}
with $\delta=6$.
Then, $\delta = z[0]+y_2[0]+y_1[0]$ and $y_1[0] >0$,
and we have
Case $(N,2M+2)$ with $N=0$ and $M=2$.
From (\ref{xy-2})--(\ref{xy-4}),
we choose all the elements from $Z[0] \sqcup Y_2[0] \sqcup Y_1[0]$.
Applying the suitable spread procedure in accordance with this choice,
we obtain
\begin{eqnarray*}
\begin{array}{rcl}
\nu&=&
z_3^2 \cdot (y_2^4\Delta^{-1}) \cdot y_1 \cdot 
(y_2^4\Delta^{-1})
 \cdot z_3^{-1} \cdot y_2^{-2} \cdot (z_3^4\Delta^{-1}) 
\cdot y_2^2 \cdot (y_1^2\Delta^{-1}) \cdot (z_3^4\Delta^{-1})\cdot
  (y_2^4\Delta^{-1}) \cdot z_3^2
\\
&&=
z_3^2 \cdot y_2^{-3} \cdot y_1 \cdot y_2^{-3}
 \cdot z_3^{-1} \cdot y_2^{-2} 
\cdot z_3^{-2}  \cdot y_2^2 \cdot y_1^{-1} \cdot z_3^{-2} \cdot
y_2^{-3} \cdot z_3^2
=:\widetilde{\nu}_1.
\\
\end{array}
\end{eqnarray*}
From Proposition \ref{pr-sc-T3},
for $\pi(\nu)$,
there is a single geodesic representative, which is denoted $\widetilde{\nu}_1$ here
and is the output $\widehat{\nu}$ in this case. $\Box$
\\

}
\end{ex}

\begin{ex}
{\rm

Here, we choose the element $\nu$ given in Example \ref{example-unique-geodesic}
with $\delta=7$.
Then, $\delta = (z[0]+y[0])+z[1]+y_2[1]+y_1[1]$ and $y_1[1] >0$,
and we have
Case $(N,2M+2)$ with $N=1$ and $M=2$.
From (\ref{xy-2})--(\ref{xy-4}),
we choose all the elements from $(Z[0] \sqcup Y[0]) \sqcup Z[1] \sqcup Y_2[1] \sqcup Y_1[1]$.
Applying the suitable spread procedure in accordance with this choice,
we obtain
\begin{eqnarray*}
\begin{array}{rcl}
\nu&=&
z_3^2 \cdot (y_2^4\Delta^{-1}) \cdot (y_1\Delta^{-1})  \cdot 
(y_2^4\Delta^{-1})
 \cdot z_3^{-1} \cdot y_2^{-2} \cdot (z_3^4\Delta^{-1})
\cdot y_2^2 \cdot (y_1^2\Delta^{-1}) \cdot (z_3^4\Delta^{-1})\\
&&
\cdot
  (y_2^4\Delta^{-1}) \cdot z_3^2
\\
&&=
z_3^2 \cdot y_2^{-3} \cdot y_1^{-2} \cdot y_2^{-3}
 \cdot z_3^{-1} \cdot y_2^{-2} 
\cdot z_3^{-2}  \cdot y_2^2 \cdot y_1^{-1} \cdot z_3^{-2} \cdot
y_2^{-3} \cdot z_3^2
=:\widetilde{\nu}_1.
\\
\end{array}
\end{eqnarray*}
From Proposition \ref{pr-sc-T3},
for $\pi(\nu)$,
there is a single geodesic representative, which is denoted as $\widetilde{\nu}_1$ here
and is the output $\widehat{\nu}$ in this case. $\Box$
\\

}
\end{ex}

\subsubsection{Case: {\boldmath $p_k-p_1$} is even for every $k \in \{ 1, \dots, n\}$}\label{section:even}

In this section, we consider the case $m=n$,
i.e., the case in which $p_k-p_1$ is even for every $k \in \{1, \dots,n\}$.
Then, as in the previous section,
from (\ref{xy-5}), we have one of the following cases.
\vspace{0.3cm}

\noindent
$\bullet$ Case  $(N,2M+1)$, in which $0 \leq N \leq p_1-2$ and 
$0 \leq M \leq n-1$:
\begin{eqnarray*}
{\displaystyle \sum_{l=0}^{N-1} y[l] \quad 
+ \sum_{k=n-(M-1)}^{n} y_k[N]}
\quad < \quad \delta \quad
< \quad {\displaystyle \sum_{l=0}^{N-1} y[l]
\quad +\sum_{k=n-M}^{n} y_k[N]};
\end{eqnarray*}

\noindent
$\bullet$ Case  $(N,2M+2)$, in which $0 \leq N \leq p_1-3$ and $0 \leq M \leq n-1$,
or $N=p_1-2$ and $0 \leq M \leq n-2$:
\begin{eqnarray*}
{\displaystyle \delta = \sum_{l=0}^{N-1} y[l] \quad
+\sum_{k=n-M}^{n} y_k[N]}
\quad
\mbox{with} \quad y_{n-M}[N] > 0.
\end{eqnarray*}

\vspace{0.2cm}

\noindent
{\bf [Case  {\boldmath $(N,2M+1)$}, in which {\boldmath $0 \leq N \leq p_1-2$} and {\boldmath $0 \leq M \leq n-1$}]}

\noindent
This case is obtained by setting
$m=n$ 
in Case  $(N,2M+1)$ (with $0 \leq N \leq p_1-2$, $n-m \leq M \leq n-1$)
considered in the previous section.
Let us define ${\bar\xi}^{(4)}$, ${\bar\xi}^{(5)}$ and ${\bar\xi}^{(6)}$ as
\begin{eqnarray}\label{e3-1}
\left\{
\begin{array}{l}
{\bar\xi}^{(4)}:=\xi^{(4)}\quad \mbox{with}\quad m:=n,\\
{\bar\xi}^{(5)}:=\xi^{(5)}\quad \mbox{with}\quad m:=n,\\
{\bar\xi}^{(6)}:=\xi^{(6)}\quad \mbox{with}\quad m:=n,
\end{array}
\right.
\end{eqnarray}
where $\xi^{(4)}$, $\xi^{(5)}$ and $\xi^{(6)}$ are given in (\ref{3-2})
with the conditions (\ref{3-3})--(\ref{3-5}).

Now, for each pair $(N,M)$ 
with $0 \leq N \leq p_1-2$ and $0 \leq M \leq n-1$,
we define the following set:
\begin{eqnarray*}
\begin{array}{rl}
{\bar{\Gamma}}^{(N,2M+1)}_{3^0}:=
\{
\bar{\widehat{\nu}} \in \Sigma^* \ ;
&
\bar{\widehat{\nu}}
\equiv
{\bar\xi}^{(4)} \cdot y_{n-M}^{-(q^-+N+1)}
\cdot {\bar\xi}^{(5)} \cdot y_{n-M}^{q^+ -N}
\cdot {\bar\xi}^{(6)},\\
&{\bar\xi}^{(4)}, {\bar\xi}^{(5)} \ \mbox{and}\ {\bar\xi}^{(6)} 
\ \mbox{are given in (\ref{e3-1})} \}.
\end{array}
\end{eqnarray*}

\vspace{0.2cm}

\noindent
 {\bf [Case {\boldmath $(N,2M+2)$}, in which {\boldmath $0 \leq N \leq p_1-3$} and {\boldmath $0 \leq M \leq n-1$},
or {\boldmath $N=p_1-2$ and $0 \leq M \leq n-2$}]}

\noindent
This case is obtained by setting $m=n$ 
in Case $(N,2M+2)$ (with $0 \leq N \leq p_1-2$, $n-m \leq M \leq n-1$)
considered in the previous section 
and appropriately by modifying the ranges of $N$ and $M$.
Let us define $\bar{\widehat{\nu}}$ as
\begin{eqnarray}\label{e4-1}
\bar{\widehat{\nu}}:=\widehat{\nu}\quad \mbox{with}\quad m:=n,
\end{eqnarray}
where $\widehat{\nu}$ is given in (\ref{4-1}) with the conditions in (\ref{4-2})--(\ref{4-4}).

Now, for each pair $(N,M)$ 
with $0 \leq N \leq p_1-3$ and $0 \leq M \leq n-1$,
or $N =p_1-2$ and $0 \leq M \leq n-2$,
we define the following set:
\begin{eqnarray*}
{\bar{\Gamma}}^{(N,2M+2)}_{3^0}:=
\{
\bar{\widehat{\nu}} \in \Sigma^* \ ; \
\mbox{$\bar{\widehat{\nu}}$ is given in (\ref{e4-1})} \}.
\end{eqnarray*}

\vspace{0.2cm}

\subsubsection{A set consisting of unique geodesic representatives
 of elements of $G(p_1,\dots,p_n)$}\label{theorem4.4}

Collecting the cases considered in Sections \ref{section:even and odd} and \ref{section:even},
we define ${\Gamma}_{3^0}$ as follows:

\noindent
$\bullet$ If $m<n$, then 
\begin{eqnarray}\label{partition-3^0}
\begin{array}{rl} 
{\Gamma}_{3^0}:=&{\displaystyle
\bigsqcup_{N=0}^{p_1-1}
\bigsqcup_{M=0}^{n-m-1}
{\Gamma}^{(N,2M+1)}_{3^0}
\sqcup
\bigsqcup_{N=0}^{p_1-2}
\bigsqcup_{M=0}^{n-m-1}
{\Gamma}^{(N,2M+2)}_{3^0}
\sqcup
\bigsqcup_{M=0}^{n-m-2}
{\Gamma}^{(p_1-1,2M+2)}_{3^0}}\\
&
{\displaystyle \sqcup
\bigsqcup_{N=0}^{p_1-2}
\bigsqcup_{M=n-m}^{n-1}
{\Gamma}^{(N,2M+1)}_{3^0}
\sqcup
\bigsqcup_{N=0}^{p_1-2}
\bigsqcup_{M=n-m}^{n-1}
{\Gamma}^{(N,2M+2)}_{3^0}}.
\end{array}
\end{eqnarray}

\noindent
$\bullet$ If $m=n$, then
\begin{eqnarray}\label{epartition-3^0}
{\Gamma}_{3^0}:={\displaystyle
\bigsqcup_{N=0}^{p_1-2}
\bigsqcup_{M=0}^{n-1}
{\bar{\Gamma}}^{(N,2M+1)}_{3^0}
\sqcup
\bigsqcup_{N=0}^{p_1-3}
\bigsqcup_{M=0}^{n-1}
{\bar{\Gamma}}^{(N,2M+2)}_{3^0}
\sqcup
\bigsqcup_{M=0}^{n-2}
{\bar{\Gamma}}^{(p_1-2,2M+2)}_{3^0}}.
\end{eqnarray}

\noindent
Then, applying the arguments given in the previous sections,
we obtain the following theorem: 

\begin{thm}\label{th-ugr}
Let $g$ be an element of ${G(p_1,\dots,p_n)}_{3^0}$.
Then
there exists a unique geodesic representative of $g$ in the set ${\Gamma}_{3^0}$.
Moreover, the restriction of $\pi$ to ${\Gamma}_{3^0}$ is a bijective map to 
${G(p_1,\dots,p_n)}_{3^0}$.
\end{thm}

Finally, 
 we define ${\Gamma}$ to be the following disjoint union:
\begin{eqnarray}\label{Gamma}
{\Gamma}:=
{\Gamma}_{1} \sqcup 
{\Gamma}_{2} \sqcup
{\Gamma}_{3^+ \cup 3^-} \sqcup
{\Gamma}_{3^0}.
\end{eqnarray}
From Theorem \ref{th-ugr} and the argument given in Section \ref{section-1,2,3^+,3^-},
we conclude that
the restriction map of $\pi$ to ${\Gamma}$ is a bijective map to
$G(p_1,\dots,p_n)$.

\section{Computation of growth series}\label{section-growth}

In this section,
we compute the spherical growth series for the sets ${\Gamma}_I$
$(I \in \{1,2,3^+,3^-,3^+\cap 3^-,3^+\cup 3^-,3^0\})$
defined in Section \ref{section-ugr}
and obtain a rational function expression for the spherical growth
series of the group $G(p_1,\dots,p_n)$.

Let ${\Lambda}$ be a subset of $\Sigma^*$.
Assume that
each element $\xi$ of ${\Lambda}$ is geodesic
and that there is no other element $\xi' \in {\Lambda}$ satisfying $\pi(\xi')=\pi(\xi)$.
We define the spherical growth series for the set ${\Lambda}$ as follows:
\begin{eqnarray}
{\cal S}_{\Lambda}(t):=
\sum_{l=0}^{\infty}
\# \{\xi \in {\Lambda} \ ; \ |\xi|=l \} \ t^l.
\end{eqnarray}
From the argument presented in the previous section,
we know that the above assumption holds for
the sets appearing on the right-hand sides of
 (\ref{partition-3^0}) and (\ref{epartition-3^0}).
Hence, we can define their spherical growth series.

\vspace{0.2cm}

In order to simplify the presentation of the growth series,
we introduce 
\begin{eqnarray}\label{def-T}
T_0:=0 \quad \mbox{and}\quad
T_u:=t+t^2+\cdots +t^u \ \ \mbox{for each}\ u \in {\bf N}.
\end{eqnarray}
Next, for $u,v\in {\bf N} \cup \{0\}$,
we define 
\begin{eqnarray}\label{def-f}
f(u,v):=T_{u}+T_{v}.
\end{eqnarray}
Also, for $k \in {\bf N}$ and $u_1,\dots, u_k,v_1,\dots,v_k \in {\bf N} \cup \{0\}$,
we define
\begin{eqnarray}\label{def-gk}
g_k(u_1,v_1;\dots;u_{k},v_{k})
:={\displaystyle \frac{1}{1-
\sum_{i=2}^{k} (i-1) \cdot F_k^i(f(u_1,v_1),\dots,f(u_{k},v_{k}))
}},\label{def-ghk}
\end{eqnarray}
\begin{eqnarray}\label{def-hk}
h_k(u_1,v_1;\dots;u_{k},v_{k})
:={\displaystyle \frac{\{1+f(u_1,v_1)\} \cdots  \{1+f(u_{k},v_{k})\}}
{1-
\sum_{i=2}^{k} (i-1) \cdot F_k^i(f(u_1,v_1),\dots,f(u_{k},v_{k}))}},
\end{eqnarray}
where $F_k^i(X_1,\dots,X_k)$ is the elementary symmetric polynomial of degree $i$ 
in the $k$ variables $X_1,\dots,X_k$.
In the case $k=n$, for conciseness, we use the following notation:
\begin{eqnarray}\label{def-F^i}
F^i(X_1,\dots,X_n):=F_n^i(X_1,\dots,X_n),
\end{eqnarray}
and
\begin{eqnarray}\label{def-g}
g(u_1,v_1;\dots;u_{n},v_{n})
:=g_n(u_1,v_1;\dots;u_{n},v_{n}),
\end{eqnarray}
\begin{eqnarray}\label{def-h}
h(u_1,v_1;\dots;u_{n},v_{n})
:=h_n(u_1,v_1;\dots;u_{n},v_{n}).
\end{eqnarray}

We have the following lemma.

\begin{lm}\label{lm-growth}
Let $n$ be an integer greater than 1.
For each $k \in \{1,\dots,n\}$,
choose non-negative integers $A_k$ and $B_k$.
Let $K$ be an integer in $\{1,\dots,n\}$,
and let
${\Lambda}_{(A_1,B_1;\dots;A_n,B_n)}$,
${\Lambda}_{(A_1,B_1;\dots;A_n,B_n)}^{K}$ and
${\Lambda}_{(A_1,B_1;\dots;A_n,B_n)}^{KK}$ 
be subsets of $\Sigma^*$ defined as 
\begin{eqnarray*}
\begin{array}{l}
\begin{array}{rl}
{\Lambda}_{(A_1,B_1;\dots;A_n,B_n)}:=\{ \xi \in \Sigma^* \ ;& \
\xi \equiv x_{i_1}^{s_1} \cdots x_{i_\tau}^{s_\tau},\
\mbox{$\xi$ satisfies the condition in (\ref{tau})},\\
&\mbox{ if $x_{i_j}=x_k$, then
$-A_k \leq s_j \leq B_k\ \mbox{and} \ s_j\neq 0\ (1 \leq j \leq \tau)$}
\},\\
\end{array}\\
\begin{array}{rl}
{\Lambda}_{(A_1,B_1;\dots;A_n,B_n)}^{K}:=\{ \xi \in {\Lambda}_{(A_1,B_1;\dots;A_n,B_n)}
 \ ;& \
\xi \equiv x_{i_1}^{s_1} \cdots x_{i_\tau}^{s_\tau},\ 
x_{i_1} \neq x_K
\},\\
\end{array}\\
\begin{array}{rl}
{\Lambda}_{(A_1,B_1;\dots;A_n,B_n)}^{KK}:=\{ \xi \in {\Lambda}_{(A_1,B_1;\dots;A_n,B_n)}\ 
 ;& \
\xi \equiv x_{i_1}^{s_1} \cdots x_{i_\tau}^{s_\tau},\ 
\tau \geq 1,\ x_{i_1} \neq x_K,\ x_{i_{\tau}} \neq x_K
\}.\\
\end{array}
\end{array}
\end{eqnarray*}
If $\tau=0$, we stipulate that $x_{i_1}^{s_1} \cdots x_{i_\tau}^{s_\tau}$ is the null word, $\varepsilon$.
For a subset $\Lambda$ of $\Sigma^*$,
we define the following formal power series:
\begin{eqnarray*}
\Omega_{\Lambda}(t):=
{\displaystyle \sum_{l=0}^{\infty}}
\# \{\xi \in {\Lambda} \ ; \ |\xi|=l \}\  t^l.\\
\end{eqnarray*}
Then we have
\begin{eqnarray*}
\begin{array}{lcl}
\Omega_{{\Lambda}_{(A_1,B_1;\dots;A_n,B_n)}}(t)&=&
h(A_1,B_1;\dots;A_n,B_n),\\
\Omega_{{\Lambda}_{(A_1,B_1;\dots;A_n,B_n)}^{K}}(t)&=&
{\displaystyle \prod_{k=1}^{K-1} }
\{
1+f(A_k,B_k)
\}
\cdot
{\displaystyle \prod_{k=K+1}^{n}}
\{1+f(A_k,B_k)\}
\cdot g(A_1,B_1;\dots;A_n,B_n),\\
\Omega_{{\Lambda}_{(A_1,B_1;\dots;A_n,B_n)}^{KK}}(t)&=&
\{
{\displaystyle \sum_{i=1}^{n-1}
}
i \cdot
F^{i}(f(A_1,B_1),\dots,f(A_{K-1},B_{K-1}),0,\\
&&f(A_{K+1},B_{K+1}),\dots,
f(A_n,B_n))
\}
\cdot g(A_1,B_1;\dots;A_n,B_n)
\end{array}
\end{eqnarray*}
for each $t$ in a sufficiently small neighborhood of the origin, $0$.
(If each element $\xi$ of ${\Lambda}_{(A_1,B_1;\dots;A_n,B_n)}$ is geodesic
and there is no other element $\xi' \in {\Lambda}_{(A_1,B_1;\dots;A_n,B_n)}$ satisfying $\pi(\xi')=\pi(\xi)$,
then we have $\Omega_{{\Lambda}_{(A_1,B_1;\dots;A_n,B_n)}}(t) = {\cal S}_{{\Lambda}_{(A_1,B_1;\dots;A_n,B_n)}}(t)$, by definition.
The same holds for ${\Lambda}_{(A_1,B_1;\dots;A_n,B_n)}^{K}$ and
${\Lambda}_{(A_1,B_1;\dots;A_n,B_n)}^{KK}$.)
\end{lm}

\noindent
{\it Proof.}
First,
we prove the assertion for ${\Lambda}_{(A_1,B_1;\dots;A_n,B_n)}$
by induction with respect to $n$.

The case $n=2$ is easily verified. (For example, see Lemma 5.1 in \cite{Fujii2}.)
Now, assume that the assertion holds in the case $n=k$.
Then we have
\begin{eqnarray*}
\Omega_{{\Lambda}_{(A_1,B_1;\dots;A_{k},B_{k})}}(t)=
h_{k}(A_1,B_1;\dots;A_{k},B_{k}).
\end{eqnarray*}

It is readily seen that
any word $\xi$ in 
${\Lambda}_{(A_1,B_1;\dots;A_{k+1},B_{k+1})}$ is expressed as
\begin{eqnarray}\label{word-Lambda}
\xi \equiv x_{k+1}^{s_0} \cdot
(w_1 x_{k+1}^{s_1} \cdot
w_2 x_{k+1}^{s_2} \cdots
w_{\tau} x_{k+1}^{s_{\tau}})
\cdot
w_{\tau+1},
\end{eqnarray} 
where 
\begin{eqnarray*}
\left\{
\begin{array}{l}
\tau \in {\bf N} \cup \{0\},\\
w_j \in {\Lambda}_{(A_1,B_1;\dots;A_{k},B_{k})} \setminus \{\varepsilon\}
\quad (1 \leq j \leq \tau),\\
w_{\tau+1} \in {\Lambda}_{(A_1,B_1;\dots;A_{k},B_{k})},\\
\end{array}
\right.
\end{eqnarray*}
and
\begin{eqnarray*}
\left\{
\begin{array}{l}
-A_{k+1} \leq s_j \leq B_{k+1}\ (0 \leq j \leq \tau),\\
s_j \neq 0
\ (1 \leq j \leq \tau).\\
\end{array}
\right.
\end{eqnarray*}
If $\tau=0$, we stipulate that $w_1 x_{k+1}^{s_1} \cdot
w_2 x_{k+1}^{s_2} \cdots
w_{\tau} x_{k+1}^{s_{\tau}}$ is the null word, $\varepsilon$.
Then, we obtain 
\begin{eqnarray*}
\begin{array}{l}
\Omega_{{\Lambda}_{(A_1,B_1;\dots;A_{k+1},B_{k+1})}}(t)\\
\begin{array}{r}
=\{1+f(A_{k+1},B_{k+1})\}
\cdot
{\displaystyle \frac{1}{1-\{h_{k}(A_1,B_1;\dots;A_{k},B_{k})-1\}\cdot f(A_{k+1},B_{k+1})}}\\
\\
\cdot
h_{k}(A_1,B_1;\dots;A_{k},B_{k})\\
\end{array}\\
\\
={\displaystyle \frac{
\{1+f(A_1,B_1)\}\cdots
\{1+f(A_{k+1},B_{k+1})\}}
{1-\sum_{i=2}^{k+1}(i-1) \cdot
F_{k+1}^i(f(A_1,B_1),\dots,f(A_{k+1},B_{k+1}))}}\\
\\
=
h_{k+1}(A_1,B_1;\dots;A_{k+1},B_{k+1}).
\end{array}
\end{eqnarray*}
Thus, the assertion holds in the case $n=k+1$.
Therefore, the assertion is proved for ${\Lambda}_{(A_1,B_1;\dots;A_n,B_n)}$.

Next, we prove the assertion for ${\Lambda}_{(A_1,B_1;\dots;A_n,B_n)}^{K}$.
Any word $\xi^K$ in this set is expressed as
\begin{eqnarray}\label{word-Lambda-K}
\xi ^K\equiv 
(w_1 x_{K}^{s_1} \cdot
w_2 x_{K}^{s_2} \cdots
w_{\tau} x_{K}^{s_{\tau}})
\cdot
w_{\tau+1},
\end{eqnarray} 
where 
\begin{eqnarray*}
\left\{
\begin{array}{l}
\tau \in {\bf N} \cup \{0\},\\
w_j \in {\Lambda}_{(A_1,B_1;\dots;A_{K-1},B_{K-1};0,0;
A_{K+1},B_{K+1};\dots;A_{n},B_{n})} \setminus \{\varepsilon\}
\quad (1 \leq j \leq \tau),\\
w_{\tau+1} \in {\Lambda}_{(A_1,B_1;\dots;A_{K-1},B_{K-1};0,0;
A_{K+1},B_{K+1};\dots;A_{n},B_{n})},\\
\end{array}
\right.
\end{eqnarray*}
and
\begin{eqnarray*}
-A_{K} \leq s_j \leq B_{K}\ (1 \leq j \leq \tau),\ 
s_j \neq 0
\ (1 \leq j \leq \tau).
\end{eqnarray*}
Then, using the assertion for ${\Lambda}_{(A_1,B_1;\dots;A_n,B_n)}$,
we obtain 
\begin{eqnarray*}
\begin{array}{l}
\Omega_{{\Lambda}_{(A_1,B_1;\dots;A_{n},B_{n})}^{K}}(t)\\
\begin{array}{r}
=
{\displaystyle \frac{1}{1-\{h(A_1,B_1;\dots;
A_{K-1},B_{K-1};0,0;
A_{K+1},B_{K+1};\dots;
A_{n},B_{n})-1\}\cdot f(A_{K}B_{K})}}\\
\\
\cdot
h(A_1,B_1;\dots;
A_{K-1},B_{K-1};0,0;
A_{K+1},B_{K+1};\dots;
A_{n},B_{n})
\\
\end{array}\\
\\
={\displaystyle \frac{
{\displaystyle \prod_{k=1}^{K-1} }
\{
1+f(A_k,B_k)
\}
\cdot
{\displaystyle \prod_{k=K+1}^{n}}
\{1+f(A_k,B_k)\}
}
{1-{\displaystyle \sum_{i=2}^{n}(i-1)} \cdot
F^i(f(A_1,B_1),\dots,f(A_{n},B_{n}))}}\\
\\
=
{\displaystyle \prod_{k=1}^{K-1} }
\{
1+f(A_k,B_k)
\}
\cdot
{\displaystyle \prod_{k=K+1}^{n}}
\{1+f(A_k,B_k)\}
\cdot g(A_1,B_1;\dots;A_n,B_n).
\end{array}
\end{eqnarray*}

Finally, we prove the assertion for ${\Lambda}_{(A_1,B_1;\dots;A_n,B_n)}^{KK}$.
Any word $\xi^{KK}$ in this set is expressed as
\begin{eqnarray}\label{word-Lambda-K}
\xi ^{KK}\equiv 
(w_1 x_{K}^{s_1} \cdot
w_2 x_{K}^{s_2} \cdots
w_{\tau} x_{K}^{s_{\tau}})
\cdot
w_{\tau+1},
\end{eqnarray} 
where 
\begin{eqnarray*}
\left\{
\begin{array}{l}
\tau \in {\bf N} \cup \{0\},\\
w_j \in {\Lambda}_{(A_1,B_1;\dots;A_{K-1},B_{K-1};0,0;
A_{K+1},B_{K+1};\dots;A_{n},B_{n})} \setminus \{\varepsilon\}
\quad (1 \leq j \leq \tau+1).\\
\end{array}
\right.
\end{eqnarray*}
and
\begin{eqnarray*}
-A_{K} \leq s_j \leq B_{K}\ (1 \leq j \leq \tau),\ 
s_j \neq 0
\ (1 \leq j \leq \tau).
\end{eqnarray*}
Then, using the assertion for ${\Lambda}_{(A_1,B_1;\dots;A_n,B_n)}$,
we obtain 
\begin{eqnarray*}
\begin{array}{l}
\Omega_{{\Lambda}_{(A_1,B_1;\dots;A_{n},B_{n})}^{KK}}(t)\\
\begin{array}{r}
=
{\displaystyle \frac{1}{1-\{h(A_1,B_1;\dots;
A_{K-1},B_{K-1};0,0;
A_{K+1},B_{K+1};\dots;
A_{n},B_{n})-1\}\cdot f(A_{K}B_{K})}}\\
\\
\cdot
\{h(A_1,B_1;\dots;
A_{K-1},B_{K-1};0,0;
A_{K+1},B_{K+1};\dots;
A_{n},B_{n})-1\}
\\
\end{array}\\
={\displaystyle 
\frac{
{\displaystyle 
\sum_{i=1}^{n-1}
}
i \cdot F^{i}
(f(A_1,B_1),\dots,f(A_{K-1},B_{K-1}),0,
f(A_{K+1},B_{K+1}),\dots,f(A_n,B_n))}
{
1-{\displaystyle \sum_{i=2}^{n}(i-1)} \cdot
F^i(f(A_1,B_1),\dots,f(A_{n},B_{n}))}}\\
=
\begin{array}{r}
\{
{\displaystyle 
\sum_{i=1}^{n-1}
}
i \cdot F^{i}
(f(A_1,B_1),\dots,f(A_{K-1},B_{K-1}),0,
f(A_{K+1},B_{K+1}),\dots,f(A_n,B_n))\}\\
\cdot
g(A_1,B_1;\dots;A_n,B_n).\ \Box
\end{array}
\end{array}
\end{eqnarray*}

\subsection{Types {\boldmath $1$, $2$ \mbox{and} $3^+ \cup 3^-$}}

\noindent
[{\bf Types {\boldmath $1$} and {\boldmath $2$}}]
By Propositions \ref{pr-mnf} and \ref{pr-T1,2,3^+,3^-} and Lemma \ref{lm-growth},
the spherical growth series for the set
\begin{eqnarray}\label{def-gamma1'}
\begin{array}{rl}
{\Gamma}'_{1}:=
\{
\nu' \in \Sigma^* \ 
;&
\nu' \equiv 
x_{i_1}^{a_1} \cdots x_{i_\tau}^{a_\tau},\\
&\mbox{$\nu'$ satisfies the conditions in
(\ref{tau}) and (\ref{cond-T1-geod})} \}\\
\end{array}
\end{eqnarray}
has the following rational function expression:
\begin{eqnarray}\label{gamma1'}
h(p_1^-,p_1^+; \cdots; p_n^-, p_n^+).
\end{eqnarray}
(Here, note that if $\tau=0$, the sub-word
$x_{i_1}^{a_1} \cdots
x_{i_\tau}^{a_\tau}$ is regarded as being the null word.)
Moreover, from (\ref{def-gamma1}),
we know that any element $\nu \in {\Gamma}_1$ 
is a product of $\nu' \in {\Gamma}'_1$ 
and $\Delta^c$ ($c \geq 1$).
Hence, with $|\Delta|=p_1$,
we have
\begin{eqnarray}\label{gr-1}
\begin{array}{rcl}
{\cal S}_{\Gamma_{1}}(t)
&=&
h(p_1^-,p_1^+; \cdots; p_n^-, p_n^+)
\cdot
(t^{p_1}+t^{2p_1}+t^{3p_1}+\cdots)\\
&=&
h(p_1^-,p_1^+; \cdots; p_n^-, p_n^+)
\cdot
\frac{t^{p_1}}{1-t^{p_1}}\\
\end{array}
\end{eqnarray}
for each $t$ in a sufficiently small neighborhood of the origin.
Similarly, from (\ref{def-gamma2}),
we have
\begin{eqnarray}\label{gr-2}
{\textstyle {\cal S}_{\Gamma_{2}}(t)
=
h(p_1^-,p_1^+; \cdots; p_n^-, p_n^+)
\cdot
\frac{t^{p_1}}{1-t^{p_1}}}
\end{eqnarray}
for each $t$ in a sufficiently small neighborhood of the origin.

\vspace{0.2cm}

\noindent
[{\bf Type {\boldmath $3^+ \cup 3^-$}}]
From (\ref{def-gamma3^+,3^-}) and (\ref{def-gamma1'}),
we have ${\Gamma}_{3^+}={\Gamma}'_1$.
Hence, from (\ref{gamma1'}),
we have
\begin{eqnarray}\label{gr-3^+}
{\textstyle {\cal S}_{\Gamma_{3^+}}(t)
=h(p_1^-,p_1^+; \cdots; p_n^-, p_n^+)}
\end{eqnarray}
for each $t$ in a sufficiently small neighborhood of the origin.
Similarly, we have
\begin{eqnarray}\label{gr-3^-}
{\textstyle {\cal S}_{\Gamma_{3^-}}(t)
=h(p_1^-,p_1^+; \cdots; p_n^-, p_n^+)}
\end{eqnarray}
for each $t$ in a sufficiently small neighborhood of the origin.
Hence, 
from (\ref{notation-gamma3^+,3^-}), (\ref{gr-3^+}) and (\ref{gr-3^-}), we obtain
\begin{eqnarray}\label{gr-3^+or3^-}
\begin{array}{rcl}
{\cal S}_{\Gamma_{3^+ \cup 3^-}}(t)
&=&{\cal S}_{\Gamma_{3^+}}(t)+{\cal S}_{\Gamma_{3^-}}(t)-{\cal S}_{\Gamma_{3^+ \cap 3^-}}(t)\\
&=&2 {\cal S}_{\Gamma_{3^+}}(t)-{\cal S}_{\Gamma_{3^+ \cap 3^-}}(t).\\
\end{array}
\end{eqnarray}
Also, from (\ref{3^+and3^-}), we have
\begin{eqnarray}\label{gr-3^+and3^-}
{\textstyle {\cal S}_{\Gamma_{3^+ \cap 3^-}}(t)
=h(p_1^-,p_1^-; \cdots; p_n^-, p_n^-)}
\end{eqnarray}
for each $t$ in a sufficiently small neighborhood of the origin.

\subsection{Type {\boldmath $3^0$}}

\subsubsection{Case: {\boldmath $p_k-p_1$} is odd for some {\boldmath $k \in \{ 1, \dots, n\}$}}\label{section-growth:even and odd}

In this section, consider the case $m<n$, i.e.,
the case in which there exists $k \in \{1, \dots,n\}$ such that $p_k-p_1$ is odd.
\vspace{0.3cm}

From Theorem \ref{th-ugr} and the partition (\ref{partition-3^0}),
we find that the spherical growth series for ${\Gamma}_{3^0}$ 
is obtained as the following summation of the growth series:
\begin{eqnarray}\label{sum-3^0}
\begin{array}{rl}
{\cal S}_{\Gamma_{3^0}}(z)
=&
{\displaystyle \sum_{N=0}^{p_1-1}
\sum_{M=0}^{n-m-1}
{\cal S}_{\Gamma_{3^0}^{(N,2M+1)}}(t)
+\sum_{N=0}^{p_1-2}
\sum_{M=0}^{n-m-1}
{\cal S}_{\Gamma_{3^0}^{(N,2M+2)}}(t)
+\sum_{M=0}^{n-m-2}
{\cal S}_{\Gamma_{3^0}^{(p_1-1,2M+2)}}(t)}\\
&{\displaystyle +\sum_{N=0}^{p_1-2}
\sum_{M=n-m}^{n-1}
{\cal S}_{\Gamma_{3^0}^{(N,2M+1)}}(t)
+\sum_{N=0}^{p_1-2}
\sum_{M=n-m}^{n-1}
{\cal S}_{\Gamma_{3^0}^{(N,2M+2)}}(t)}.
\end{array}
\end{eqnarray}

Now, we give rational function expressions for all the terms on
the right-hand side of (\ref{sum-3^0}).

\vspace{0.2cm}

\noindent
{\bf [Case {\boldmath $(N,2M+1)$} : {\boldmath $0 \leq N \leq p_1-1$},
{\boldmath $0 \leq M \leq n-m-1$}]}

\noindent
By Proposition \ref{pr-sc-T3} and Lemma \ref{lm-growth},
the spherical growth series for the three sets
\begin{eqnarray*}
\begin{array}{l}
\Gamma_{\xi^{(1)}}:= \{\xi^{(1)} \in \Sigma^* \ ; \ 
\mbox{$\xi^{(1)}$ is a word given in (\ref{1-2}) with
the condition in (\ref{1-3})} \},\\
\Gamma_{\xi^{(2)}}:=\{\xi^{(2)} \in \Sigma^* \ ; \ 
\mbox{$\xi^{(2)}$ is a word given in (\ref{1-2}) with
the condition in (\ref{1-4})} \},\\
\Gamma_{\xi^{(3)}}:=\{\xi^{(3)} \in \Sigma^* \ ; \ 
\mbox{$\xi^{(3)}$ is a word given in (\ref{1-2}) with
the condition in (\ref{1-5})} \}
\end{array}
\end{eqnarray*}
have the following rational function expressions:
\begin{eqnarray*}
\begin{array}{rcl}
&{\cal S}_{\Gamma_{\xi^{(1)}}}(t)=&
\prod_{k=1}^m 
\{
1+f(q_k^- +N,q_k^+ -N)
\}
\cdot
\prod_{k=m+1}^{n-M-1}
\{1+f(r_k^- +N,r_k^+ -N)\}
\\
&&\cdot
\prod_{k=n-M+1}^{n}
\{
1+f(r_k^- +N+1,r_k^+ -N-1)
\}\\
&&\cdot
g(q_1^- +N,q_1^+ -N; \dots; q_m^- +N,q_m^+ -N;\\
&&
r_{m+1}^- +N,r_{m+1}^+ -N; \dots; r_{n-M-1}^- +N,r_{n-M-1}^+ -N
;\\
&&
r_{n-M}^- +N,r_{n-M}^+ -N-1;\\
&&
r_{n-M+1}^- +N+1,r_{n-M+1}^+ -N-1;
\dots;
r_{n}^- +N+1,r_{n}^+ -N-1
),\\
&{\cal S}_{\Gamma_{\xi^{(2)}}}(t)=&
\{
\sum_{i=1}^{n-1} 
i \cdot F^i(f(q_1^- +N,q_1^+ -N),\dots,f(q_m^- +N,q_m^+ -N),\\
&&
f(r_{m+1}^- +N,r_{m+1}^+ -N),\dots,f(r_{n-M-1}^- +N,r_{n-M-1}^+ -N),\\
&&
0,f(r_{n-M+1}^- +N+1,r_{n-M+1}^+ -N-1),\dots,f(r_{n}^- +N+1,r_{n}^+ -N-1))
\}
\\
&&\cdot
g(q_1^- +N,q_1^+ -N; \dots; q_m^- +N,q_m^+ -N;\\
&&
r_{m+1}^- +N,r_{m+1}^+ -N; \dots; r_{n-M-1}^- +N,r_{n-M-1}^+ -N
;\\
&&
r_{n-M}^- +N+1,r_{n-M}^+ -N-1;\\
&&
r_{n-M+1}^- +N+1,r_{n-M+1}^+ -N-1;
\dots;
r_{n}^- +N+1,r_{n}^+ -N-1
),\\
&{\cal S}_{\Gamma_{\xi^{(3)}}}(t)=&
\prod_{k=1}^m 
\{
1+f(q_k^- +N,q_k^+ -N)
\}
\cdot
\prod_{k=m+1}^{n-M-1}
\{1+f(r_k^- +N,r_k^+ -N)\}
\\
&&\cdot
\prod_{k=n-M+1}^{n}
\{
1+f(r_k^- +N+1,r_k^+ -N-1)
\}\\
&&\cdot
g(q_1^- +N,q_1^+ -N; \dots; q_m^- +N,q_m^+ -N;\\
&&
r_{m+1}^- +N,r_{m+1}^+ -N; \dots; r_{n-M-1}^- +N,r_{n-M-1}^+ -N;\\
&&
r_{n-M}^- +N,r_{n-M}^+ -N;\\
&&
r_{n-M+1}^- +N+1,r_{n-M+1}^+ -N-1;
\dots;
r_{n}^- +N+1,r_{n}^+ -N-1
).\\
\end{array}
\end{eqnarray*}
Hence, using the equality
$(r_{n-M}^- +N+1)+(r_{n-M}^+ -N)=r_{n-M}$,
we obtain
\begin{eqnarray}\label{gr-3^0,1}
{\cal S}_{\Gamma_{3^0}^{(N,2M+1)}}(t)
=t^{r_{n-M}}\cdot 
{\cal S}_{\Gamma_{\xi^{(1)}}}(t)
\cdot {\cal S}_{\Gamma_{\xi^{(2)}}}(t)
\cdot {\cal S}_{\Gamma_{\xi^{(3)}}}(t)
\end{eqnarray}
for each $t$ in a sufficiently small neighborhood of the origin.

\noindent
{\bf [Case {\boldmath $(N,2M+2)$} :
{\boldmath $0 \leq N \leq p_1-2$} and {\boldmath $0 \leq M \leq n-m-1$},
or {\boldmath $N=p_1-1$} and {\boldmath $0 \leq M \leq n-m-2$}]}

\noindent
We define the following four sets:
\begin{eqnarray*}
&&\begin{array}{rcl}
{\Phi}_1:=\{\xi \in \Sigma^* ;\
\xi \equiv
x_{i_1}^{a_1}  x_{i_2}^{a_2} \cdots
x_{i_{\kappa}}^{a_\kappa},\
\mbox{$\kappa$ and $a_j$ satisfy 
the condition in (\ref{2-2})} \},
\end{array}\\
&&\begin{array}{rcl}
{\Phi}_2:=\{\xi \in \Sigma^* ;\
\xi \equiv
x_{i_1}^{a_1}  x_{i_2}^{a_2} \cdots
x_{i_{\kappa}}^{a_\kappa},\ 
\mbox{$\kappa$ and $a_j$ satisfy 
the condition in (\ref{2-5})} \},
\end{array}\\
&&\begin{array}{rcl}
{\Phi}_3:=\{\xi \in \Sigma^* ;\
\xi \equiv
x_{i_1}^{a_1}  x_{i_2}^{a_2} \cdots
x_{i_{\kappa}}^{a_\kappa},\
\mbox{$\kappa$ and $a_j$ satisfy 
the condition in (\ref{2-6})} \},
\end{array}\\
&&\begin{array}{rcl}
{\Phi}_4:=\{\xi \in \Sigma^* ;\
\xi \equiv
x_{i_1}^{a_1}  x_{i_2}^{a_2} \cdots
x_{i_{\kappa}}^{a_\kappa},\
\mbox{$\kappa$ and $a_j$ satisfy 
the condition in (\ref{2-7})} \},
\end{array}
\end{eqnarray*}
where the conditions (\ref{2-5})--(\ref{2-7}) 
referred to are the following:
\begin{eqnarray}
&&
\bullet 
\left\{
\begin{array}{rl}
\bullet& \kappa \in {\bf N} \cup \{0\},\\
\bullet&
x_{i_j}=z_k\ \mbox{with}\ n-M+1 \leq k \leq n\\
&\Longrightarrow 
\begin{array}{l}
-(r_k^- +N+1) \leq a_j \leq  
r_k^+ -N-1,\ 
a_j \neq 0,\\
\end{array}
\\
\bullet&
x_{i_j}=z_k\ \mbox{with}\ k = n-M \\
&\Longrightarrow 
\begin{array}{l}
-(r_{n-M}^- +N) \leq a_j \leq  
r_{n-M}^+ -N-1,\ 
a_j \neq 0,\\
\end{array}
\\
\bullet&
x_{i_j}=z_k\ \mbox{with}\ m+1 \leq k \leq n-M-1\\
&\Longrightarrow 
\begin{array}{l}
-(r_k^- +N) \leq a_j \leq  
r_k^+ -N,\ 
a_j \neq 0,\\
\end{array}
\\
\bullet&
x_{i_j}=y_k\ \mbox{with}\ 1 \leq k \leq m\\
&\Longrightarrow 
\begin{array}{l}
-(q_k^- +N) \leq a_j \leq  
q_k^+ -N,\
a_j \neq 0.\\
\end{array}
\\
\end{array}
\right.
\label{2-5}\\
&&\bullet
\left\{
\begin{array}{rl}
\bullet& \kappa \in {\bf N} \cup \{0\},\\
\bullet&
x_{i_j}=z_k\ \mbox{with}\ n-M \leq k \leq n\\
&\Longrightarrow 
\begin{array}{l}
-(r_k^- +N+1) \leq a_j \leq  
r_k^-,\ 
a_j \neq 0,\\
\end{array}
\\
\bullet&
x_{i_j}=z_k\ \mbox{with}\ m+1 \leq k \leq n-M-1\\
&\Longrightarrow 
\begin{array}{l}
-(r_k^- +N) \leq a_j \leq  
r_k^-,\ 
a_j \neq 0,\\
\end{array}
\\
\bullet&
x_{i_j}=y_k\ \mbox{with}\ 1 \leq k \leq m\\
&\Longrightarrow 
\begin{array}{l}
-(q_k^- +N) \leq a_j \leq  
q_k^-,\
a_j \neq 0.\\
\end{array}
\\
\end{array}
\right.
\label{2-6}\\
&&\bullet 
\left\{
\begin{array}{rl}
\bullet& \kappa \in {\bf N} \cup \{0\},\\
\bullet&
x_{i_j}=z_k\ \mbox{with}\ n-M+1 \leq k \leq n\\
&\Longrightarrow 
\begin{array}{l}
-(r_k^- +N+1) \leq a_j \leq  
r_k^-,\ 
a_j \neq 0,\\
\end{array}
\\
\bullet&
x_{i_j}=z_k\ \mbox{with}\ m+1 \leq k \leq n-M\\
&\Longrightarrow 
\begin{array}{l}
-(r_k^- +N) \leq a_j \leq  
r_k^-,\ 
a_j \neq 0,\\
\end{array}
\\
\bullet&
x_{i_j}=y_k\ \mbox{with}\ 1 \leq k \leq m\\
&\Longrightarrow 
\begin{array}{l}
-(q_k^- +N) \leq a_j \leq  
q_k^-,\
a_j \neq 0.\\
\end{array}
\\
\end{array}
\right.
\label{2-7}
\end{eqnarray}
Then, we have
\begin{eqnarray}\label{gamma-N2}
{\Gamma}_{3^0}^{(N,2M+2)}
=({\Phi}_1 \setminus {\Phi}_2 )
\setminus ({\Phi}_3 \setminus {\Phi}_4 ).
\end{eqnarray}
By Proposition \ref{pr-sc-T3} and Lemma \ref{lm-growth},
the growth series for the sets
${\Phi}_1,{\Phi}_2,{\Phi}_3$ and ${\Phi}_4$
have the following rational function expressions:
\begin{eqnarray*}
\begin{array}{rcl}
{\cal S}_{\Phi_1}(t)&=&
h(q_1^- +N,q_1^+ -N; \dots; q_m^- +N,q_m^+ -N;\\
&&
r_{m+1}^- +N,r_{m+1}^+ -N; \dots; r_{n-M-1}^- +N,r_{n-M-1}^+ -N
;\\
&&
r_{n-M}^- +N+1,r_{n-M}^+ -N-1;\\
&&
r_{n-M+1}^- +N+1,r_{n-M+1}^+ -N-1;
\dots;
r_{n}^- +N+1,r_{n}^+ -N-1
),\\
{\cal S}_{\Phi_2}(t)&=&
h(q_1^- +N,q_1^+ -N; \dots; q_m^- +N,q_m^+ -N;\\
&&
r_{m+1}^- +N,r_{m+1}^+ -N; \dots; r_{n-M-1}^- +N,r_{n-M-1}^+ -N;\\
&&
r_{n-M}^- +N,r_{n-M}^+ -N-1;\\
&&
r_{n-M+1}^- +N+1,r_{n-M+1}^+ -N-1;
\dots;
r_{n}^- +N+1,r_{n}^+ -N-1
),\\
{\cal S}_{\Phi_3}(t)&=&
h(q_1^- +N,q_1^-; \dots; q_m^- +N,q_m^-;\\
&&
r_{m+1}^- +N,r_{m+1}^-; \dots; r_{n-M-1}^- +N,r_{n-M-1}^-;\\
&&
r_{n-M}^- +N+1,r_{n-M}^-;\\
&&
r_{n-M+1}^- +N+1,r_{n-M+1}^-;
\dots;
r_{n}^- +N+1,r_{n}^-
),\\
{\cal S}_{\Phi_4}(t)&=&
h(q_1^- +N,q_1^-; \dots; q_m^- +N,q_m^-;\\
&&
r_{m+1}^- +N,r_{m+1}^-; \dots; r_{n-M-1}^- +N,r_{n-M-1}^-;\\
&&
r_{n-M}^- +N,r_{n-M}^-;\\
&&
r_{n-M+1}^- +N+1,r_{n-M+1}^-;
\dots;
r_{n}^- +N+1,r_{n}^-
).\\
\end{array}
\end{eqnarray*}
Hence, from (\ref{gamma-N2}), we obtain
\begin{eqnarray}\label{gr-3^0,2}
{\cal S}_{\Gamma_{3^0}^{(N,2M+2)}}(t)
=
\{{\cal S}_{\Phi_1}(t)-{\cal S}_{\Phi_2}(t)\}-\{{\cal S}_{\Phi_3}(t)-{\cal S}_{\Phi_4}(t)\}
\end{eqnarray}
for each $t$ in a sufficiently small neighborhood of the origin.


\noindent
{\bf [Case {\boldmath $(N,2M+1)$} : {\boldmath $0 \leq N \leq p_1-2$},
{\boldmath $n-m \leq M \leq n-1$}]}

\noindent
By Proposition \ref{pr-sc-T3} and Lemma \ref{lm-growth},
the spherical growth series for the three sets
\begin{eqnarray*}
\begin{array}{l}
\Gamma_{\xi^{(4)}}:= \{\xi^{(4)} \in \Sigma^* \ ; \ 
\mbox{$\xi^{(4)}$ is a word given in (\ref{3-2}) with
the condition in (\ref{3-3})} \},\\
\Gamma_{\xi^{(5)}}:=\{\xi^{(5)} \in \Sigma^* \ ; \ 
\mbox{$\xi^{(5)}$ is a word given in (\ref{3-2}) with
the condition in (\ref{3-4})} \},\\
\Gamma_{\xi^{(6)}}:=\{\xi^{(6)} \in \Sigma^* \ ; \ 
\mbox{$\xi^{(6)}$ is a word given in (\ref{3-2}) with
the condition in (\ref{3-5})} \}
\end{array}
\end{eqnarray*}
have the following rational function expressions:
\begin{eqnarray*}
\begin{array}{rcl}
&{\cal S}_{\Gamma_{\xi^{(4)}}}(t)=&
\prod_{k=1}^{n-M-1} 
\{
1+f(q_k^- +N,q_k^+ -N)
\}
\cdot
\prod_{k=n-M+1}^{m}
\{1+f(q_k^- +N+1,q_k^+ -N-1)\}
\\
&&\cdot
\prod_{k=m+1}^{n}
\{
1+f(r_k^- +N+1,r_k^+ -N-1)
\}\\
&&\cdot
g(q_1^- +N,q_1^+ -N; \dots; q_{n-M-1}^- +N,q_{n-M-1}^+ -N;\\
&&
q_{n-M}^- +N,q_{n-M}^+ -N-1;\\
&&
q_{n-M+1}^- +N+1,q_{n-M+1}^+ -N-1; \dots; q_{m}^- +N+1,q_{m}^+ -N-1;\\
&&
r_{m+1}^- +N+1,r_{m+1}^+ -N-1; \dots; r_{n}^- +N+1,r_{n}^+ -N-1
),\\
&{\cal S}_{\Gamma_{\xi^{(5)}}}(t)=&
\{
\sum_{i=1}^{n-1} 
i \cdot F^i(f(q_1^- +N,q_1^+ -N),\dots,f(q_{n-M-1}^- +N,q_{n-M-1}^+ -N),\\
&&
0,
f(q_{n-M+1}^- +N+1,q_{n-M+1}^+ -N-1),\dots,f(q_{m}^- +N+1,q_{m}^+ -N-1),\\
&&
f(r_{m+1}^- +N+1,r_{m+1}^+ -N-1),\dots,f(r_{n}^- +N+1,r_{n}^+ -N-1))
\}
\\
&&\cdot
g(q_1^- +N,q_1^+ -N; \dots; q_{n-M-1}^- +N,q_{n-M-1}^+ -N;\\
&&
q_{n-M}^- +N+1,q_{n-M}^+ -N-1;\\
&&
q_{n-M+1}^- +N+1,q_{n-M+1}^+ -N-1; \dots; q_{m}^- +N+1,q_{m}^+ -N-1;\\
&&
r_{m+1}^- +N+1,r_{m+1}^+ -N-1; \dots; r_{n}^- +N+1,r_{n}^+ -N-1
),\\
&{\cal S}_{\Gamma_{\xi^{(6)}}}(t)=&
\prod_{k=1}^{n-M-1} 
\{
1+f(q_k^- +N,q_k^+ -N)
\}
\cdot
\prod_{k=n-M+1}^{m}
\{1+f(q_k^- +N+1,q_k^+ -N-1)\}
\\
&&\cdot
\prod_{k=m+1}^{n}
\{
1+f(r_k^- +N+1,r_k^+ -N-1)
\}\\
&&\cdot
g(q_1^- +N,q_1^+ -N; \dots; q_{n-M-1}^- +N,q_{n-M-1}^+ -N;\\
&&
q_{n-M}^- +N,q_{n-M}^+ -N;\\
&&
q_{n-M+1}^- +N+1,q_{n-M+1}^+ -N-1; \dots; q_{m}^- +N+1,q_{m}^+ -N-1;\\
&&
r_{m+1}^- +N+1,r_{m+1}^+ -N-1; \dots; r_{n}^- +N+1,r_{n}^+ -N-1
).\\
\end{array}
\end{eqnarray*}
Hence, using the equality
$(q_{n-M}^- +N+1)+(q_{n-M}^+ -N)=q_{n-M}$,
we obtain
\begin{eqnarray}\label{gr-3^0,3}
{\cal S}_{\Gamma_{3^0}^{(N,2M+1)}}(t)
=t^{q_{n-M}}\cdot 
{\cal S}_{\Gamma_{\xi^{(4)}}}(t)
\cdot {\cal S}_{\Gamma_{\xi^{(5)}}}(t)
\cdot {\cal S}_{\Gamma_{\xi^{(6)}}}(t)
\end{eqnarray}
for each $t$ in a sufficiently small neighborhood of the origin.

\vspace{0.2cm}

\noindent
{\bf [Case {\boldmath $(N,2M+2)$} :
{\boldmath $0 \leq N \leq p_1-2$}, {\boldmath $n-m \leq M \leq n-1$}]}

\noindent
We define the following four sets:
\begin{eqnarray*}
&&\begin{array}{rcl}
{\Psi}_1:=\{\xi \in \Sigma^* ;\
\xi \equiv
x_{i_1}^{a_1}  x_{i_2}^{a_2} \cdots
x_{i_{\kappa}}^{a_\kappa},\
\mbox{$\kappa$ and $a_j$ satisfy 
the condition in (\ref{4-2})} \},
\end{array}\\
&&\begin{array}{rcl}
{\Psi}_2:=\{\xi \in \Sigma^* ;\
\xi \equiv
x_{i_1}^{a_1}  x_{i_2}^{a_2} \cdots
x_{i_{\kappa}}^{a_\kappa},\ 
\mbox{$\kappa$ and $a_j$ satisfy 
the condition in (\ref{4-5})} \},
\end{array}\\
&&\begin{array}{rcl}
{\Psi}_3:=\{\xi \in \Sigma^* ;\
\xi \equiv
x_{i_1}^{a_1}  x_{i_2}^{a_2} \cdots
x_{i_{\kappa}}^{a_\kappa},\
\mbox{$\kappa$ and $a_j$ satisfy 
the condition in (\ref{4-6})} \},
\end{array}\\
&&\begin{array}{rcl}
{\Psi}_4:=\{\xi \in \Sigma^* ;\
\xi \equiv
x_{i_1}^{a_1}  x_{i_2}^{a_2} \cdots
x_{i_{\kappa}}^{a_\kappa},\
\mbox{$\kappa$ and $a_j$ satisfy 
the condition in (\ref{4-7})} \},
\end{array}
\end{eqnarray*}
where the conditions (\ref{4-5})--(\ref{4-7}) 
referred to are the following:
\begin{eqnarray}
&&
\bullet 
\left\{
\begin{array}{rl}
\bullet& \kappa \in {\bf N} \cup \{0\},\\
\bullet&
x_{i_j}=z_k\ \mbox{with}\ m+1 \leq k \leq n\\
&\Longrightarrow 
\begin{array}{l}
-(r_k^- +N+1) \leq a_j \leq  
r_k^+ -N-1,\ 
a_j \neq 0,\\
\end{array}
\\
\bullet&
x_{i_j}=y_k\ \mbox{with}\ n-M+1 \leq k \leq m \\
&\Longrightarrow 
\begin{array}{l}
-(q_k^- +N+1) \leq a_j \leq  
q_k^+ -N-1,\ 
a_j \neq 0,\\
\end{array}
\\
\bullet&
x_{i_j}=y_k\ \mbox{with}\ k = n-M\\
&\Longrightarrow 
\begin{array}{l}
-(q_{n-M}^- +N) \leq a_j \leq  
q_{n-M}^+ -N-1,\ 
a_j \neq 0,\\
\end{array}
\\
\bullet&
x_{i_j}=y_k\ \mbox{with}\ 1 \leq k \leq n-M-1\\
&\Longrightarrow 
\begin{array}{l}
-(q_k^- +N) \leq a_j \leq  
q_k^+ -N,\
a_j \neq 0.\\
\end{array}
\\
\end{array}
\right.
\label{4-5}\\
&&\bullet
\left\{
\begin{array}{rl}
\bullet& \kappa \in {\bf N} \cup \{0\},\\
\bullet&
x_{i_j}=z_k\ \mbox{with}\ m+1 \leq k \leq n\\
&\Longrightarrow 
\begin{array}{l}
-(r_k^- +N+1) \leq a_j \leq  
r_k^-,\ 
a_j \neq 0,\\
\end{array}
\\
\bullet&
x_{i_j}=y_k\ \mbox{with}\ n-M \leq k \leq m\\
&\Longrightarrow 
\begin{array}{l}
-(q_k^- +N+1) \leq a_j \leq  
q_k^-,\ 
a_j \neq 0,\\
\end{array}
\\
\bullet&
x_{i_j}=y_k\ \mbox{with}\ 1 \leq k \leq n-M-1\\
&\Longrightarrow 
\begin{array}{l}
-(q_k^- +N) \leq a_j \leq  
q_k^-,\
a_j \neq 0.\\
\end{array}
\\
\end{array}
\right.
\label{4-6}\\
&&\bullet 
\left\{
\begin{array}{rl}
\bullet& \kappa \in {\bf N} \cup \{0\},\\
\bullet&
x_{i_j}=z_k\ \mbox{with}\ m+1 \leq k \leq n\\
&\Longrightarrow 
\begin{array}{l}
-(r_k^- +N+1) \leq a_j \leq  
r_k^-,\ 
a_j \neq 0,\\
\end{array}
\\
\bullet&
x_{i_j}=y_k\ \mbox{with}\ n-M+1 \leq k \leq m\\
&\Longrightarrow 
\begin{array}{l}
-(q_k^- +N+1) \leq a_j \leq  
q_k^-,\ 
a_j \neq 0,\\
\end{array}
\\
\bullet&
x_{i_j}=y_k\ \mbox{with}\ 1 \leq k \leq n-M\\
&\Longrightarrow 
\begin{array}{l}
-(q_k^- +N) \leq a_j \leq  
q_k^-,\
a_j \neq 0.\\
\end{array}
\\
\end{array}
\right.
\label{4-7}
\end{eqnarray}
Then, we have
\begin{eqnarray}\label{gamma-N4}
{\Gamma}_{3^0}^{(N,2M+2)}
=({\Psi}_1 \setminus {\Psi}_2 )
\setminus ({\Psi}_3 \setminus {\Psi}_4 ).
\end{eqnarray}
By Proposition \ref{pr-sc-T3} and Lemma \ref{lm-growth},
the growth series for the sets
${\Psi}_1,{\Psi}_2,{\Psi}_3$ and ${\Psi}_4$
have the following rational function expressions:
\begin{eqnarray*}
\begin{array}{rcl}
{\cal S}_{\Psi_1}(t)&=&
h(q_1^- +N,q_1^+ -N; \dots; q_{n-M-1}^- +N,q_{n-M-1}^+ -N;\\
&&
q_{n-M}^- +N+1,q_{n-M}^+ -N-1;\\
&&
q_{n-M+1}^- +N+1,q_{n-M+1}^+ -N-1; \dots; q_{m}^- +N+1,q_{m}^+ -N-1;\\
&&
r_{m+1}^- +N+1,r_{m+1}^+ -N-1;
\dots;
r_{n}^- +N+1,r_{n}^+ -N-1
),\\
{\cal S}_{\Psi_2}(t)&=&
h(q_1^- +N,q_1^+ -N; \dots; q_{n-M-1}^- +N,q_{n-M-1}^+ -N;\\
&&
q_{n-M}^- +N,q_{n-M}^+ -N-1;\\
&&
q_{n-M+1}^- +N+1,q_{n-M+1}^+ -N-1; \dots; q_{m}^- +N+1,q_{m}^+ -N-1;\\
&&
r_{m+1}^- +N+1,r_{m+1}^+ -N-1;
\dots;
r_{n}^- +N+1,r_{n}^+ -N-1
),\\
{\cal S}_{\Psi_3}(t)&=&
h(q_1^- +N,q_1^-; \dots; q_{n-M-1}^- +N,q_{n-M-1}^-;\\
&&
q_{n-M}^- +N+1,q_{n-M}^-;\\
&&
q_{n-M+1}^- +N+1,q_{n-M+1}^-; \dots; q_{m}^- +N+1,q_{m}^-;\\
&&
r_{m+1}^- +N+1,r_{m+1}^-;
\dots;
r_{n}^- +N+1,r_{n}^-
),\\
{\cal S}_{\Psi_4}(t)&=&
h(q_1^- +N,q_1^-; \dots; q_{n-M-1}^- +N,q_{n-M-1}^-;\\
&&
q_{n-M}^- +N,q_{n-M}^-;\\
&&
q_{n-M+1}^- +N+1,q_{n-M+1}^-; \dots; q_{m}^- +N+1,q_{m}^-;\\
&&
r_{m+1}^- +N+1,r_{m+1}^-;
\dots;
r_{n}^- +N+1,r_{n}^-
).\\
\end{array}
\end{eqnarray*}
Hence, from (\ref{gamma-N4}), we obtain
\begin{eqnarray}\label{gr-3^0,4}
{\cal S}_{\Gamma_{3^0}^{(N,2M+2)}}(t)
=
\{{\cal S}_{\Psi_1}(t)-{\cal S}_{\Psi_2}(t)\}-\{{\cal S}_{\Psi_3}(t)-{\cal S}_{\Psi_4}(t)\}
\end{eqnarray}
for each $t$ in a sufficiently small neighborhood of the origin.

\subsubsection{Case: {\boldmath $p_k-p_1$} is even for every {\boldmath $k \in \{ 1, \dots, n\}$}}\label{section-growth:all even}

In this section, we consider the case $m=n$, i.e.,
that in which $p_k-p_1$ is even for every $k \in \{1, \dots,n\}$.
\vspace{0.3cm}

From Theorem \ref{th-ugr} and the partition (\ref{epartition-3^0}),
we find that the spherical growth series for ${\Gamma}_{3^0}$ 
is obtained as the following summation of the growth series:
\begin{eqnarray}\label{esum-3^0}
{\cal S}_{\Gamma_{3^0}}(z)
=
{\displaystyle \sum_{N=0}^{p_1-2}
\sum_{M=0}^{n-1}
{\cal S}_{\bar{\Gamma}_{3^0}^{(N,2M+1)}}(t)
+\sum_{N=0}^{p_1-3}
\sum_{M=0}^{n-1}
{\cal S}_{\bar{\Gamma}_{3^0}^{(N,2M+2)}}(t)
+\sum_{M=0}^{n-2}
{\cal S}_{\bar{\Gamma}_{3^0}^{(p_1-2,2M+2)}}(t)}.
\end{eqnarray}

Now, we give rational function expressions 
for all the terms on 
the right-hand side of (\ref{esum-3^0}).

\noindent
{\bf [Case {\boldmath $(N,2M+1)$} : {\boldmath $0 \leq N \leq p_1-2$},
{\boldmath $0 \leq M \leq n-1$}]}

\noindent
By Proposition \ref{pr-sc-T3} and Lemma \ref{lm-growth},
for each $i \in \{4,5,6\}$,
the spherical growth series of the set
\begin{eqnarray*}
\begin{array}{l}
\bar{\Gamma}_{\bar\xi^{(i)}}:= \{\bar\xi^{(i)} \in \Sigma^* \ ; \ 
\mbox{$\bar\xi^{(i)}$ is a word given in (\ref{e3-1})} \}\\
\end{array}
\end{eqnarray*}
satisfies the following:
\begin{eqnarray*}
{\cal S}_{\bar{\Gamma}_{\bar\xi^{(i)}}}(t)=
{\cal S}_{\Gamma_{\xi^{(i)}}}(t)\quad \mbox{with}\quad m=n.
\end{eqnarray*}
Hence, 
we obtain
\begin{eqnarray}\label{gr-3^0,e1}
{\cal S}_{\bar{\Gamma}_{3^0}^{(N,2M+1)}}(t)
={\cal S}_{\Gamma_{3^0}^{(N,2M+1)}}(t)\quad \mbox{with}\quad m=n
\end{eqnarray}
for each $t$ in a sufficiently small neighborhood of the origin.

\vspace{0.2cm}

\noindent
{\bf [Case {\boldmath $(N,2M+2)$} :
{\boldmath $0 \leq N \leq p_1-3$} and {\boldmath $0 \leq M \leq n-1$},
or 
{\boldmath $N=p_1-2$} and {\boldmath $0 \leq M \leq n-2$}]}

\noindent
Here, as in the previous case, we obtain 
\begin{eqnarray}\label{gr-3^0,4}
{\cal S}_{\bar{\Gamma}_{3^0}^{(N,2M+2)}}(t)
={\cal S}_{\Gamma_{3^0}^{(N,2M+2)}}(t)\quad \mbox{with}\quad m=n
\end{eqnarray}
for each $t$ in a sufficiently small neighborhood of the origin.

\subsection{Growth series for $G(p_1,\dots,p_n)$}

From the partition (\ref{Gamma}), we have
\begin{eqnarray}\label{sum-1,2,3}
{\cal S}_{G(p_1,\dots,p_n)}(t)
= {\cal S}_{\Gamma_1}(t)+{\cal S}_{\Gamma_2}(t)+{\cal S}_{\Gamma_{3^+ \cup 3^-}}(t)
+{\cal S}_{\Gamma_{3^0}}(t),
\end{eqnarray}
and from (\ref{gr-1})--(\ref{gr-3^+}),
we have
\begin{eqnarray}\label{gr-12-2}
{\cal S}_{\Gamma_1}(t)={\cal S}_{\Gamma_2}(t)={\cal S}_{\Gamma_{3^+}}(t) \cdot 
\frac{t^{p_1}}{1-t^{p_1}}.
\end{eqnarray}

Combining (\ref{sum-1,2,3}), (\ref{gr-12-2}), (\ref{gr-3^+or3^-}), (\ref{sum-3^0})
and (\ref{esum-3^0}),
we obtain a rational function expression for the spherical growth series
of $G(p_1,\dots,p_n)$ as follows:

\begin{thm}\label{th-growth}
Let $n$ be an integer greater than 1, let
$p_1,p_2, \dots,p_n$ be integers satisfying $2 \leq p_1 \leq p_2 \leq \cdots \leq p_n$,
and let $G(p_1,\dots,p_n)$ be the group presented as in (\ref{group-presentation}).
Let $m$ be the number of $p_k$'s such that $p_k-p_1$ is even.
Then the spherical growth series of the group $G(p_1,\dots,p_n)$ 
with respect to the generating set 
$\{x_1,\dots,x_n,x_1^{-1},\dots,x_n^{-1} \}$
has the following
rational function expression:

\begin{eqnarray}
{\cal S}_{G(p_1,\dots,p_n)}(t)
&=&2 {\cal S}_{\Gamma_{3^+}}(t) \cdot 
{\displaystyle \frac{1}{1-t^{p_1}}}-{\cal S}_{\Gamma_{3^+ \cap 3^-}}(t)
+{\cal S}_{\Gamma_{3^0}}(t).
\end{eqnarray}
Here, we note the following:

\noindent
$\bullet$
${\cal S}_{\Gamma_{3^+}}(t)$ and ${\cal S}_{\Gamma_{3^+ \cap 3^-}}(t)$
are given in (\ref{gr-3^+}) and (\ref{gr-3^+and3^-}), repsectively.

\noindent
$\bullet$
If $m<n$, then the last term, ${\cal S}_{\Gamma_{3^0}}(t)$, is given in
(\ref{sum-3^0}), and all the terms appearing in (\ref{sum-3^0}) are given in 
Section \ref{section-growth:even and odd}.

\noindent
$\bullet$
If $m=n$, then ${\cal S}_{\Gamma_{3^0}}(t)$ is given in
(\ref{esum-3^0}), and all the terms appearing in (\ref{esum-3^0}) are given in 
Sections \ref{section-growth:even and odd} and \ref{section-growth:all even}.

\noindent
$\bullet$
The rational functions
$T_u$, $f(u,v)$, $g(u_1,v_1;\dots;u_n,v_n)$
and $h(u_1,v_1;\dots;u_n,v_n)$ are given in 
(\ref{def-T}), (\ref{def-f}), (\ref{def-g}) and  (\ref{def-h}), respectively.
\end{thm}

\vspace{0.1cm}

\begin{ex}\label{ex-growth}
\begin{eqnarray*}
\begin{array}{l}

{\cal S}_{G(2,2,2)}(t)
=
\frac{
(1 + t) (-1 + 2 t^2)}
{(-1 + t) (-1 + 2 t)^2}.\\
\\
{\cal S}_{G(2,2,3)}(t)
=\frac{
(1 + t) (1 + t - 3 t^2 - 15 t^3 - 10 t^4 + 30 t^5 + 28 t^6 + 
   16 t^7)}
{(-1 + t) (-1 + t + 4 t^2)^2 (-1 + t + 2 t^2 + 2 t^3)}.\\
\\

{\cal S}_{G(2,2,4)}(t)
=\frac{
-1 - t + 6 t^3 + 6 t^4 + 2 t^5}
{(-1 + t) (-1 + 2 t + 2 t^2)^2}.\\
\\

{\cal S}_{G(2,3,3)}(t)
=\frac{
(1 + t) (-1 + 2 t) (-1 - t + t^2 + 16 t^3 + 32 t^4 + 20 t^5 + 
   8 t^6)}
{(-1 + t) (-1 + 2 t + 4 t^2)^2 (-1 + t + 3 t^2 + 2 t^3)}.\\
\\

{\cal S}_{G(2,3,4)}(t)
=\frac{
\{1 + 2 t - 6 t^2 - 29 t^3 - 51 t^4 + 7 t^5 + 220 t^6 + 445 t^7 + 
 463 t^8 + 284 t^9 + 100 t^{10} + 
 16 t^{11}\}}
{(-1 + t) (-1 + t + 7 t^2 + 4 t^3)^2 (-1 + t + 4 t^2 + 
   5 t^3 + 2 t^4)}.\\
\\

{\cal S}_{G(2,3,5)}(t)
=\frac{
\{1 + 2 t - 6 t^2 - 29 t^3 - 51 t^4 + 7 t^5 + 220 t^6 + 445 t^7 + 
 463 t^8 + 284 t^9 + 100 t^{10} + 
 16 t^{11}\}}
{(-1 + t) (-1 + t + 7 t^2 + 4 t^3)^2 (-1 + t + 4 t^2 + 
   5 t^3 + 2 t^4)}.\\
\\

{\cal S}_{G(2,3,6)}(t)
=\frac{
\{(1 + t) (1 - 8 t^2 - 18 t^3 - 12 t^4 + 48 t^5 + 166 t^6 + 274 t^7 + 
   265 t^8 + 172 t^9 + 68 t^{10} + 16 t^{11})\}}
{(-1 + t) (-1 + 2 t + 
   3 t^2 + 2 t^3) (-1 + t + 7 t^2 + 7 t^3 + 4 t^4)^2}.\\
\\

{\cal S}_{G(2,3,7)}(t)
=
\frac{
\begin{array}{r}
{\scriptstyle 
\{(1 + t) (1 + 4 t - t^2 - 43 t^3 - 138 t^4 - 193 t^5 + 75 t^6 + 
     1056 t^7 + 2930 t^8 + 5284 t^9 + 7160 t^{10}}\\
{\scriptstyle + 7638 t^{11} + 
     6544 t^{12} + 4524 t^{13} + 2500 t^{14} + 1088 t^{15} + 336 t^{16} + 
     64 t^{17})\}}
\end{array}}
{(-1 + t) (-1 + 8 t^2 + 14 t^3 + 14 t^4 + 
     8 t^5)^2 (-1 + 5 t^2 + 11 t^3 + 12 t^4 + 9 t^5 + 5 t^6 + 2 t^7)}.\\
\\

{\cal S}_{G(3,3,3)}(t)
=
\frac{
(1 + t) (-1 + 2 t) (1 + 2 t^2)}
{(-1 + t) (-1 + 4 t) (-1 + 2 t + 
  2 t^2)}.\\
\\

{\cal S}_{G(3,3,4)}(t)
=
\frac{
\begin{array}{r}
{\scriptstyle \{-1 + t + 17 t^2 + 11 t^3 - 67 t^4 - 167 t^5 - 246 t^6 - 130 t^7 + 
   638 t^8 + 2152 t^9 + 3672 t^{10} + 4272 t^{11} + 3704 t^{12}}\\
 {\scriptstyle + 
   2376 t^{13} + 1040 t^{14} + 272 t^{15} + 
   32 t^{16}\}}
\end{array}}
{(-1 + t) (-1 + 2 t + 8 t^2 + 4 t^3)^2 (-1 + t + 5 t^2 + 
     6 t^3 + 2 t^4) (-1 + t + 3 t^2 + 4 t^3 + 4 t^4 + 2 t^5)}.\\
\\

{\cal S}_{G(3,6,7)}(t)
=
\frac{
\begin{array}{r}
{\scriptstyle \{-1 - 4 t + 21 t^2 + 172 t^3 + 301 t^4 - 1070 t^5 - 7231 t^6 - 
   18462 t^7 - 16780 t^8 + 58393 t^9 + 331780 t^{10}}\\
{\scriptstyle + 992604 t^{11} + 
   2257404 t^{12} + 4289278 t^{13} + 7087793 t^{14} + 10394402 t^{15}
+ 
   13680190 t^{16}}\\
{\scriptstyle + 16262820 t^{17} + 17530098 t^{18} + 17171612 t^{19}
+ 
   15299344 t^{20} + 12395068 t^{21}}\\
{\scriptstyle + 9117638 t^{22} + 6071726 t^{23}
+ 
   3643664 t^{24} + 1957008 t^{25} + 931500 t^{26}}\\
{\scriptstyle + 387360 t^{27}
+ 
   137776 t^{28} + 40544 t^{29} + 9328 t^{30} + 1504 t^{31} + 
   128 t^{32}\}}
\end{array}}
{
\begin{array}{r}
{\scriptstyle 
\{(-1 + t) (-1 + t + 8 t^2 + 15 t^3 + 15 t^4 + 11 t^5 + 
     6 t^6 + 2 t^7) (-1 + 12 t^2 + 32 t^3 + 48 t^4 + 46 t^5 + 
     26 t^6 + 8 t^7)^2}\\
{\scriptstyle \cdot (-1 + 8 t^2 + 21 t^3 + 31 t^4 + 35 t^5 + 
     32 t^6 + 24 t^7 + 13 t^8 + 6 t^9 + 2 t^{10})\}}
\end{array}
}.\\
\\

{\cal S}_{G(2,2,2,2)}(t)
=
\frac{
(1 + t) (-1 + 3 t^2)}
{(-1 + t) (-1 + 3 t)^2}.\\
\\

{\cal S}_{G(2,3,4,5)}(t)
=
\frac{
\begin{array}{r}
{\scriptstyle \{1 + 4 t - 16 t^2 - 149 t^3 - 393 t^4 + 45 t^5 + 3879 t^6 + 
   16001 t^7 + 40715 t^8}\\
{\scriptstyle  + 75854 t^9 + 109176 t^{10} + 124076 t^{11} + 
   112301 t^{12} + 80936 t^{13}}\\
{\scriptstyle  + 45998 t^{14} + 20136 t^{15} + 6464 t^{16} + 
   1368 t^{17} + 
   144 t^{18}\}}
\end{array}}
{(-1 + t) (-1 + t + 17 t^2 + 38 t^3 + 36 t^4 + 
     12 t^5)^2 (-1 + t + 12 t^2 + 25 t^3 + 30 t^4 + 22 t^5 + 11 t^6 + 
     3 t^7)}.\\
\\

\end{array}
\end{eqnarray*}
\end{ex}

\begin{eqnarray*}
\begin{array}{l}

{\cal S}_{G(2,3,7,8)}(t)
=
\frac{
\begin{array}{r}
{\scriptstyle \{1 + 4 t - 16 t^2 - 163 t^3 - 516 t^4 - 320 t^5 + 4296 t^6 + 
   24213 t^7 + 81073 t^8 + 206772 t^9 + 434218 t^{10} + 778907 t^{11}}\\
 {\scriptstyle + 
   1218441 t^{12} + 1683284 t^{13} + 2070095 t^{14} + 2277157 t^{15} + 
   2246173 t^{16}}\\
{\scriptstyle  + 1987792 t^{17} + 1576248 t^{18} + 1116383 t^{19} + 
   702382 t^{20}}\\
{\scriptstyle + 389314 t^{21} + 187754 t^{22} + 77324 t^{23} + 
   26408 t^{24}}\\
{\scriptstyle  + 7112 t^{25} + 1368 t^{26} + 
   144 t^{27}\}}
\end{array}}
{
\begin{array}{r}
{\scriptstyle \{(-1 + t) (-1 + t + 17 t^2 + 43 t^3 + 69 t^4 + 80 t^5 + 
     62 t^6 + 36 t^7 + 12 t^8)^2}\\
{\scriptstyle \cdot (-1 + t + 12 t^2 + 33 t^3 + 55 t^4 + 
     68 t^5 + 64 t^6 + 47 t^7 + 26 t^8 + 11 t^9 + 3 t^{10})\}}
\end{array}
}.\\
\\

{\cal S}_{G(2,2,2,2,2)}(t)
=
\frac{
(1 + t) (-1 + 2 t) (1 + 2 t)}
{(-1 + t) (-1 + 4 t)^2}.\\
\\

{\cal S}_{G(3,3,3,3,3)}(t)
=
\frac{
(1 + t) (-1 + 4 t - 4 t^2 + 8 t^3)}
{(-1 + t) (-1 + 8 t) (-1 + 4 t + 
   4 t^2)}.\\
\\

{\cal S}_{G(3,4,5,6,7)}(t)
=
\frac{
\begin{array}{r}
{\scriptstyle \{-1 - 4 t + 88 t^2 + 839 t^3 + 829 t^4 - 29238 t^5 - 227292 t^6 - 
   729448 t^7 + 537936 t^8 + 19529418 t^9 + 126701958 t^{10}}\\
{\scriptstyle + 
   554358269 t^{11} + 1923140887 t^{12} + 5629278332 t^{13} + 
   14371787530 t^{14}
+ 32675409006 t^{15} + 67112438033 t^{16}}\\
{\scriptstyle + 
   125827368129 t^{17} + 217031712921 t^{18} + 346438522281 t^{19}
+ 
   514115396075 t^{20} + 711769397882 t^{21}}\\
{\scriptstyle + 921744217990 t^{22} + 
   1118743167770 t^{23} + 1274436479732 t^{24}
 + 1363938860927 t^{25} + 
   1372181491114 t^{26}}\\
{\scriptstyle + 1297977449990 t^{27} + 1154299044726 t^{28} + 
   964677120102 t^{29}
 + 757061517889 t^{30} + 557294933886 t^{31}}\\
{\scriptstyle + 
   384230982310 t^{32} + 247632202450 t^{33} + 148818119668 t^{34}
   + 83136635312 t^{35} + 43006809888 t^{36}}\\
{\scriptstyle + 20501800616 t^{37} + 
   8952019048 t^{38} + 3552861968 t^{39}
+ 1268966672 t^{40} + 
   402581120 t^{41}}\\
{\scriptstyle + 111449184 t^{42} + 26255936 t^{43} + 5070080 t^{44} + 
   754944 t^{45} + 77312 t^{46} + 
   4096 t^{47}\}}
\end{array}}
{
\begin{array}{r}
{\scriptstyle \{(-1 + t) (-1 + t + 39 t^2 + 177 t^3 + 441 t^4 + 
     734 t^5 + 872 t^6 + 748 t^7 + 448 t^8 + 172 t^9 + 
     32 t^{10})^2}\\
{\scriptstyle \cdot (-1 + 2 t + 30 t^2 + 99 t^3 + 188 t^4 + 249 t^5 + 
     251 t^6 + 198 t^7 + 123 t^8 + 58 t^9 + 20 t^{10} + 4 t^{11})}\\
\begin{array}{r}
{\scriptstyle \cdot (-1 + 
     t + 24 t^2 + 90 t^3 + 209 t^4 + 364 t^5 + 510 t^6 + 594 t^7 + 
     589 t^8 + 500 t^9}\\
{\scriptstyle + 366 t^{10} + 227 t^{11} + 119 t^{12} + 51 t^{13} + 
     16 t^{14} + 4 t^{15})\}}
\end{array}\\
\end{array}
}.
\end{array}
\end{eqnarray*}


\noindent
{\Large \bf Acknowledgments}

\vspace{0.2cm}

\noindent
The first author would like to thank Yasushi Yamashita
for his useful comments concerning growth series of the fundamental groups of 
Seifert fiber spaces,
Laura Ciobanu for her interesting and instructive comments 
concerning growth series of both the Artin groups of dihedral type
and the fundamental groups of Seifert fiber spaces,
Yohei Komori and Yuriko Umemoto for their
interesting and useful talks on growth series of finitely generated groups
given in a workshop held at Osaka City University,
and Ruth Kellerhals for her useful comments and kind hospitality during his stays at University of Fribourg.



\end{document}